\documentclass[11pt]{article}

\usepackage{amssymb}
\usepackage{relsize}
\usepackage{amsthm}
\usepackage{mathtools}
\usepackage[sc]{mathpazo}
\usepackage{hyperref} 
\allowdisplaybreaks

\newcommand{\smatris}[1]{\begin{smallmatrix} #1 \end{smallmatrix}}
\newcommand{\matris}[1]{\begin{matrix} #1 \end{matrix}}
\newcommand{\wrum}{\mathcal{W}}
\newcommand{\real}{\mathbb{R}}
\newcommand{\integ}{\mathbb{Z}}
\newcommand{\ratio}{\mathbb{Q}}
\newcommand{\rmd}{\mathrm{d}}
\renewcommand{\d}{\partial}

\renewcommand{\tilde}[1]{\widetilde{#1}}
\renewcommand{\bar}[1]{\overline{#1}}
\newcommand{\vecb}[1]{\mathbf{#1}}
\renewcommand{\hat}[1]{\widehat{#1}}
\newcommand{\grad}{\nabla}
\newcommand{\eps}{\varepsilon}
\newcommand{\vecx}{\vecb{x}}
\newcommand{\vect}{\vecb{t}}
\newcommand{\vecy}{\vecb{y}}
\newcommand{\vecs}{\vecb{s}}
\newcommand{\sety}{\mathcal{Y}}
\newcommand{\ZeToT}{0,T}
\newcommand{\tint}{( \ZeToT )}
\newcommand{\smint}{\mathsmaller{\int}}
\newcommand{\ssum}{\mathsmaller{\sum\limits}}
\newcommand{\boldproof}{\textbf{\emph{Proof.}}}
\newcommand{\cont}{\mathcal{C}}
\newcommand{\inter}[1]{\shortintertext{#1}}
\newcommand{\scaleconv}[1]{\overset{{}_{#1}}{-\!\!-\!\!\!\rightharpoonup}}
\newcommand{\distr}{\mathcal{D}}
\newcommand{\aoper}{\mathcal{A}}
\newcommand{\scope}[1]{[ \mspace{-3.5mu} [ #1 ] \mspace{-3.5mu} ]}
\newcommand{\preq}{\preccurlyeq}
\newcommand{\sep}{\mathrm{sep}}
\newcommand{\wsep}{\mathrm{wsep}}
\newcommand{\jointly}{\mathcal{J}}
\newcommand{\ellring}{\ell \, \mathring{ }}
\newcommand{\hatellring}{\hat{\ell} \, \mathring{ }}
\newcommand{\iring}{i \, \mathring{ }}

\theoremstyle{theorem} \newtheorem{defin}{Definition}
\theoremstyle{theorem} \newtheorem{theo}[defin]{Theorem}
\theoremstyle{theorem} \newtheorem{prop}[defin]{Proposition}
\theoremstyle{theorem} \newtheorem{coro}[defin]{Corollary}
\theoremstyle{theorem} \newtheorem{rema}[defin]{Remark}
\theoremstyle{theorem} \newtheorem{lemma}[defin]{Lemma}
\theoremstyle{theorem} \newtheorem{exam}[defin]{Example}

\begin{document}

\title{\textbf{Homogenisation of Monotone Parabolic} \\ \textbf{Problems with Several Temporal Scales} \\
		----------------- \\ \Large{\textbf{The Detailed arXiv e-Print Version}}}

\author{Jens Persson\footnote{E-mail: \texttt{jens.persson@miun.se}} \\ \vspace{0.25cm} \\ \emph{Department
		of Engineering and Sustainable Development,} \\ \emph{Mid Sweden University, SE-831~25 \"{O}stersund, Sweden}}
		
\date{}

\maketitle

\begin{abstract}
In this paper we homogenise monotone parabolic problems with two spatial scales and finitely many temporal scales.
Under a certain well-separatedness assumption on the spatial and temporal scales as explained in the paper,
we show that there is an H-limit defined by at most four distinct sets of local problems corresponding to
slow temporal oscillations, slow resonant spatial and temporal oscillations (the ``slow'' self-similar case), rapid temporal
oscillations, and rapid resonant spatial and temporal oscillations (the ``rapid'' self-similar case), respectively.
\end{abstract}
					${} \qquad \qquad \quad \; \;$ monotone, several temporal scales 

\section{Introduction}

We will give here a brief survey---with some important references---of homogenisation theory and two-scale convergence
techniques which is followed by a statement of the research objective of the present paper. Finally in this section we
give a list of notations employed in the paper.

\textbf{Homogenisation theory.} Homogenisation theory is the study of the convergence of---in some suitable
sense---sequences of equations involving sequences of operators and (possibly) source functions and the responding
sequences of solutions. The main applications involve the study of the convergence of sequences of partial differential
equations described by heterogeneous coefficients which become more and more refined such that the problem tends to a
homogenised limit. In the case of parabolic partial differential equations the convergence modes used to achieve homogenised
limits are the so called G- and H-convergences, where the former is employed when the coefficients can be arranged as a
symmetric matrix (see \cite{Spag67,Spag68}), and the latter is the generalisation which includes non-symmetric matrices
(see \cite{Mura78,Mura78b,Tart78,Tart79}) and even non-linear problems (see \cite{Tart77}). ``Homogenising'' a problem means
in this context to find the limit in the G- or H-convergence process.

\textbf{Two-scale convergence.} The theory of homogenisation experienced a quantum leap in the late 1980's when the two-scale
convergence technique was introduced (see \cite{Ngue89,Alla92})---effectively replacing Tartar's method of oscillating test
functions (see \cite{Tart77,Tart78}) as the main tool to achieve G- or H-convergence---and the technique has subsequently
improved since then. Two-scale convergence (with generalisations such as multiscale convergence \cite{AllBri96},
``generalised'' two-scale convergence \cite{Holm96,HSSW06}, scale convergence \cite{MasToa01}, $\lambda$-scale
convergence \cite{HolSil06,Silf07}, $\Sigma$-convergence \cite{Ngue03,Ngue04} etc.) is today an indispensable tool
to the modern homogenisation theorist.

\textbf{Aims in the present paper.} The main purpose of this paper is to perform homogenisation of monotone, possibly non-linear,
parabolic problems of the type
\begin{equation}
	\label{eq:monoparaprob}
	\left\{
		\begin{aligned}
			\tfrac{\d}{\d t} u_{\eps} (x,t) - \grad \cdot a ( x,t \, , \,  \tfrac{x}{\eps} \, , \, \tfrac{t}{\eps_1'},
				\ldots, \tfrac{t}{\eps_m'} \, ; \, \grad u_{\eps})
							&= f(x,t)	&	&\textrm{in } \Omega \times (0,T), \\
			u_{\eps}(x,0)	&= u_0(x)	&	&\textrm{in } \Omega, \\
			u_{\eps}(x,t)	&= 0		&	&\textrm{on }\d \Omega \times (0,T),
		\end{aligned}
	\right.
\end{equation}
i.e., having two spatial and $m+1$ temporal scales, where $\Omega$ is an open bounded set in $\real^N$ and $T > 0$.
As $\eps$ tends to $0$ we get a sequence of equations given by \eqref{eq:monoparaprob} above and  the objective is
to find the homogenised problem, i.e., to find the homogenised limit $b$ of the flux $a$ which defines a homogenised
equation which admits a limit $u$ of the sequence of solutions $\{ u_{\eps} \}$. In order to homogenise
\eqref{eq:monoparaprob} we impose a certain separatedness restriction on the scale functions $\eps, \eps_1',\ldots,\eps_m'$.
The homogenised limit $b$ will not contain any fast spatial or temporal oscillations and (if considered as a function
of $\grad u$) is given in terms of an integral over the local variables $y,s_1,\ldots,s_m$ involving the flux $a$ and
a function $u_1$ which is the unique solution of some local problems depending on the behaviour of the scale functions.
We discern four distinct cases giving different local problems for $u_1$, namely
the cases (i)~$\eps^2 / \eps_m' \rightarrow 0$ as $\eps \rightarrow 0$, (ii)~$\eps_m'
\sim \eps^2$, and (iii)~$\eps_i' / \eps^2 \rightarrow 0$ but $\eps_{i-1}' /
\eps^2 \rightarrow \infty$ as $\eps \rightarrow 0$ for some $\eps_i'$ tending more rapidly to $0$
than $\eps$ does, and (iv)~$\eps_{\ellring-1}' \sim \eps^2$ for some $\eps_{\ellring-1}'
\neq \eps_m'$ tending more rapidly to $0$ than $\eps$ does. Case (i) corresponds to slow temporal
oscillations (compared to the spatial one), (ii) is the so-called ``slow'' self-similar case where the spatial and temporal
oscillations are in resonance, (iii) corresponds to rapid temporal oscillations, and (iv) is the ``rapid'' self-similar case.

\textbf{Notations and conventions.} The following notations and conventions are used in this paper:

\emph{Spatial and temporal domains.} Throughout the paper, $\Omega$ defining the spatial domain is a non-empty open bounded
set in $\real^N$ with Lipschitz boundary, and $T > 0$ is the maximal time defining the temporal domain $(0,T)$.

\emph{Sets of positive integers.} We define the following convenient subsets of $\integ$: for any $0 < i \leqslant j$ in $\real$,
$\scope{i,j} = [i,j] \cap \integ$ (the integers between $i$ and $j$); in particular, $\scope{j} = \scope{1,j}$ (the positive
integers up to at most $j$). Moreover, if $i < j$, we define $\scope{0},\scope{j,i} = \emptyset$ (empty sets of positive
integers); note that we employ the convention that statements over the empty set are by default always trivially true.
Examples: $\scope{2,4} = \{ 2,3,4 \}$, $\scope{\tfrac{10}{3}} = \{ 1,2,3 \}$, $\scope{4,2} = \emptyset$,
$\scope{\tfrac{3}{10}} = \scope{1,\tfrac{3}{10}} = \emptyset$, and $\ell > \ell$ for all $\ell \in \emptyset$.

\emph{Functions with mean value zero and periodic functions.} Let $\mathcal{F}(A)/\mathbb{R}$ denote all functions in
$\mathcal{F}(A)$ with mean value zero over $A \subset \mathbb{R}^M$, and let $\mathcal{F}_{\#}(Z)$ denote all locally
$\mathcal{F}$ functions over $\mathbb{R}^M$ that are periodical repetitions of some functions in $\mathcal{F}(Z)$ where
$Z = (0,1)^M$. In particular, $\mathcal{F}_{\#}(Z)/\mathbb{R}$ is the set of locally $\mathcal{F}$ functions over
$\mathbb{R}^M$ with mean value zero over $Z$ which are periodic repetitions of some functions in $\mathcal{F}(Z)$.

\emph{Tensor product sets.} The subset $\mathcal{F}_1 (A_1) \odot \cdots \odot \mathcal{F}_k (A_k)$ of the tensor product
$\mathcal{F}_1 (A_1) \otimes \cdots \otimes \mathcal{F}_k (A_k)$ of function spaces $\mathcal{F}_1 (A_1), \ldots,
\mathcal{F}_k (A_k)$ is the set of all functions $f$ that can be written as the tensor product
\begin{align*}
	f					&= f_1 \otimes \cdots \otimes f_k,
\inter{i.e.,}
	f(z_1,\ldots,z_k)	&= f_1(z_1) \cdots f_k(z_k) \qquad \qquad (z_i \in A_i, \, i \in \scope{k}),
\end{align*}
for some $f_i \in \mathcal{F}_i$, $i \in \scope{k}$. We say that $\mathcal{F}_1 (A_1) \odot \cdots \odot \mathcal{F}_k (A_k)$
is a tensor product set (which we note spans the tensor product space). Example:  Any function $\psi$ in the tensor product
set $\distr (\Omega) \, \odot \, \distr (0,T)$ defined on $\Omega \times (0,T)$ can be written as
\begin{align*}
\psi			&= v \otimes c,
\inter{i.e.,}
	\psi(x,t)	&= v(x) \, c(t) \qquad \qquad (x \in \Omega, \, t \in (0,T)),
\end{align*}
for some $v \in \distr (\Omega)$, $c \in \distr (0,T)$. (Note that $\distr = \cont_0^{\infty}$, i.e., infinitely
differentiable functions with compact support on the set argument.) 
			
\emph{Placement of $\eps$-indices.} When $\eps$ is an upper index it refers to an explicit construction like, e.g.,
\begin{equation}
	\label{eq:upperindex}
	\psi^{\eps} (x,t) = \psi(x,t \, , \,  \tfrac{x}{\eps_1}, \ldots, \tfrac{x}{\eps_n}
		\, , \, \tfrac{t}{\eps_1'}, \ldots, \tfrac{t}{\eps_m'}),
\end{equation}
for functions $\psi$ defined on, in this case, $\Omega \times (0,T) \, \times \,(0,1)^{nN} \times (0,1)^m$.
A lower index form refers to an implicit construction not based on \eqref{eq:upperindex}; see, e.g.,
the solution $u_{\eps}$ to \eqref{eq:monoparaprob} where $\eps$ only indirectly defines the function.

\emph{Partial derivatives.} There are two kinds of partial derivatives. The partial derivatives of the first kind,
$\grad = \bigl( \tfrac{\d}{\d x_1}, \ldots, \tfrac{\d}{\d x_N} \bigr)$ and
$\tfrac{\d}{\d t}$, only discern whether one differentiates with respect to the space variable
$x = (x_1,\ldots,x_N)$ or the time variable $t$, respectively. The partial derivatives of the second kind,
$\grad_x = \bigl( \d_{x_1}, \ldots, \d_{x_N} \bigr)$ and $\d_t$ (i.e., with the variable
as a subscript) are proper partial derivatives with respect to space and time, respectively. Note that partial
derivatives of the local variables will always be of the proper, second kind. Example: Let $\psi = \psi(x,t,y,s)$
be a weakly differentiable real-valued function with respect to the global space and time variables $x$ and $t$
and the local space and time variables $y$ and $s$. Suppose $y = \eta x$ and $s = \sigma t$ for some real constants
$\eta$ and $\sigma$, then the chain rule and the conventions above give
\begin{equation*}
	\grad \psi	= \grad_x \psi + \eta \grad_y \psi \quad \textrm{and}
		\quad \tfrac{\d}{\d t} \psi = \d_t \psi + \sigma \d_s \psi;
\end{equation*}
these differentiation rules will be important to keep in mind later in this paper.

\emph{Hilbert space structure.} We use the convention that we work solely in $L^2$ and derivations such as $H^1$,
$H^1_{\#}/\real$ etc.~rather than in the more general $L^p$, $p \geqslant 1$, with derivations $W^{1,p}$,
$W_{\#}^{1,p}/\real$ etc. The reason we work in $L^2$ is of course due to the fact that it is a Hilbert space which
means that the topological dual is also $L^2$. Heuristically speaking, Hilbert spaces such as e.g.~$L^2$ are more
``natural'' than non-Hilbert spaces since they are generalisations of finite-dimensional vector spaces. The drawback
is that we lose some important examples of non-linear problems such as e.g.~the evolution $p$-Laplacian equation
(with $p \neq 2$) which describes non-linear diffusion phenomena and employed in e.g. image processing \cite{Weic98}.

\section{Multiscale Convergence}

The concept of two-scale convergence was introduced in 1989 by Nguetseng (see \cite{Ngue89}) and further developed
by Allaire in 1992 (see \cite{Alla92}). In words, two-scale convergence is a kind of weak convergence mode for a
sequence of functions of a global variable where the limit is a function of both the global (or macroscopic) and
the local (or microscopic) variable. (For an excellent review on two-scale convergence, see \cite{LuNgWa02}.) By
using the periodic unfolding (or two-scale transform) technique (see \cite{ArDoHo90,CiDaGr02}) or alternatively
the inverse two-scale transform technique (see \cite{Nech04}), this peculiar convergence mode is realised to be
equivalent to an ordinary weak convergence for sequences of functions which depends both on the global and the
local variable.

The rigorous definition of two-scale convergence is given below. (If nothing else is stated, in this paper we let
$y \in Y$ where $Y = (0,1)^N$).
\begin{defin}
A sequence $\{ u_ {\eps} \}$ of functions in $L^2(\Omega)$ is said to two-scale converge to a limit $u_0 \in L^2(\Omega \times Y)$
if, as $\eps \rightarrow 0$ (from above),
\begin{equation}
	\label{eq:two-scaleconvdef}
	\smint_{\Omega} u_{\eps} (x) \, v(x,\tfrac{x}{\eps}) \, \rmd x
		\rightarrow \smint_{\Omega} \smint_{Y} u_0(x,y) \, v(x,y) \, \rmd y \, \rmd x
\end{equation}
for all $v \in L^2 \bigl( \Omega; \cont_{\#}(Y) \bigr)$, and we write $u_{\eps} \overset{2}{\rightharpoonup} u_0$
as $\eps \rightarrow 0$.
\end{defin}
\begin{rema}
Alternatively one can write ``$\rightharpoonup \rightharpoonup$'' instead of ``$\overset{2}{\rightharpoonup}$''.
Note also that instead of using the positive scale parameter $\eps$ tending to zero it is possible to employ a
perhaps more fundamental scale parameter $h$ tending to positive infinity. (This means that
$\lim\limits_{h \rightarrow \infty} \eps = 0$; in the remainder of the paper this can at any point be achieved by
substituting $\eps = 1/h$. The substitution would, e.g., give $hx$ instead of $\tfrac{x}{\eps}$ everywhere.)
\end{rema}
From now on we assume that all limits are taken as $\eps \rightarrow 0$ (from above) if nothing else is stated.

In Definition~\ref{def:scalefcn} below we introduce the notion of scale functions which are functions with respect
to the scale parameter.
\begin{defin}
\label{def:scalefcn}
A scale function $\eps_{\ast} \, : \, \real_+ \rightarrow \real$ is a real-valued function of the scale parameter
$\eps$ for which $\eps_{\ast} (\eps) \rightarrow 0$ (i.e., $\eps_{\ast}$ is microscopic), and for which there
exists $\delta > 0$ such that $\eps_{\ast}(\eps) > 0$ for all $0 < \eps < \delta$ (i.e., $\eps_{\ast}$ is
ultimately positive).
\end{defin}
Note that the scale parameter $\eps$ itself (i.e., $\eps_{\ast}(\eps) = \eps$) is a trivial example of a
scale function. An example of a function $\eps_{\ast}$ of $\eps$ that is not a scale function is, e.g.,
$\eps_{\ast}(\eps) = \eps \, \sin \tfrac{1}{\eps}$ since $\eps_{\ast}$ in this case---though being
microscopic---is not ultimately positive.

The concept of scale functions leads to the notion of multiscale convergence which was introduced in 1996
by Allaire and Briane (see \cite{AllBri96}) as a generalisation of two-scale convergence in order to be
able to perform homogenisation of problems with multiple scales. This convergence mode is defined below.
(If nothing else is stated, in this paper we let $y_i \in Y_i$, where $Y_i = (0,1)^N$, $i \in \scope{n}$.)
\begin{defin}
A sequence $\{ u_ {\eps} \}$ of functions in $L^2(\Omega)$ is said to $(n+1)$-scale converge to a limit
$u_0 \in L^2(\Omega \times Y_1 \times \cdots \times Y_n)$ if
\begin{align}
	\label{eq:multiscaleconvdef}
	&\smint_{\Omega} u_{\eps} (x) \, v(x,\tfrac{x}{\eps_1},\ldots,\tfrac{x}{\eps_n}) \, \rmd x \nonumber \\
	&\phantom{ \smint_{\Omega} u_{\eps} (x) }
		\rightarrow \smint_{\Omega} \smint_{Y_1} \cdots \smint_{Y_n} u_0(x,y_1,\ldots,y_n)
			\, v(x,y_1,\ldots,y_n) \, \rmd y_n \cdots \rmd y_1 \, \rmd x
\end{align}
for all $v \in L^2 \bigl( \Omega; \cont_{\#}(Y_1 \times \cdots \times Y_n) \bigr)$, and we write
$u_{\eps} \scaleconv{n+1} u_0$.
\end{defin}
In order to simplify the notation, from now on we will write $\vecy_n = (y_1,\ldots,y_n)$ and
$Y^n = Y_1 \times \cdots \times Y_n$ so that $\vecy_n \in Y^n$ which collects the local (spatial)
variables and local (spatial) sets under one roof. (Naturally, the Lebesgue measure on $Y^n$ is denoted
$\rmd \vecy_n$.) We also write $\vecx_n^{\eps} = (\tfrac{x}{\eps_1},\ldots,\tfrac{x}{\eps_n})$ in the same
spirit where we note that $\vecx_n^{\eps}$ actually depends on the particular choice of scale functions
$\eps_1,\ldots,\eps_n$. Of course, multiscale convergence is highly dependent on the behaviour of the (spatial)
scale functions. For ordered lists of scale functions we have the following definitions:
\begin{defin}
The list $\{ \eps_i \}_{i=1}^n$ of scale functions is said to be separated if $\tfrac{\eps_{k+1}}{\eps_k}
\rightarrow 0$ for all $k \in \scope{n-1}$.
\end{defin}
\begin{defin}
The list $\{ \eps_i \}_{i=1}^n$ of scale functions is said to be well-separated if there exists a positive
integer $\ell$ such that $\tfrac{1}{\eps_k} \bigl( \tfrac{\eps_{k+1}}{\eps_k} \bigr)^{\ell} \rightarrow 0$ for
all $k \in \scope{n-1}$.
\end{defin}
\begin{rema}
Note that well-separatedness is a stronger requirement than separatedness.
\end{rema}
Homogenisation for linear parabolic problems with several temporal scales using the multiscale convergence
technique was first achieved by Flod\'{e}n and Olsson in 2007 (see \cite{FloOls07}). This was a further development
of the work by Holmbom in 1996 and 1997 (see \cite{Holm96} and \cite{Holm97}, respectively) where two-scale
convergence was employed to homogenise linear parabolic problems with both a spatial and a temporal microscale.
General $(n+1,m+1)$-scale convergence can be expressed according to the definition below. (If nothing else is stated,
in this paper we let $s_j \in S_j$, where $S_j = (0,1)$, $j \in \scope{m}$.)
\begin{defin}
A sequence $\{ u_ {\eps} \}$ in $L^2(\Omega \times \tint)$ is said to $(n+1,m+1)$-scale converge to a limit
$u_0 \in L^2 \bigl( \Omega \times \tint \, \times \, Y^n \times S_1 \times \cdots \times S_m \bigr)$ if
\begin{align}
	\label{eq:nm-scaleconvdef}
	&\smint_0^T \smint_{\Omega} u_{\eps} (x,t) \, v(x,t,\vecx_n^{\eps},\tfrac{t}{\eps_1'},
		\ldots,\tfrac{t}{\eps_m'}) \, \rmd x \rmd t \nonumber \\
	&\phantom{ \smint_0^T \smint_{\Omega} u_{\eps} }
		\rightarrow \smint_0^T \smint_{\Omega} \smint_{Y^n} \smint_{S_1} \cdots \smint_{S_m}
			u_0(x,t,\vecy_n,s_1,\ldots,s_m) \nonumber \\
	&\phantom{ \smint_0^T \smint_{\Omega} u_{\eps}
				\rightarrow \smint_{\tint} \smint_{\Omega} \smint_{Y^n} \smint_{S_1} \cdots \smint_{S_m}}
		\times v(x,t,\vecy_n,s_1,\ldots,s_m) \, \rmd s_m \cdots \rmd s_1 \, \rmd \vecy_n \rmd x \rmd t
\end{align}
for all $v \in L^2 \bigl( \Omega \times \tint; \cont_{\#}(Y^n \times S_1 \times \cdots \times S_m) \bigr)$,
and we write $u_{\eps} \scaleconv{(n+1,m+1)} u_0$.
\end{defin}
Trivially, this definition also works for vector valued functions where the product becomes a dot product, or mixed
scalar and vector valued functions which would give vector valued integrals above. All results below concerning the notion of
$(n+1,m+1)$-scale convergence can of course be formulated for such functions as well. In particular, gradient
functions will later be of interest.

In order to simplify the notation, from now on we will write $\vecs_m = (s_1,\ldots,s_m)$ and $S^m = S_1
\times \cdots \times S_m$ so that $\vecs_m \in S^m$. (The Lebesgue measure on $S^m$ will of course be
denoted $\rmd \vecs_m$.) Moreover, $\vect_m^{\eps} = (\tfrac{t}{\eps_1'},\ldots,\tfrac{t}{\eps_m'})$ which is
noted to depend on the particular choice of temporal scale functions $\{ \eps_j' \}_{j=1}^m$. Furthermore, introduce
$\Omega_T = \Omega \times \tint$ so that $(x,t) \in \Omega_T$, and $\sety_{nm} = Y^n \times S^m$ so that $(\vecy_n,\vecs_m)
\in \sety_{nm}$.

It is clear that we need to introduce some convenient restrictions on the spatial and temporal scale functions
$\{ \eps_i \}_{i=1}^n$ and $\{ \eps_j' \}_{j=1}^m$ in order for them to collaborate in a meritorious manner.
In Definition~\ref{def:jointsep} below we define a certain set of pairs of lists of such meritoriously collaborating
spatial and temporal scale functions.
\begin{defin}
\label{def:jointsep}
Suppose we have a list $\{ \eps_i \}_{i=1}^n$ of $n$ spatial scale functions
and a list $\{ \eps_j' \}_{j=1}^m$ of $m$ temporal scale functions. We say that the pair $\bigl( \{ \eps_i \}_{i=1}^n ,
\{ \eps_j' \}_{j=1}^m \bigr)$ belongs to the set $\jointly_{\sep}^{nm}$ if  $\{ \eps_i \}_{i=1}^n$ and
$\{ \eps_j' \}_{j=1}^m$ are both separated and that the following two conditions hold:
\begin{itemize}
	\item[(i)] There exist possibly empty subsets $A \subset \scope{n}$ and $A' \subset \scope{m}$ with $|A| = |A'| = k$ such that
				there exist bijections $\beta \, : \, \scope{k} \rightarrow A$ and $\beta' \, : \, \scope{k} \rightarrow A'$,
				respectively, such that $\eps_{\beta(i)} = \eps_{\beta'(i)}'$ for all $i \in \scope{k}$. (We have no requirement
				in the empty case $k = 0$.)
	\item[(ii)] There exists a permutation $\pi$ on the set $\scope{n+m-2k}$ such that the permutation
				$\{ \eps_{\pi(\ell)}'' \}_{\ell = 1}^{n+m-2k}$ of the list
				\begin{equation*}
					\{ \eps_{\ell}'' \}_{\ell=1}^{n+m-2k} = \bigl\{ \{ \eps_i \}_{i \not\in A} , \{ \eps_j' \}_{j \not\in A'} \bigr\}
				\end{equation*}
				of the remaining $n+m-2k$ scale functions is separated. (We have no requirement in the empty case $n + m -2k = 0$.)
\end{itemize}
If we require well-separatedness instead of mere separatedness we can define the corresponding set $\jointly_{\wsep}^{nm}$.
\end{defin}
Note that $\jointly_{\wsep}^{nm} \subset \jointly_{\sep}^{nm}$. The idea of the definition above is that we can
localise all the spatial and temporal scale functions in two disjoint categories, (i) and (ii), where the former category consists
of those that are equal and the latter category consists of those that are jointly (well-)separated. Note also that since neither
$n$ nor $m$ vanishes, it can not be the case that both categories (i) and (ii) of Definition~\ref{def:jointsep} are empty.
\begin{exam}
As examples of pairs of lists that do and do not belong to $\jointly_{\wsep}^{nm}$, consider $(e_1,e_1')$, $(e_2,e_2')$
and $(e_3,e_3')$ defined by
\begin{equation*}
	\left\{
		\begin{aligned}
			e_1		&= \{ \eps, \eps^3 \},		&\qquad	e_1'	&= \{ \eps^2, \eps^3, \eps^4 \}, \\
			e_2		&= \{ \eps, \eps^3 \},		&\qquad	e_2'	&= \{ \eps^2, \tfrac{\eps^2}{|\log \eps|}, \eps^3 \}, \\
			e_3		&= \{ \eps, \eps^3 \},		&\qquad	e_3'	&= \{ \eps, \eps^2, \tfrac{\eps^3}{|\log \eps|} \}.
		\end{aligned}
	\right.
\end{equation*}

Clearly, the first pair $(e_1,e_1')$ belongs to $\jointly_{\wsep}^{2 \, 3}$ since both $e_1$ and $e_1'$ are well-separated
lists and the combined list $\{ \eps, \eps^2, \eps^4 \}$ where we have removed the common scale function $\eps^3$ is well-separated.

We have that the middle pair $(e_2,e_2')$ does not belong to $\jointly_{\wsep}^{2 \, 3}$ since $e_2'$ is not well-separated.

The last pair $(e_3,e_3')$ does not belong to $\jointly_{\wsep}^{2 \, 3}$. Indeed, the combined list
$\{ \eps^2, \eps^3, \tfrac{\eps^3}{|\log \eps|} \}$ (with removed common scale function $\eps$) is not well-separated.
\end{exam}
In Proposition~\ref{prop:asymptnmscaleconv} below we recall that if $q$, $f$ and $g$ are functions of $\eps$ where $f = qg$
and $q \rightarrow 1$, then we say that $f \sim g$, i.e., $f$ and $g$ are asymptotically equal.
\begin{prop}
\label{prop:asymptnmscaleconv}
Suppose $u_{\eps} \scaleconv{(n+1,m+1)} u_0$ and that $r = r(\eps)$ satisfies $r \sim r_0$, $r_0 \in \real$.
Then $r (\eps) u_{\eps} \scaleconv{(n+1,m+1)} r_0 u_0$.
\end{prop}
\begin{proof}[\boldproof]
Clearly,
\begin{align*}
	&\smint_{\Omega_T} \bigl( r(\eps) u_{\eps}(x,t) \bigr) \, v(x,t,\vecx_n^{\eps},\vect_m^{\eps}) \, \rmd x \rmd t \\
	&\phantom{ \smint_{\Omega_T} \bigl( r(\eps) u_{\eps}(x,t) \bigr) \, }
		= r(\eps) \smint_{\Omega_T} u_{\eps}(x,t) \, v(x,t,\vecx_n^{\eps},\vect_m^{\eps}) \, \rmd x \rmd t \\
	&\phantom{ \smint_{\Omega_T} \bigl( r(\eps) u_{\eps}(x,t) \bigr) \, }
		\rightarrow r_0 \smint_{\Omega_T} \smint_{\sety_{nm}} u_0 (x,t,\vecy_n,\vecs_m)
			\, v(x,t,\vecy_n,\vecs_m) \, \rmd \vecs_m \rmd \vecy_n \rmd x \rmd t \\
	&\phantom{ \smint_{\Omega_T} \bigl( r(\eps) u_{\eps}(x,t) \bigr) \, }
		= \smint_{\Omega_T} \smint_{\sety_{nm}} \bigl( r_0 u_0 (x,t,\vecy_n,\vecs_m) \bigr)
			\, v(x,t,\vecy_n,\vecs_m) \, \rmd \vecs_m \rmd \vecy_n \rmd x \rmd t
\end{align*}
for all $v \in L^2 \bigl( \Omega_T ; \cont_{\#} (\sety_{nm}) \bigr)$, which precisely means that $r (\eps) u_{\eps}
\scaleconv{(n+1,m+1)} r_0 u_0$.
\end{proof}
Under certain restrictions it can be shown that \eqref{eq:nm-scaleconvdef} only has to hold for a certain class of
smooth functions in order to get $(n+1,m+1)$-scale convergence; see the proposition below.
\begin{prop}
\label{prop:convinsmallspace}
Let $\{ u_{\eps} \}$ be a bounded sequence in $L^2 (\Omega_T)$ and let $u_0 \in L^2(\Omega_T \times \sety_{nm})$.
Furthermore, suppose \eqref{eq:nm-scaleconvdef} holds for all
$v \in \distr \bigl( \Omega_T; \cont_{\#}^{\infty} (\sety_{nm}) \bigr)$. Then $u_{\eps} \scaleconv{(n+1,m+1)} u_0$.
\end{prop}
\begin{proof}[\boldproof]
Let $w \in L^2 \bigl( \Omega_T; \cont_{\#} (\sety_{nm}) \bigr)$ be arbitrary. Furthermore, let $\{ v_{\mu} \}$
be a sequence in $\distr \bigl( \Omega_T; \cont_{\#}^{\infty} (\sety_{nm}) \bigr)$ that converges to $w$ in
$L^2 \bigl( \Omega_T; \cont_{\#} (\sety_{nm}) \bigr)$ as $\mu \rightarrow \infty$. It is trivial that
\begin{multline}
	\label{eq:trivialconv}
	\lim_{\eps \rightarrow 0} \smint_{\Omega_T}	u_{\eps} (x,t) \, w (x,t,\vecx_n^{\eps},\vect_m^{\eps}) \, \rmd x \rmd t \\
	= \lim_{\mu \rightarrow \infty} \lim_{\eps \rightarrow 0}
			\Bigl( \smint_{\Omega_T} u_{\eps} (x,t) \, (w-v_{\mu}) (x,t,\vecx_n^{\eps},\vect_m^{\eps}) \, \rmd x \rmd t \\
				+ \smint_{\Omega_T} u_{\eps} (x,t) \, v_{\mu} (x,t,\vecx_n^{\eps},\vect_m^{\eps}) \, \rmd x \rmd t \Bigr)
\end{multline}
holds.

By assumption, for the second term in the right-hand side of \eqref{eq:trivialconv} we have
\begin{align*}
	&\lim_{\mu \rightarrow \infty} \lim_{\eps \rightarrow 0} \smint_{\Omega_T} u_{\eps} (x,t)
		\, v_{\mu} (x,t,\vecx_n^{\eps},\vect_m^{\eps}) \, \rmd x \rmd t \\
	&\phantom{ \lim_{\mu \rightarrow \infty} \lim_{\eps \rightarrow 0} \smint_{\Omega_T} u_{\eps} (x,t) }
		= \lim_{\mu \rightarrow \infty} \smint_{\Omega_T} \smint_{\sety_{nm}} u_0(x,t,\vecy_n,\vecs_m) \,
			v_{\mu} (x,t,\vecy_n,\vecs_m) \, \rmd \vecs_m \rmd \vecy_n \, \rmd x \rmd t \\
	&\phantom{ \lim_{\mu \rightarrow \infty} \lim_{\eps \rightarrow 0} \smint_{\Omega_T} u_{\eps} (x,t) }
		= \smint_{\Omega_T} \smint_{\sety_{nm}} u_0(x,t,\vecy_n,\vecs_m) \,
			w (x,t,\vecy_n,\vecs_m) \, \rmd \vecs_m \rmd \vecy_n \, \rmd x \rmd t.
\end{align*}
The second equality comes from the fact that
\begin{align*}
	&\Bigl| \smint_{\Omega_T} \smint_{\sety_{nm}} u_0 (x,t,\vecy_n,\vecs_m) \,
				(v_{\mu} - w) (x,t,\vecy_n,\vecs_m) \, \rmd \vecs_m \rmd \vecy_n \, \rmd x \rmd t \Bigr| \\
	&\phantom{ \Bigl| \smint_{\Omega_T} \smint_{\sety_{nm}} u_0 (x,t,\vecy_n,\vecs_m) \, }
		\leqslant \bigl\| u_0 \, (v_{\mu} - w) \bigr\|_{L^1 (\Omega_T \times \sety_{nm})}
			\leqslant C_1 \bigl\| v_{\mu} - w \bigr\|_{L^2 (\Omega_T \times \sety_{nm})} \\
	&\phantom{ \Bigl| \smint_{\Omega_T} \smint_{\sety_{nm}} u_0 (x,t,\vecy_n,\vecs_m) \, }
		\leqslant C_1 \bigl\| v_{\mu} - w \bigr\|_{L^2 \bigl( \Omega_T; \cont_{\#} (\sety_{nm}) \bigr)}
			\rightarrow 0
\end{align*}
as $\mu \rightarrow \infty$, where we have used H\"{o}lder's inequality in the second inequality.

It remains to treat the first term in the right-hand side of \eqref{eq:trivialconv}; we want it to vanish. Indeed,
\begin{align*}
	&\lim_{\mu \rightarrow \infty} \lim_{\eps \rightarrow 0} \smint_{\Omega_T} u_{\eps} (x,t)
			\, (w-v_{\mu}) (x,t,\vecx_n^{\eps},\vect_m^{\eps}) \, \rmd x \rmd t \\
	&\phantom{ \lim_{\mu \rightarrow \infty} \lim_{\eps \rightarrow 0} \smint_{\Omega_T} u_{\eps} (x,t) }
		\leqslant \lim_{\mu \rightarrow \infty} \lim_{\eps \rightarrow 0} \bigl\| u_{\eps}
			\, (w^{\eps}-v_{\mu}^{\eps}) \bigr\|_{L^1 (\Omega_T)}
				\leqslant \lim_{\mu \rightarrow \infty} \lim_{\eps \rightarrow 0}
					C_2 \bigl\| w^{\eps} - v_{\mu}^{\eps} \bigr\|_{L^2 (\Omega_T)} \\
	&\phantom{ \lim_{\mu \rightarrow \infty} \lim_{\eps \rightarrow 0} \smint_{\Omega_T} u_{\eps} (x,t) }
		\leqslant \lim_{\mu \rightarrow \infty} \lim_{\eps \rightarrow 0}
			C_2 \bigl\| w - v_{\mu} \bigr\|_{L^2 \bigl( \Omega_T; \cont_{\#} (\sety_{nm}) \bigr)} = 0,
\end{align*}
where we have used H\"{o}lder's inequality in the second inequality and employed that $\{ u_{\eps} \}$ is bounded
in $L^2 (\Omega_T)$. (The last inequality follows from the fact that the
$L^2 \bigl( \Omega_T; \cont_{\#} (\sety_{nm}) \bigr)$-norm involves a maximum with respect to the local variables.)

To conclude, \eqref{eq:trivialconv} becomes
\begin{align*}
	&\lim_{\eps \rightarrow 0} \smint_{\Omega_T} u_{\eps} (x,t) \, w (x,t,\vecx_n^{\eps},\vect_m^{\eps}) \, \rmd x \rmd t \\
	&\phantom{ \lim_{\eps \rightarrow 0} \smint_{\Omega_T} u_{\eps} (x,t) }
		= \smint_{\Omega_T} \smint_{\sety_{nm}} u_0(x,t,\vecy_n,\vecs_m) \,
			w (x,t,\vecy_n,\vecs_m) \, \rmd \vecs_m \rmd \vecy_n \, \rmd x \rmd t
\end{align*}
for all $w \in L^2 \bigl( \Omega_T; \cont_{\#} (\sety_{nm}) \bigr)$; we have in fact shown that
$u_{\eps} \scaleconv{(n+1,m+1)} u_0$.
\end{proof}
We have the following important compactness result.
\begin{theo}
\label{th:nm-scaleconvres}
Suppose that the pair $\bigl( \{ \eps_i \}_{i=1}^n , \{ \eps_j' \}_{j=1}^m \bigr)$ of lists of spatial and
temporal scale functions belongs to $\jointly_{\sep}^{nm}$. Furthermore, let $\{ u_{\eps} \}$ be
a bounded sequence in $L^2(\Omega_T)$. Then there is a function $u_0 \in L^2(\Omega_T \times \sety_{nm})$
such that, up to a subsequence, $u_{\eps} \scaleconv{(n+1,m+1)} u_0$.
\end{theo}
\begin{proof}[\boldproof]
(We assume here that both categories (i) and (ii) of Definition~\ref{def:jointsep} are non-empty, i.e.,
that $k \in \scope{ \bigl\lfloor \tfrac{1}{2}(n+m) \bigr\rfloor }$. The cases when exactly one category is empty would be even more
straightforward to analyse and are thus left out from the discussion for brevity.)

Without loss of generality we can assume that the labelling of the indices is such that $\eps_i = \eps_i'$, $i \in \scope{k}$.
(If not, simply relabel the scale functions.) Let us introduce the $k$ number of ($N+1$)-dimensional local variables
$\tilde{y}_i = (y_i,s_i)$ and corresponding scale functions $\tilde{\eps}_i = \eps_i = \eps_i'$, $i \in \scope{k}$.
In category (ii) there are now $n+m-2k$ separated scales to take care of. Introduce the $n+m-2k$ local ``ghost'' variables
$\{ y_i \}_{i=k+1}^{n+m-k}$ and $\{ s_j \}_{j=k+1}^{n+m-k}$ such that one can form the $n+m-2k$ number of ($N+1$)-dimensional
local variables $\tilde{y}_i = (y_j,s_j)$ and scale parameters $\tilde{\eps}_i = \eps_j$ (if $s_j$ where $j \in \scope{k+1,m}$
is the ``ghost'') or $\tilde{\eps}_i = \eps_j'$ (if $y_j$ where $j \in \scope{k+1,n}$ is the ``ghost'') for
$i \in \scope{k+1,n+m-k}$. (Of course, here it is assumed that $k \in \scope{\min \{n,m \}-1}$. If this is not true we simply
introduce ``ghosts'' of only spatial type (i.e., if $k = m < n$) or temporal type (i.e., if $k = n < m$).)

In total we have introduced a local variable
\begin{equation*}
	\tilde{\vecy}_{n+m-k} = (\underbrace{\tilde{y}_1,\ldots,\tilde{y}_k}_{\smatris{\textrm{contains} \\ \textrm{no ``ghosts''}}} \, ,
		\, \underbrace{\tilde{y}_{k+1},\ldots,\tilde{y}_{n+m-k}}_{\smatris{ \textrm{contains} \\ n+m-2k \textrm{ ``ghosts''} }}).
\end{equation*}
which belongs to $\tilde{\sety}^{n+m-k} = (Y_1 \times S_1) \times \cdots \times (Y_{n+m-k} \times S_{n+m-k})$. Define
$\tilde{x} = (x,t)$ and $\tilde{\Omega} = \Omega_T$ such that $\tilde{x} \in \tilde{\Omega}$ for $(x,t) \in \Omega_T$,
and $\tilde{\vecx}_{n+m-k}^{\eps} = (\tfrac{\tilde{x}}{\tilde{\eps}_1},\ldots,\tfrac{\tilde{x}}{\tilde{\eps}_{n+m-k}})$.
Furthermore, given an arbitrary test function $v \in L^2 \bigl( \Omega_T; \cont_{\#} (\sety_{nm}) \bigr)$, let
\begin{equation*}
	\tilde{u}_{\eps}(\tilde{x}) = u_{\eps}(x,t)
		\quad \textrm{and} \quad
	\tilde{v} (\tilde{x},\tilde{\vecy}_{n+m-k}) = v (x,t,\vecy_n,\vecs_m)
\end{equation*}
for all $\tilde{\Omega} \ni \tilde{x} = (x,t) \in \Omega_T$ and all $\tilde{\sety}^{n+m-k} \ni \tilde{\vecy}_{n+m-k} = (\vecy_n,\vecs_m) \in \sety_{nm}$.
We realise that since $v$ is independent of the $n+m-2k$ local ``ghost'' variables, $\tilde{v}$ is too, and we equivalently have that
$\tilde{v} \in L^2 \bigl( \tilde{\Omega}; \cont_{\#} (\tilde{\sety}^{n+m-k}) \bigr)$.

We have by definition
\begin{equation*}
	\smint_{\Omega_T} u_{\eps}(x,t) \, v(x,t,\vecx_n^{\eps},\vect_m^{\eps}) \, \rmd x \rmd t
			= \smint_{\tilde{\Omega}} \tilde{u}_{\eps}(\tilde{x}) \, \tilde{v} (\tilde{x},\tilde{\vecx}_{n+m-k}^{\eps}) \, \rmd \tilde{x}.
\end{equation*}
According to Theorem~2.4 in \cite{AllBri96}, up to a subsequence, $\{ \tilde{u}_{\eps} \}$ ($n+m-k+1$)-converges to a limit
$\tilde{u}_0 \in L^2 (\tilde{\Omega} \times \tilde{\sety}^{n+m-k})$, i.e.,
\begin{equation*}
	\smint_{\tilde{\Omega}} \tilde{u}_{\eps}(\tilde{x}) \, \tilde{v} (\tilde{x},\tilde{\vecx}_{n+m-k}^{\eps}) \, \rmd \tilde{x}
		\rightarrow \smint_{\tilde{\Omega}} \smint_{\tilde{\sety}^{n+m-k}} \tilde{u}_0 (\tilde{x},\tilde{\vecy}_{n+m-k}) \,
			\tilde{v} (\tilde{x},\tilde{\vecy}_{n+m-k}) \, \rmd \tilde{\vecy}_{n+m-k} \rmd \tilde{x}.
\end{equation*}
It is clear that $\tilde{u}_0$ does not depend on the local ``ghost'' variables which implies that there exists
$u_0 \in L^2 (\Omega_T \times \sety_{nm})$ such that
\begin{equation*}
	\tilde{u}_0 (\tilde{x},\tilde{\vecy}_{n+m-k}) = u_0 (x,t,\vecy_n,\vecs_m)
\end{equation*}
for all $\tilde{\Omega} \ni \tilde{x} = (x,t) \in \Omega_T$ and all $\tilde{\sety}^{n+m-k} \ni \tilde{\vecy}_{n+m-k} = (\vecy_n,\vecs_m) \in \sety_{nm}$.
If $\vecy_{\textrm{gh}}$ collects the local ``ghost'' variables and $\sety_{\textrm{gh}}$ is the corresponding local set such that
$\vecy_{\textrm{gh}} \in \sety_{\textrm{gh}}$,
\begin{align*}
	&\smint_{\tilde{\Omega}} \smint_{\tilde{\sety}^{n+m-k}} \tilde{u}_0 (\tilde{x},\tilde{\vecy}_{n+m-k}) \,
		\tilde{v} (\tilde{x},\tilde{\vecy}_{n+m-k}) \, \rmd \tilde{\vecy}_{n+m-k} \rmd \tilde{x} \\
	&\phantom{ \smint_{\tilde{\Omega}} \smint_{\tilde{\sety}^{n+m-k}} \tilde{u}_0 (\tilde{x}, }
		= \smint_{\Omega_T} \smint_{\sety_{nm}} \; \smint_{\sety_{\textrm{gh}}} \! u_0 (x,t,\vecy_n,\vecs_m) \rmd \vecy_{\textrm{gh}} \;
			v (x,t,\vecy_n,\vecs_m) \rmd \vecs_m \rmd \vecy_n \, \rmd x \rmd t \\
	&\phantom{ \smint_{\tilde{\Omega}} \smint_{\tilde{\sety}^{n+m-k}} \tilde{u}_0 (\tilde{x}, }
		= \smint_{\Omega_T} \smint_{\sety_{nm}} u_0 (x,t,\vecy_n,\vecs_m) \, v (x,t,\vecy_n,\vecs_m) \rmd \vecs_m \rmd \vecy_n \, \rmd x \rmd t.
\end{align*}
To conclude, we have shown that
\begin{equation*}
	\smint_{\Omega_T} u_{\eps}(x,t) \, v(x,t,\vecx_n^{\eps},\vect_m^{\eps}) \, \rmd x \rmd t
			\rightarrow  \smint_{\Omega_T} \smint_{\sety_{nm}} u_0 (x,t,\vecy_n,\vecs_m)
					\, v (x,t,\vecy_n,\vecs_m) \rmd \vecs_m \rmd \vecy_n \, \rmd x \rmd t
\end{equation*}
for all $v \in L^2 \bigl( \Omega_T; \cont_{\#} (\sety_{nm}) \bigr)$ where $u_0 \in L^2 (\Omega_T \times \sety_{nm})$.
This means precisely that, for the extracted subsequence, $u_{\eps} \scaleconv{(n+1,m+1)} u_0$, and we are done.
\end{proof}
The proposition below states that under certain restrictions for $v$ defined on $\Omega_T \times \sety_{nm}$,
the sequence $\{ v^{\eps} \}$ converges weakly in $L^2 (\Omega_T)$ to the average over the local variables. 
\begin{prop}
\label{prop:weakconveps}
Suppose that the pair $\bigl( \{ \eps_i \}_{i=1}^n , \{ \eps_j' \}_{j=1}^m \bigr)$ of lists of spatial and temporal
scale functions belongs to $\jointly_{\sep}^{nm}$. Then
\begin{equation}
	\label{eq:weakconvveps}
	v^{\eps} \rightharpoonup \smint_{\sety_{nm}} v(\cdot,\vecy_n,\vecs_m)
		\, \rmd \vecs_m \rmd \vecy_n \qquad \qquad \textrm{in } L^2 (\Omega_T)
\end{equation}
for every $v \in \cont \bigl( \bar{\Omega}_T; \cont_{\#} (\sety_{nm}) \bigr)$.
\end{prop}
\begin{proof}[\boldproof]
Proceed as in the first part of the proof of Theorem~\ref{th:nm-scaleconvres}---i.e., introducing quantities expressed
with tilde---but letting $v \in \cont \bigl( \bar{\Omega}_T; \cont_{\#} (\sety_{nm}) \bigr)$ instead. Now we have
introduced a collection of $n+m-2k$ local ``ghost'' variables collected in the variable $\vecy_{\textrm{gh}} \in \sety_{\textrm{gh}}$.
For every $\lambda > 0$, let $\{ \tilde{K}_{\mu}^{\lambda} \}_{\mu=1}^M$ be a covering of $\tilde{\Omega}$ where
$\tilde{K}_{\mu}^{\lambda}$ are cubes of side length $\tfrac{1}{\lambda}$ such that $\tilde{K}_{\mu}^{\lambda}
\cap \tilde{\Omega} \neq \emptyset$. Moreover, introduce $\tilde{x}_{\mu}^{\lambda} \in \tilde{K}_{\mu}^{\lambda}$,
$\mu \in \scope{M}$. According to the convergence result of Lemma~4.2.2 in \cite{Pank97}, for any given
$v \in \cont \bigl( \bar{\Omega}_T; \cont_{\#} (\sety_{nm}) \bigr)$ and fixed $\mu \in \scope{M}$, it holds that
\begin{equation}
	\label{eq:convfirstterm}
	\smint_{\tilde{\Omega}} \tilde{v}(\tilde{x}_{\mu}^{\lambda},\tilde{\vecx}_{n+m-k}^{\eps})
		\, \tilde{\phi} (\tilde{x}) \, \rmd \tilde{x} \rightarrow  \smint_{\tilde{\Omega}} \smint_{\tilde{\sety}^{n+m-k}}
			\tilde{v}(\tilde{x}_{\mu}^{\lambda},\tilde{\vecy}_{n+m-k}) \, \tilde{\phi} (\tilde{x}) \, \rmd \tilde{\vecy}_{n+m-k} \rmd \tilde{x}
\end{equation}
for all $\tilde{\phi} \in L^2 (\tilde{\Omega})$ since $\tilde{v} (\tilde{x}_{\mu}^{\lambda},\cdot) \in
\cont_{\#} (\tilde{\sety}^{n+m-k}) \subset L_{\#}^2(\tilde{\sety}^{n+m-k})$.

Now, define the simple function (with respect to $\tilde{x} \in \tilde{\Omega}$)
\begin{equation*}
	\tilde{v}^{\lambda} (\tilde{x},\tilde{\vecy}_{n+m-k})
		= \ssum_{\mu=1}^M \tilde{v} (\tilde{x}_{\mu}^{\lambda},\tilde{\vecy}_{n+m-k})
			\, \chi_{\tilde{K}_{\mu}^{\lambda} \cap \tilde{\Omega}} (\tilde{x})
				\qquad \qquad (\tilde{x} \in \tilde{\Omega}, \tilde{\vecy}_{n+m-k} \in \tilde{\sety}^{n+m-k}),
\end{equation*}
where $\chi_A$ is the characteristic function on $A \subset \real^{N+1}$, and
\begin{equation*}
	\tilde{\delta}^{\lambda} (\tilde{x}) = \sup_{\tilde{\sety}^{n+m-k}}
		\bigl| (\tilde{v} - \tilde{v}^{\lambda}) (\tilde{x},\tilde{\vecy}_{n+m-k}) \bigr|.
\end{equation*}
Note that for every fixed $\tilde{x} \in \tilde{\Omega}$, the difference $(\tilde{v} - \tilde{v}^{\lambda}) (\tilde{x},\cdot)$
is uniformly continuous on $\tilde{\sety}^{n+m-k}$. This means in particular that the supremum above can be taken
over any countable dense subset of $\tilde{\sety}^{n+m-k}$ like, e.g., $\tilde{\sety}^{n+m-k} \cap \ratio^{(n+m-k)(N+1)}$.
We observe that $\tilde{\delta}^{\lambda}$ is the supremum of a countable family of measurable functions,
and in virtue of claim (9a) on p.~1012 in \cite{ZeidIIB} this implies that $\tilde{\delta}^{\lambda}$ itself is measurable as well.
The strong regularity of $\tilde{v}$ guarantees that
\begin{equation*}
	\tilde{\delta}^{\lambda} (\tilde{x}) \rightarrow 0
\end{equation*}
as $\lambda \rightarrow \infty$ for every fixed $\tilde{x} \in \tilde{\Omega}$. Furthermore,
we clearly have a majoriser
\begin{align*}
	\bigl| \tilde{\delta}^{\lambda} (\tilde{x})	\bigr|	&\leqslant \sup_{\tilde{\Omega}
		\times \tilde{\sety}^{n+m-k}} \bigl| \tilde{v} (\tilde{x},\tilde{\vecy}_{n+m-k}) \bigr| + \sup_{\tilde{\Omega}
			\times \tilde{\sety}^{n+m-k}} \bigl| \tilde{v}^{\lambda} (\tilde{x},\tilde{\vecy}_{n+m-k}) \bigr| \\
			&\leqslant 2 \sup_{\tilde{\Omega} \times \tilde{\sety}^{n+m-k}} \bigl| \tilde{v} (\tilde{x},\tilde{\vecy}_{n+m-k}) \bigr|
\end{align*}
(i.e., a constant majoriser). Hence, according to Lebesgue's dominated convergence theorem, we have shown that
\begin{equation*}
	\smint_{\tilde{\Omega}} \tilde{\delta}^{\lambda} (\tilde{x})
		\, \rmd \tilde{x} \rightarrow \smint_{\tilde{\Omega}} 0 \, \rmd \tilde{x} = 0.
\end{equation*}

We get the estimation
\begin{align*}
	\Bigl| \smint_{\tilde{\Omega}}	&\tilde{v} (\tilde{x},\tilde{\vecx}_{n+m-k}^{\eps}) \, \phi (\tilde{x}) \, \rmd \tilde{x}
		- \smint_{\tilde{\Omega}} \smint_{\tilde{\sety}^{n+m-k}} \tilde{v} (\tilde{x},\tilde{\vecy}_{n+m-k}) \,
		\phi (\tilde{x}) \, \rmd \tilde{\vecy}_{n+m-k} \rmd \tilde{x} \Bigr| \\
									&\leqslant
		\Bigl| \smint_{\tilde{\Omega}} \tilde{v}^{\lambda} (\tilde{x},\tilde{\vecx}_{n+m-k}^{\eps}) \, \phi (\tilde{x}) \, \rmd \tilde{x}
		- \smint_{\tilde{\Omega}} \smint_{\tilde{\sety}^{n+m-k}} \tilde{v}^{\lambda} (\tilde{x},\tilde{\vecy}_{n+m-k}) \,
		\phi (\tilde{x}) \, \rmd \tilde{\vecy}_{n+m-k} \rmd \tilde{x} \Bigr| \\
									&\quad
		+ \Bigl| \smint_{\tilde{\Omega}} (\tilde{v} - \tilde{v}^{\lambda}) (\tilde{x},\tilde{\vecx}_{n+m-k}^{\eps})
			\, \phi(\tilde{x}) \, \rmd \tilde{x} \Bigr| + \Bigl| \smint_{\tilde{\Omega}} (\tilde{v}^{\lambda} - \tilde{v})
				(\tilde{x},\tilde{\vecy}_{n+m-k}) \, \phi(\tilde{x}) \, \rmd \tilde{\vecy}_{n+m-k} \rmd \tilde{x} \Bigr|
\end{align*}
for every $\phi \in \distr (\tilde{\Omega})$. The convergence result \eqref{eq:convfirstterm} implies that the first
term tends to zero. For any fixed $\eps > 0$, the middle and last terms are both majorised by $\tilde{\delta}^{\lambda}$,
which in the limit $\lambda \rightarrow \infty$ means that these terms vanish. Thus, we have proven that for every given
$v \in \cont \bigl( \bar{\Omega}_T; \cont_{\#} (\sety_{nm}) \bigr)$,
\begin{equation}
	\label{eq:convrestilde}
	\smint_{\tilde{\Omega}}	\tilde{v} (\tilde{x},\tilde{\vecx}_{n+m-k}^{\eps}) \, \tilde{\phi} (\tilde{x}) \, \rmd \tilde{x}
		\rightarrow \smint_{\tilde{\Omega}} \smint_{\tilde{\sety}^{n+m-k}} \tilde{v} (\tilde{x},\tilde{\vecy}_{n+m-k}) \,
			\tilde{\phi} (\tilde{x}) \, \rmd \tilde{\vecy}_{n+m-k} \rmd \tilde{x}
\end{equation}
for all $\tilde{\phi} \in \distr (\tilde{\Omega})$. Since $\tilde{v}(\cdot,\tilde{\vecy}_{n+m-k})$ is a bounded function
in $L^2 (\tilde{\Omega})$ for every $\tilde{\vecy}_{n+m-k} \in \tilde{\sety}^{n+m-k}$, the convergence \eqref{eq:convrestilde}
also holds for all $\tilde{\phi} \in L^2 (\tilde{\Omega})$.

Define $\phi$ by
\begin{equation*}
	\phi(x,t) = \tilde{\phi}(\tilde{x}) \qquad \qquad (\Omega_T \ni (x,t) = \tilde{x} \in \tilde{\Omega}).
\end{equation*}
Then $\tilde{\phi} \in L^2 (\tilde{\Omega})$ is equivalent to saying that $\phi \in L^2 (\Omega_T)$. The convergence result
\eqref{eq:convrestilde} is thus realised to mean that for every given
$v \in \cont \bigl( \bar{\Omega}_T; \cont_{\#} (\sety_{nm}) \bigr)$,
\begin{equation*}
	\smint_{\Omega_T} v (x,t,\vecx_n^{\eps},\vect_m^{\eps}) \, \phi (x,t) \, \rmd x \rmd t
		\rightarrow \smint_{\Omega_T} \smint_{\sety_{nm}} v (x,t,\vecy_n,\vecs_m) \,
			\phi (x,t) \, \rmd \vecs_m \rmd \vecy_n \rmd x \rmd t
\end{equation*}
for all $\phi \in L^2 (\Omega_T)$. Hence, we have shown \eqref{eq:weakconvveps}, and the proof is complete.
\end{proof}
\begin{prop}
\label{prop:strongtonm}
Suppose that the pair $\bigl( \{ \eps_i \}_{i=1}^n , \{ \eps_j' \}_{j=1}^m \bigr)$ of lists of spatial and temporal
scale functions belongs to $\jointly_{\sep}^{nm}$. Moreover, assume that $\{ u_{\eps} \}$ converges strongly
to $u$ in $L^2 (\Omega_T)$. Then $u_{\eps} \scaleconv{(n+1,m+1)} u$.
\end{prop}
\begin{proof}[\boldproof]
From Proposition~\ref{prop:weakconveps} we have
\begin{equation}
	v^{\eps} \rightharpoonup \smint_{\sety_{nm}} v(\cdot,\vecy_n,\vecs_m)
		\, \rmd \vecs_m \rmd \vecy_n \qquad \qquad \textrm{in } L^2 (\Omega_T)
\end{equation}
for every $v \in \cont \bigl( \bar{\Omega}_T; \cont_{\#} (\sety_{nm}) \bigr)$. This combined with the assumption
\begin{equation*}
	u_{\eps} \rightarrow u \qquad \qquad \textrm{in } L^2 (\Omega_T)
\end{equation*}
implies
\begin{equation*}
	\smint_{\Omega_T} u_{\eps} (x,t) \, v (x,t,\vecx_n^{\eps},\vect_m^{\eps}) \, \rmd x \rmd t
		\rightarrow \smint_{\Omega_T} \smint_{\sety_{nm}} u(x,t)
			\, v(x,t,\vecy_n,\vecs_m) \, \rmd \vecs_m \rmd \vecy_n \rmd x \rmd t
\end{equation*}
for every $v \in \cont \bigl( \bar{\Omega}_T; \cont_{\#} (\sety_{nm}) \bigr)
\subset L^2 \bigl( \Omega_T; \cont_{\#} (\sety_{nm}) \bigr)$, where we have used the weak--strong convergence
theorem with respect to $L^2 (\Omega_T)$. Due to Proposition~\ref{prop:convinsmallspace} this convergence
in fact holds for all functions $v \in L^2 \bigl( \Omega_T; \cont_{\#} (\sety_{nm}) \bigr)$ due to the inclusion
$\cont \bigl( \bar{\Omega}_T; \cont_{\#} (\sety_{nm}) \bigr)
\supset \distr \bigl( \Omega_T; \cont_{\#}^{\infty} (\sety_{nm}) \bigr)$. Hence, $u_{\eps} \scaleconv{(n+1,m+1)} u$,
and we are done.
\end{proof}
For the next theorem, Theorem~\ref{th:gradchar} concerning multiscale convergence of gradient sequences, we need the
two lemmas below. Note first that we introduce the following notations. We write $Y^{\scope{i_1,i_2}} = Y_{i_1}
\times \cdots \times Y_{i_2}$ and $S^{\scope{j_1,j_2}} = S_{j_1} \times \cdots \times S_{j_2}$.
Moreover, $\vecy_{\scope{i_1,i_2}} \in Y^{\scope{i_1,i_2}}$ and $\vecs_{\scope{j_1,j_2}} \in S^{\scope{j_1,j_2}}$
are the corresponding local variables. The Lebesgue measures on the introduced local sets are written accordingly.
Furthermore, we define $\wrum_k = H^1_{\#}(Y_k)/\real$, $k \in \scope{n}$, for brevity. It should be emphasised that all
derivatives are taken in the weak (or distributional or generalised) sense.
\begin{lemma}
\label{lm:DefOfH}
Let $\mathcal{H}$ be the subspace of generalised divergence-free functions in $L^2 ( \Omega \times Y^n )^N$ defined by
\begin{align*}
	\mathcal{H} = \Bigl\{	&\psi \in L^2 ( \Omega \times Y^n )^N \, : \, \grad_{y_n} \cdot \psi = 0
							\, \textrm{ and } \, \smint_{Y^{\scope{k+1,n}}} \grad_{y_k} \cdot \psi (x,\vecy_n) \,
								\rmd \vecy_{\scope{k+1,n}} = 0 \\
							&\textrm{ for all } x \in \Omega, \, \vecy_k \in Y^k \textrm{ and all } k \in \scope{n-1} \Bigr\}.
\end{align*}
Then the subspace $\mathcal{H}$ has the following properties:
\begin{itemize}
	\item[(i)] The intersection $\distr \bigl( \Omega; \cont_{\#}^{\infty} (Y^n)^N \bigr) \cap \mathcal{H}$
				is dense in $\mathcal{H}$;
	\item[(ii)] The orthogonal complement $\mathcal{H}^{\perp}$ in $L^2 ( \Omega \times Y^n )^N$ of
				$\mathcal{H}$ is
				\begin{equation*}
					\mathcal{H}^{\perp} = \Bigl\{ \ssum_{k=1}^n \grad_{y_k} u_k \, :
						\, u_k \in L^2 ( \Omega \times Y^{k-1}; \wrum_k ) \Bigr\}.
				\end{equation*}
\end{itemize}
\end{lemma}
\begin{proof}[\boldproof]
See Lemma~3.7 in \cite{AllBri96}.
\end{proof}
\begin{lemma}
\label{lm:phiinebound}
Let $k \in \scope{n}$ and suppose that the list $\{ \eps_i \}_{i=1}^n$ is well-separated. Furthermore, introduce
\begin{equation*}
	\mathcal{E}_k = \Bigl\{ \phi \in \distr \bigl( \Omega; \cont_ {\#}^{\infty} (Y^n) \bigr) \, :
						\, \smint_{Y^{\scope{k,n}}} \phi ( x,\vecy_n ) \, \rmd \vecy_{\scope{k,n}} = 0
							\textrm{ for all } x \in \Omega, \, \vecy_{k-1} \in Y^{k-1} \Bigl\}.
\end{equation*}
Then, for any function $\phi \in \mathcal{E}_k$, the sequence $\bigl\{ \tfrac{1}{\eps_k} \phi^{\eps} \bigr\}$
is bounded in $H^{-1}(\Omega)$.
\end{lemma}
\begin{proof}[\boldproof]
See Corollary~3.4 in \cite{AllBri96}.
\end{proof}
For the $(n+1,m+1)$-scale convergence of sequences of gradients we have the important Theorem~\ref{th:gradchar} below.
\begin{theo}
\label{th:gradchar}
Suppose that the pair $\bigl( \{ \eps_i \}_{i=1}^n , \{ \eps_j' \}_{j=1}^m \bigr)$ of lists of spatial and temporal
scale functions belongs to $\jointly_{\wsep}^{nm}$. Moreover, assume that $\{ u_{\eps} \}$ is a bounded
sequence in $H^1 \bigl( \ZeToT;H^1_0(\Omega),H^{-1}(\Omega) \bigr)$. Then, up to a subsequence, we have
\begin{align*}
	u_{\eps}		&\rightarrow u 			\qquad \qquad \textrm{in } L^2(\Omega_T), \\
	u_{\eps}		&\rightharpoonup  u 	\qquad \qquad \textrm{in } L^2 \bigl( \ZeToT;H^1_0(\Omega) \bigr), \\
	\inter{and}
	\grad u_{\eps}	&\scaleconv{(n+1,m+1)} \grad u + \ssum_{k=1}^n \grad_{y_k} u_k,
\end{align*}
where $u \in L^2 \bigl( \ZeToT;H^1_0(\Omega) \bigr)$ and $u_k \in L^2 ( \Omega_T \times \sety_{(k-1)m}; \wrum_k )$
for all $k \in \scope{n}$.
\end{theo}
\begin{proof}[\boldproof]
Since $\{ u_{\eps} \}$ is bounded in $H^1 \bigl( \ZeToT; H_0^1 (\Omega), H^{-1} (\Omega) \bigr)$,
(i) $\{ u_{\eps} \}$ is also bounded in $L^2 \bigl( \ZeToT; H_0^1 (\Omega) \bigr)$,
(ii) $\{ \tfrac{\d}{\d t} u_{\eps}\}$ is bounded in $L^2 \bigl( \ZeToT; H^{-1} (\Omega) \bigr)$ and
(iii) $\{ \grad u_{\eps} \}$ is bounded in $L^2 ( \Omega_T )^N$. The first statement (i) implies,
up to a subsequence,
\begin{equation*}
	u_{\eps} \rightharpoonup u \qquad \qquad \textrm{in } L^2 \bigl( \ZeToT; H_0^1 (\Omega) \bigr)
\end{equation*}	
for some unique $u \in L^2 \bigl( \ZeToT;H^1_0(\Omega) \bigr)$. By Lemmas~8.2~and~8.4 in
\cite{ConFoi88}, the statements (i) and (ii) imply, up to a subsequence, that
\begin{equation}
	\label{eq:strongconv}
	u_{\eps} \rightarrow u 			\qquad \qquad \textrm{in } L^2(\Omega_T).
\end{equation}
Hence, we have proven the convergences for $u_{\eps}$.

From Theorem~\ref{th:nm-scaleconvres} and (i) and (iii) we know that , up to a subsequence,
\begin{equation}
	\label{eq:gradnm-conv}
	\grad u_{\eps} \scaleconv{(n+1,m+1)} w_0
\end{equation}
for some limit function $w_0 \in L^2 (\Omega_T \times \sety_{nm})^N$.

We will now characterise $w_0$ in terms of gradients. Using the vector valued product test
function $\psi \in L^2 \bigl( \Omega_T; \cont_{\#} (\sety_{nm}) \bigr)$ defined by
\begin{equation*}
	\psi (x,t,\vecy_n,\vecs_m) = v(x,\vecy_n) \, c(t,\vecs_m)
\end{equation*}
for all $(x,t) \in \Omega_T$ and all $(\vecy_n,\vecs_m) \in \sety_{nm}$, where $v \in \distr \bigl( \Omega; \cont_{\#}^{\infty} (Y^n) \bigr) \cap \mathcal{H}$ and
$c \in \distr (\ZeToT) \odot \cont_{\#}^{\infty} (S^m)$, in the $(n+1,m+1)$-scale convergence
result \eqref{eq:gradnm-conv} yields, up to a subsequence,
\begin{align}
	\smint_{\Omega_T} \grad u_{\eps} (x,t) \cdot
			&v(x,\vecx_n^{\eps}) \, c(t,\vect_m^{\eps}) \, \rmd x \rmd t \nonumber \\
		&\rightarrow \smint_{\Omega_T} \smint_{\sety_{nm}} w_0 (x,t,\vecy_n,\vecs_m) \cdot
			v(x,\vecy_n) \, c(t,\vecs_m) \, \rmd \vecs_m \rmd \vecy_n \rmd x \rmd t. \label{eq:gradprodconv}
\end{align}
Taking a closer look at the left-hand side of \eqref{eq:gradprodconv} we get
\begin{align*}
	&\smint_{\Omega_T} \grad u_{\eps} (x,t) \cdot v(x,\vecx_n^{\eps}) \, c(t,\vect_m^{\eps}) \, \rmd x \rmd t \\
	&\phantom{ \smint_{\Omega_T} \grad u_{\eps} (x,t) \cdot v(x,\vecx_n^{\eps}) }
		= - \smint_{\Omega_T} u_{\eps} (x,t) \, \Bigl( \grad_x
			+ \ssum_{k=1}^n \tfrac{1}{\eps_k} \grad_{y_k} \Bigr)
				\cdot v(x,\vecx_n^{\eps}) \, c(t,\vect_m^{\eps}) \, \rmd x \rmd t \\
	&\phantom{ \smint_{\Omega_T} \grad u_{\eps} (x,t) \cdot v(x,\vecx_n^{\eps}) }
		= - \smint_{\Omega_T} u_{\eps} (x,t) \, \Bigl( \grad_x
			+ \ssum_{k=1}^{n-1} \tfrac{1}{\eps_k} \grad_{y_k} \Bigr)
				\cdot v(x,\vecx_n^{\eps}) \, c(t,\vect_m^{\eps}) \, \rmd x \rmd t,
\end{align*}
where we in the first equality have have used partial integration on $\Omega$, divergence theorem on $\Omega$ and
the fact that both (though only one is necessary) $u_{\eps}$ and $v$ vanish on $\d \Omega$, and in the second
equality used the fact that $v \in \mathcal{H}$ implying $\grad_{y_n} \cdot v = 0$. We claim now that
$\grad_{y_k} \cdot v \in \mathcal{E}_{k+1}$, $k \in \scope{n-1}$. Indeed, for any $k \in \scope{n-1}$ we have $\grad_{y_k}
\cdot v \in \distr \bigl( \Omega; \cont_{\#}^{\infty} (Y^n) \bigr)$ and
\begin{equation*}
	\smint_{Y^{\scope{k+1,n}}} \grad_{y_k} \cdot v(x,\vecy_n) \, \rmd \vecy_{\scope{k+1,n}} = 0, \qquad x \in \Omega, \, \vecy_k \in Y^k,
\end{equation*}
where we have simply employed the definition of $v$ being in $\mathcal{H}$ making the
multiple integral vanish, so $\grad_{y_k} \cdot v \in \mathcal{E}_{k+1}$.
Thus, by Lemma~\ref{lm:phiinebound} we have that $\bigl\{ \tfrac{1}{\eps_{k+1}} \grad_{y_k} \cdot v^{\eps} \bigr\}$
is bounded in $H^{-1}(\Omega)$ for all $k \in \scope{n-1}$. This boundedness yields an estimation
\begin{align*}
	\Bigl| \smint_{\Omega_T} u_{\eps} (x,t) \, &\ssum_{k=1}^{n-1} \tfrac{1}{\eps_k} \grad_{y_k}
			\cdot v(x,\vecx_n^{\eps}) \, c(t,\vect_m^{\eps}) \, \rmd x \rmd t \Bigr|^2 \\
												&\leqslant
			T \smint_0^T \Bigl| \smint_{\Omega} u_{\eps} (x,t) \, \ssum_{k=1}^{n-1} \tfrac{1}{\eps_k} \grad_{y_k}
			\cdot v(x,\vecx_n^{\eps}) \, c(t,\vect_m^{\eps}) \, \rmd x \Bigr|^2 \rmd t \\
												&\leqslant
			T \smint_0^T \Bigl| \Bigl\langle \ssum_{k=1}^{n-1} \tfrac{1}{\eps_k} \grad_{y_k} \cdot v^{\eps},
				u_{\eps}(t) c(t,\vect_m^{\eps}) \Bigr\rangle_{H^{-1}(\Omega),H_0^1(\Omega)} \Bigr|^2 \rmd t \\
												&\leqslant
			T \smint_0^T \Bigl\| \ssum_{k=1}^{n-1} \tfrac{1}{\eps_k} \grad_{y_k} \cdot v^{\eps} \Bigr\|_{H^{-1}(\Omega)}^2
				\bigl\| u_{\eps}(t) c(t,\vect_m^{\eps}) \bigr\|_{H_0^1(\Omega)}^2 \rmd t,
\end{align*}
i.e.,
\begin{align*}
	\Bigl| \smint_{\Omega_T} u_{\eps} (x,t) \, &\ssum_{k=1}^{n-1} \tfrac{1}{\eps_k} \grad_{y_k}
			\cdot v(x,\vecx_n^{\eps}) \, c(t,\vect_m^{\eps}) \, \rmd x \rmd t \Bigr|^2 \\
												&\leqslant
			C_1 \Bigl( \ssum_{k=1}^{n-1} \tfrac{\eps_{k+1}}{\eps_k} \bigl\| \tfrac{1}{\eps_{k+1}} \grad_{y_k}
				\cdot v^{\eps} \bigr\|_{H^{-1}(\Omega)} \Bigr)^2
					\smint_0^T  \bigl\| u_{\eps}(t) \bigr\|_{H_0^1(\Omega)}^2 | c(t,\vect_m^{\eps}) | \, \rmd t \\
												&\leqslant
			C_2 \Bigl( \ssum_{k=1}^{n-1} \tfrac{\eps_{k+1}}{\eps_k} \Bigr)^2
				\smint_0^T \bigl\| u_{\eps}(t) \bigr\|_{H_0^1(\Omega)}^2 \rmd t
												=
			C_2 \Bigl( \ssum_{k=1}^{n-1} \tfrac{\eps_{k+1}}{\eps_k} \Bigr)^2 
				\bigl\| u_{\eps} \bigr\|_{L^2 \bigl( \ZeToT; H_0^1(\Omega) \bigr)}^2 \\
												&\leqslant
			C_3 \Bigl( \ssum_{k=1}^{n-1} \tfrac{\eps_{k+1}}{\eps_k} \Bigr)^2 \rightarrow 0
\end{align*}
since the scale functions are separated. We thus conclude that the left-hand side of \eqref{eq:gradprodconv}
converges according to
\begin{align*}
	&\smint_{\Omega_T} \grad u_{\eps} (x,t) \cdot v(x,\vecx_n^{\eps}) \, c(t,\vect_m^{\eps}) \, \rmd x \rmd t \\
	&\phantom{ \smint_{\Omega_T} \grad u_{\eps} (x,t) \cdot v(x,\vecx_n^{\eps}) }
		\rightarrow - \smint_{\Omega_T} \smint_{\sety_{nm}} u (x,t)
			\, \grad_x \cdot v(x,\vecy_n) \, c(t,\vecs_m) \, \rmd \vecs_m \rmd \vecy_n \rmd x \rmd t \\
	&\phantom{ \smint_{\Omega_T} \grad u_{\eps} (x,t) \cdot v(x,\vecx_n^{\eps}) }
		= \smint_{\Omega_T} \smint_{\sety_{nm}} \grad u (x,t) \, \cdot v(x,\vecy_n)
			\, c(t,\vecs_m) \, \rmd \vecs_m \rmd \vecy_n \rmd x \rmd t
\end{align*}
for all $v \in \distr \bigl( \Omega; \cont_{\#}^{\infty} (Y^n) \bigr) \cap \mathcal{H}$ and all
$c \in \distr (\ZeToT) \odot \cont_{\#}^{\infty} (S^m)$. From the right-hand side of of
\eqref{eq:gradprodconv} we thus obtain
\begin{equation*}
	\smint_{\Omega_T} \smint_{\sety_{nm}} \bigl( w_0 (x,t,\vecy_n,\vecs_m) - \grad u (x,t) \bigr)
		\cdot v(x,\vecy_n) \, c(t,\vecs_m) \, \rmd \vecs_m \rmd \vecy_n \rmd x \rmd t = 0,
\end{equation*}
or
\begin{equation*}
	\smint_0^T \smint_{S^m} \Bigl( \smint_{\Omega} \smint_{Y^n}
		\bigl( w_0 (x,t,\vecy_n,\vecs_m) - \grad u (x,t) \bigr) \cdot v(x,\vecy_n)
			\, \rmd \vecy_n \rmd x \Bigr) \, c(t,\vecs_m) \, \rmd \vecs_m \rmd t = 0.
\end{equation*}
By the Variational Lemma and utilising density (i.e., (i) in Lemma~\ref{lm:DefOfH}), for every $v \in \mathcal{H}$
it holds that
\begin{equation*}
	\smint_{\Omega} \smint_{Y^n} \bigl( w_0 (x,t,\vecy_n,\vecs_m) - \grad u (x,t) \bigr)
		\cdot v(x,\vecy_n) \, \rmd \vecy_n \rmd x = 0
\end{equation*}
a.e.~on $\tint \times S^m$. Hence,
\begin{equation*}
	w_0 - \grad u \; \perp \; v \qquad \qquad \textrm{in }
		L^2 ( \Omega \times Y^n )^N \textrm{ a.e.~on } \tint \times S^m,
\end{equation*}
i.e., $w_0 - \grad u \in \mathcal{H}^{\perp}$ a.e.~on $\tint \times S^m$. Thus, by (ii) in
Lemma~\ref{lm:DefOfH},
\begin{equation*}
	w_0 - \grad u = \ssum_{k=1}^n \grad_{y_k} u_k
		\qquad \qquad \textrm{ a.e.~on } \tint \times S^m,
\end{equation*}
where $u_k \in L^2 ( \Omega \times Y^{k-1}; \wrum_k )$ a.e.~on $\tint \times S^m$.

What remains is to prove that $u_k \in L^2 ( \Omega_T \times \sety_{(k-1)m}; \wrum_k )$,
$k \in \scope{n}$. We will perform a proof by induction accomplished in two steps: the
Base case followed by the Inductive step.

\emph{Base case.} We must show that $u_1 \in L^2 ( \Omega_T \times S^m; \wrum_1 )$. We have,
a.e.~on $\Omega_T \times \sety_{1m}$,
\begin{align}
	\grad_{y_1} u_1 (x,t,y_1,\vecs_m)
			&= \smint_{Y^{\scope{2,n}}} \grad_{y_1} u_1 (x,t,y_1,\vecs_m)
				\, \rmd \vecy_{\scope{2,n}} \nonumber \\
			&= \smint_{Y^{\scope{2,n}}} \ssum_{i=1}^n \grad_{y_i} u_i (x,t,\vecy_{i},\vecs_m)
				\, \rmd \vecy_{\scope{2,n}} \nonumber \\
			&= \smint_{Y^{\scope{2,n}}} \bigl( w_0 (x,t,\vecy_n,\vecs_m) - \grad u (x,t) \bigr)
				\, \rmd \vecy_{\scope{2,n}} \nonumber \\
			&= \smint_{Y^{\scope{2,n}}} w_0 (x,t,\vecy_n,\vecs_m)
				\, \rmd \vecy_{\scope{2,n}} - \grad u (x,t), \label{eq:basecase}
\end{align}
where the second equality follows from the fact that $u_i$ is $Y_i$-periodic.
Hence, by \eqref{eq:basecase} and the well-known characterisation of the $\wrum_1$-norm in terms of an
$L^2$-norm of the gradient (see, e.g., Proposition~3.52 in \cite{CioDon99}),
\begin{align}
	\| u_1 \|_{L^2 ( \Omega_T \times S^m; \wrum_1 )}
			&= \| \grad_{y_1} u_1 \|_{L^2(\Omega_T \times \sety_{1m})^N} \nonumber \\
			&= \bigl\| \smint_{Y^{\scope{2,n}}} w_0
				- \grad u \bigr\|_{L^2(\Omega_T \times \sety_{1m})^N} \nonumber \\
			&\leqslant \bigl\| \smint_{Y^{\scope{2,n}}} w_0 \bigr\|_{L^2(\Omega_T \times \sety_{1m})^N}
				+ \| \grad u \|_{L^2(\Omega_T \times \sety_{1m})^N}. \label{eq:basecasenorm}
\end{align}
Since $w_0 \in L^2 (\Omega_T \times \sety_{nm})^N$, we have that $\smint_{Y^{\scope{2,n}}} w_0
\in L^2 (\Omega_T \times \sety_{1m})^N$, and since $u \in L^2 \bigl( \ZeToT; H_0^1(\Omega) \bigr)$, it holds that
$\grad u \in L^2 (\Omega_T)^N \subset L^2 (\Omega_T \times \sety_{1m})^N$. Thus, by \eqref{eq:basecasenorm},
\begin{equation*}
	\| u_1 \|_{L^2 ( \Omega_T \times S^m; \wrum_1 )} < \infty, 
\end{equation*}
which means that $u_1 \in L^2 ( \Omega_T \times S^m; \wrum_1 )$ as desired; the Base case is complete.

\emph{Inductive step.} Assume that $u_j \in L^2 ( \Omega_T \times \sety_{(j-1)m}; \wrum_j )$ for all
$j \in \scope{\ell}$ where $\ell \in \scope{n-1}$ (requires $n > 1$; the case $n = 1$ is already treated
in the Base case above). We must show that this assumption implies $u_{\ell + 1}
\in L^2 \bigl( \Omega_T \times \sety_{\ell m}; \wrum_{\ell + 1} \bigr)$. If $\ell \in \scope{n-2}$ we have,
a.e.~on $\Omega_T \times \sety_{(\ell+1)m}$,
\begin{align}
	&\grad_{y_{\ell+1}} u_{\ell+1} (x,t,\vecy_{\ell+1},\vecs_m) \nonumber \\
		&\phantom{ \grad_{y_{\ell+1}} u }
			= \smint_{Y^{\scope{\ell+2,n}}} \grad_{y_{\ell+1}} u_{\ell+1} (x,t,\vecy_{\ell+1},\vecs_m)
				\, \rmd \vecy_{\scope{\ell+2,n}} \nonumber \\
		&\phantom{ \grad_{y_{\ell+1}} u }
			= \smint_{Y^{\scope{\ell+2,n}}} \ssum_{i=1}^n \grad_{y_i} u_i (x,t,\vecy_{i},\vecs_m)
				\, \rmd \vecy_{\scope{\ell+2,n}} - \smint_{Y^{\scope{\ell+2,n}}} \ssum_{i=1}^{\ell}
					\grad_{y_i} u_i (x,t,\vecy_{i},\vecs_m) \, \rmd \vecy_{\scope{\ell+2,n}} \nonumber \\
		&\phantom{ \grad_{y_{\ell+1}} u }
			= \smint_{Y^{\scope{\ell+2,n}}} \bigl( w_0 (x,t,\vecy_n,\vecs_m) - \grad u (x,t) \bigr)
				\, \rmd \vecy_{\scope{\ell+2,n}} - \ssum_{i=1}^{\ell} \grad_{y_i} u_i (x,t,\vecy_{i},\vecs_m) \nonumber \\
		&\phantom{ \grad_{y_{\ell+1}} u }
			= \smint_{Y^{\scope{\ell+2,n}}} w_0 (x,t,\vecy_n,\vecs_m) \, \rmd \vecy_{\scope{\ell+2,n}}
				- \grad u (x,t) - \ssum_{i=1}^{\ell} \grad_{y_i} u_i (x,t,\vecy_{i},\vecs_m), \label{eq:indstep}
\end{align}
where the second equality follows from the fact that $u_i$ is $Y_i$-periodic.
If we in this proof interpret integration over ``$Y^{\scope{n+1,n}}$'' as performing no integration at all
(i.e., $\int_{Y^{\scope{n+1,n}}} w_0 = w_0$ by definition), \eqref{eq:indstep} actually works for $\ell = n-1$
as well. We get the norm
\begin{align}
	\| u_{\ell+1} \|_{L^2 ( \Omega_T \times \sety_{\ell m}; \wrum_{\ell+1} )}
			&= \| \grad_{y_{\ell+1}} u_{\ell+1} \|_{L^2 ( \Omega_T \times \sety_{(\ell+1) m} )^N} \nonumber \\
			&= \bigl\| \smint_{Y^{\scope{\ell+2,n}}} w_0 - \grad u
				- \ssum_{i=1}^{\ell} \grad_{y_i} u_i \bigr\|_{L^2 ( \Omega_T \times \sety_{(\ell+1) m} )^N} \nonumber \\
			&\leqslant \bigl\| \smint_{Y^{\scope{\ell+2,n}}}
				w_0 \bigr\|_{L^2 ( \Omega_T \times \sety_{(\ell+1) m} )^N} \nonumber \\
			&\qquad \quad + \| \grad u \|_{L^2 ( \Omega_T \times \sety_{(\ell+1) m} )^N}
				+ \ssum_{i=1}^{\ell} \| \grad_{y_i} u_i \|_{L^2 ( \Omega_T \times \sety_{(\ell+1) m} )^N} \nonumber \\
			&= \bigl\| \smint_{Y^{\scope{\ell+2,n}}}
				w_0 \bigr\|_{L^2 ( \Omega_T \times \sety_{(\ell+1) m} )^N} \nonumber \\
			&\qquad \quad + \| \grad u \|_{L^2 ( \Omega_T \times \sety_{(\ell+1) m} )^N}
					+ \ssum_{i=1}^{\ell} \| u_i \|_{L^2 ( \Omega_T \times \sety_{(i-1)m}; \wrum_i )}. \label{eq:indstepnorm}
\end{align}
Since $w_0 \in L^2 (\Omega_T \times \sety_{nm})^N$, we have that $\smint_{Y^{\scope{\ell+2,n}}} w_0
\in L^2 (\Omega_T \times \sety_{(\ell+1)m})^N$, and since $u \in L^2 \bigl( \ZeToT; H_0^1(\Omega) \bigr)$, it holds that
$\grad u \in L^2 (\Omega_T)^N \subset L^2 (\Omega_T \times \sety_{(\ell+1)m})^N$. By the inductive assumption,
$u_j \in L^2 ( \Omega_T \times \sety_{(j-1)m}; \wrum_j )$ for all $j \in \scope{\ell}$. Thus, \eqref{eq:indstepnorm} gives
\begin{equation*}
	\| u_{\ell+1} \|_{L^2 ( \Omega_T \times \sety_{\ell m}; \wrum_{\ell+1} )} < \infty,
\end{equation*}
which means that $u_{\ell + 1} \in L^2 \bigl( \Omega_T \times \sety_{\ell m}; \wrum_{\ell + 1} \bigr)$ as desired;
the Inductive step is complete, and we are done.
\end{proof}
When performing the homogenisation later in this paper we will limit ourselves to two spatial scales, $n = 1$,
where the microscale is described by the single spatial scale function $\eps_1$. The scale function $\eps_1$ is, without
loss of generality, assumed to coincide with the scale parameter, i.e., $\eps_1(\eps) = \eps$.
Note that in what follows, the list $\{ \eps \}$ of the single spatial scale function will be written as $\eps$ for brevity.
In this setting we have Theorem~\ref{th:nmtypelimit} below. We first need a lemma.
\begin{lemma}
\label{eq:hypoellipticity}
Suppose $\rho \in \cont_{\#}^{\infty}(Y)/\real$. Then there exists a unique $\theta \in \cont_{\#}^{\infty}(Y)/\real$
such that $\rho = \Delta_y \theta$ where $\Delta_y$ is the Laplace operator with respect to $y$ (i.e.,
$\Delta_y = \grad_y \cdot \grad_y$).
\end{lemma}
\begin{proof}[\boldproof]
First we note that for any given $\rho \in L_{\#}^2(Y)/\real$ there exists a unique function $\theta \in H_{\#}^1(Y)/\real$
such that $\rho = \Delta_y \theta$. Then we consider only smooth source functions $\rho \in \cont_{\#}^{\infty}(Y)/\real
\subset L_{\#}^2(Y)/\real$ and utilise the hypoellipticity property of the Laplace operator to conclude that $\theta$
must also belong to $\cont_{\#}^{\infty}(Y)/\real$. (For a further discussion see, e.g., Remark~3.2 in \cite{NguWou07}.)
\end{proof}
In the remainder of the paper, let $\wrum = H_{\#}^1(Y)/\real$.
\begin{theo}
\label{th:nmtypelimit}
Suppose that the pair $\bigl( \eps , \{ \eps_i' \}_{i=1}^m \bigr)$ of lists of spatial and temporal scale functions
belongs to $\jointly_{\wsep}^{1m}$ and assume that $\{ u_{\eps} \}$ is a bounded sequence in the function space
$H^1 \bigl( \ZeToT;H^1_0(\Omega),H^{-1}(\Omega) \bigr)$. Then, up to a subsequence,
\begin{equation}
	\label{eq:veryweakconv}
	\smint_{\Omega_T} \tfrac{1}{\eps} u_{\eps}(x,t) \, \phi (x,t,\tfrac{x}{\eps},\vect_m^{\eps}) \, \rmd x \rmd t
		\rightarrow \smint_{\Omega_T} \smint_{\sety_{1m}} u_1(x,t,y,\vecs_m)
			\, \phi (x,t,y,\vecs_m) \, \rmd \vecs_m \rmd y \rmd x \rmd t 
\end{equation}
for all $\phi \in \distr(\Omega) \! \odot \! \distr (\ZeToT) \odot \bigl( \cont^{\infty}_{\#}(Y)/\real \bigr)
\! \odot \! \bigl( \prod_{i=1}^m \cont^{\infty}_{\#}(S_i) \bigr)$, where $u_1 \in L^2 (\Omega_T \times S^m; \wrum)$
is as in Theorem~\ref{th:gradchar} (with $n = 1$).
\end{theo}
\begin{proof}[\boldproof]
Fix an arbitrary $\phi \in \distr(\Omega) \! \odot \! \distr (\ZeToT) \odot
\bigl( \cont^{\infty}_{\#}(Y)/\real \bigr) \! \odot \! \bigl( \prod_{i=1}^m \cont^{\infty}_{\#}(S_i) \bigr)$.
Then there exist unique $\psi \in \distr(\Omega) \! \odot \! \distr (\ZeToT) \, \odot \,
\bigl( \prod_{i=1}^m \cont^{\infty}_{\#}(S_i) \bigr)$ and $\rho \in \cont^{\infty}_{\#}(Y)/\real$ such that
$\phi = \psi \rho$. The left-hand side of \eqref{eq:veryweakconv} can then be written
\begin{align*}
	\smint_{\Omega_T} \tfrac{1}{\eps} u_{\eps}(x,t) \, \phi (x,t,\tfrac{x}{\eps},\vect_m^{\eps}) \, \rmd x \rmd t
		&= \smint_{\Omega_T} \tfrac{1}{\eps} u_{\eps}(x,t) \psi (x,t,\vect_m^{\eps})
			\, \rho (\tfrac{x}{\eps}) \, \rmd x \rmd t \\
		&= \smint_{\Omega_T} \tfrac{1}{\eps} u_{\eps}(x,t) \psi (x,t,\vect_m^{\eps})
			\, \Delta_y \theta (\tfrac{x}{\eps}) \, \rmd x \rmd t \\
		&= \smint_{\Omega_T} u_{\eps}(x,t) \psi (x,t,\vect_m^{\eps}) \, \tfrac{1}{\eps}
			\grad_y \cdot \grad_y \theta (\tfrac{x}{\eps}) \, \rmd x \rmd t
\end{align*}
for some unique $\theta \in \cont_{\#}^{\infty}(Y)/\real$ due to Lemma~\ref{eq:hypoellipticity}. By noting that
\begin{equation*}
	\grad \cdot \sigma (\tfrac{x}{\eps}) = \tfrac{1}{\eps} \grad_y \cdot \sigma (\tfrac{x}{\eps})
\end{equation*}
for any $\sigma$ differentiable over $Y$ (here $\sigma = \grad_y \theta$),
we get by partial integration on $\Omega$ that
\begin{align*}
	\smint_{\Omega_T} \tfrac{1}{\eps} u_{\eps}(x,t) \, \phi (x,t,\tfrac{x}{\eps},\vect_m^{\eps}) \, \rmd x \rmd t
		&= \smint_{\Omega_T} u_{\eps}(x,t) \psi (x,t,\vect_m^{\eps})
			\, \grad \cdot \grad_y \theta (\tfrac{x}{\eps}) \, \rmd x \rmd t \\
		&= \smint_{\Omega_T} \Bigl( \grad \cdot \bigl( u_{\eps}(x,t) \psi (x,t,\vect_m^{\eps})
			\, \grad_y \theta (\tfrac{x}{\eps}) \bigr) \\
		&\phantom{ = \smint_{\Omega_T} \Bigl(  }
			- \grad u_{\eps}(x,t) \, \psi (x,t,\vect_m^{\eps}) \cdot \grad_y \theta (\tfrac{x}{\eps}) \\
		&\phantom{ = \smint_{\Omega_T} \Bigl( }
			- u_{\eps}(x,t) \, \grad \psi (x,t,\vect_m^{\eps})
				\cdot \grad_y \theta (\tfrac{x}{\eps}) \Bigr) \, \rmd x \rmd t \\
		&= - \smint_{\Omega_T} \grad u_{\eps}(x,t) \, \psi (x,t,\vect_m^{\eps})
			\cdot \grad_y \theta (\tfrac{x}{\eps}) \, \rmd x \rmd t \\
		&\qquad \qquad - \smint_{\Omega_T} u_{\eps}(x,t) \, \grad \psi (x,t,\vect_m^{\eps})
			\cdot \grad_y \theta (\tfrac{x}{\eps}) \, \rmd x \rmd t,
\end{align*}
where we in the last equality for the first term in the integrand have employed the divergence theorem on
$\Omega$ and used the fact that both (though only one is necessary) $u_{\eps}$ and $\psi$ vanish on $\d \Omega$.
Furthermore, by utilisation of Theorem~\ref{th:gradchar} with $n = 1$, we get (with $u \in
L^2 \bigl( \ZeToT; H^1_0 (\Omega) \bigr)$ and $u_1 \in L^2 (\Omega_T \times S^m; \wrum)$ as in
Theorem~\ref{th:gradchar} with $n = 1$, up to a subsequence,
\begin{align*}
	&\smint_{\Omega_T} \tfrac{1}{\eps} u_{\eps}(x,t) \, \phi (x,t,\tfrac{x}{\eps},\vect_m^{\eps}) \, \rmd x \rmd t \\
		&\qquad \rightarrow  - \smint_{\Omega_T} \smint_{\sety_{1m}} \bigl( \grad u(x,t)
			+ \grad_y u_1 (x,t,y,\vecs_m) \bigr) \, \psi (x,t,\vecs_m)
				\cdot \grad_y \theta (y) \, \rmd \vecs_m \rmd y \rmd x \rmd t \\
		&\phantom{ \qquad \rightarrow  - \smint_{\Omega_T} \smint_{\sety_{1m}} \bigl( \grad u(x,t)
			+ \grad_y } - \smint_{\Omega_T} \smint_{\sety_{1m}} u(x,t) \, \grad \psi (x,t,\vecs_m)
				\cdot \grad_y \theta (y) \, \rmd \vecs_m \rmd y \rmd x \rmd t \\
		&\qquad = \smint_{\Omega_T} \smint_{\sety_{1m}} u(x,t) \, \grad \psi (x,t,\vecs_m)
			\cdot \grad_y \theta (y) \, \rmd \vecs_m \rmd y \rmd x \rmd t \\
		&\qquad \qquad \qquad + \smint_{\Omega_T} \smint_{\sety_{1m}} u_1(x,t,y,\vecs_m)
			\, \psi (x,t,\vecs_m) \, \grad_y \cdot \grad_y \theta (y) \, \rmd \vecs_m \rmd y \rmd x \rmd t \\
		&\qquad \qquad \qquad \qquad \qquad \qquad \qquad - \smint_{\Omega_T} \smint_{\sety_{1m}} u(x,t)
			\, \grad \psi (x,t,\vecs_m) \cdot \grad_y \theta (y) \, \rmd \vecs_m \rmd y \rmd x \rmd t \\
		&\qquad = \smint_{\Omega_T} \smint_{\sety_{1m}} u_1(x,t,y,\vecs_m)
			\, \psi (x,t,\vecs_m) \, \rho (y) \, \rmd \vecs_m \rmd y \rmd x \rmd t \\
		&\qquad = \smint_{\Omega_T} \smint_{\sety_{1m}} u_1(x,t,y,\vecs_m)
			\, \phi (x,t,y,\vecs_m) \, \rmd \vecs_m \rmd y \rmd x \rmd t ,
\end{align*}
where we have performed a partial integration on $\Omega$ and $Y$, respectively, of the first integral followed
by using the divergence theorem and noting that $u$ and $\psi$ vanish on $\d \Omega$ and that $u_1$ and $\theta$
are $Y$-periodic (giving a vanishing surface integral over~$\d Y$); we have derived \eqref{eq:veryweakconv}.
Since $\phi \in \distr(\Omega) \! \odot \! \distr (\ZeToT) \, \odot \, \bigl( \cont^{\infty}_{\#}(Y)/\real \bigr)
\! \odot \! \bigl( \prod_{i=1}^m \cont^{\infty}_{\#}(S_i) \bigr)$ was arbitrary, the claim of the theorem follows.
\end{proof}
\begin{rema}
\label{re:nmtypelimit}
Theorem~\ref{th:nmtypelimit} is a mere variety of Lemma~3.1 in \cite{NguWou07} in the special case of periodicity
but generalised to include several temporal scales. In its turn, the result in \cite{NguWou07} is a mere variation of
Corollary~3.3 in \cite{Holm97} generalised to the non-periodic setting and with the sequence $\bigl\{ \tfrac{1}{\eps}
u_{\eps} \bigr\}$ (as in Theorem~\ref{th:nmtypelimit} above) instead of the slightly more complicated sequence
$\bigl\{ \tfrac{1}{\eps} (u_{\eps} - u) \bigr\}$ found in \cite{Holm97}.

The convergence mode in Theorem~\ref{th:nmtypelimit} can be regarded as a kind of feeble, or ``very weak'',
$(2,m+1)$-scale convergence of $\bigl\{ \tfrac{1}{\eps} u_{\eps} \bigr\}$ since the heavily restricted set of
test functions in question is more permissible compared to the larger set of test functions employed in ordinary
$(2,m+1)$-scale convergence.

Finally, we remark that a result analogous to Proposition~\ref{prop:asymptnmscaleconv} holds for sequences
of the type $\bigl\{ \tfrac{1}{\eps} r(\eps) u_{\eps} \bigr\}$ having a ``very weak'' limit $r_0 u_1$ instead
of $u_1$ if $r(\eps) \rightarrow r_0$.
\end{rema}

\section{Monotone Parabolic Operators}

Consider the operator-form evolution problem
\begin{equation}
	\label{eq:operatoreq}
	\left\{
		\begin{aligned}
			\tfrac{\rmd}{\rmd t} u + \aoper u	&= f, \\
			u(0)				&= u^0 \in H, \\
			u 					&\in H^1 (\ZeToT;V,V'),
		\end{aligned}
	\right.
\end{equation}
where $f \in L^2 \bigl( \ZeToT; V' \bigr)$ and $\aoper \, : \, L^2 \bigl( \ZeToT; V \bigr) \rightarrow
L^2 \bigl( \ZeToT; V' \bigr)$. Here $H$ is some Hilbert space and $V$ is some Banach space with topological
dual $V'$. Note that $u \in H^1 (\ZeToT;V,V')$ means $u \in L^2 \bigl( \ZeToT; V \bigr)$ and $\tfrac{\rmd}{\rmd t} u
\in L^2 \bigl( \ZeToT; V' \bigr)$, $\tfrac{\rmd}{\rmd t}$ being the weak (or distributional or generalised)
derivative with respect to the temporal variable $t \in \tint$. The definition below establishes a convenient
relation between $H$, $V$ and $V'$.
\begin{defin}
Suppose $H$ is a real and separable Hilbert space and that $V$ is a real, separable and reflexive
Banach space such that $V$ is continuously embedded and dense in $H$. We then call $(V,H,V')$ an evolution triple.
\end{defin}
\begin{rema}
Recall that $V$ is continuously embedded in $H$ if $V \subset H$ and there exists $C > 0$ such that
$\| u \|_H \leqslant C \| u \|_V$ for all $u \in V$. Also note that by Riesz's representation theorem,
$H$ can be identified by its dual $H'$ and that $H'$ is continuously embedded and dense in $V'$.
Schematically we have
\begin{equation*}
	V \overset{\smatris{\textrm{Cont.~emb.} \\ \textrm{\& dense}}}{\subset} H \overset{\smatris{\textrm{Riesz's} \\ 
		\textrm{repr.~th.}}}{\sim} H' \overset{\smatris{\textrm{Cont.~emb.} \\ \textrm{\& dense}}}{\subset} V'.
\end{equation*}
\end{rema}
Let for every fixed $t \in \tint$ the operator $\aoper (t) \, : \, V \rightarrow V'$ be defined by
\begin{equation}
	\label{eq:timedeptimeindep}
	\aoper (t) u(t) =  (\aoper u)(t) \qquad \qquad (u \in L^2(\ZeToT;V)).
\end{equation}
In order for the problem \eqref{eq:operatoreq} to have a unique solution the operator $\aoper$ should satisfy the
following five conditions:
\begin{itemize}
	\item[(A$_1$)] $\bigl\langle \aoper(t) u - \aoper(t) v \, , \, u - v \bigr\rangle_{V',V} \geqslant 0$
				for all $u,v \in V$ and all $t \in \tint$ (i.e., $\aoper (t)$ is monotone);
	\item[(A$_2$)] The $[0,1] \rightarrow \real$ function $q \mapsto \bigl\langle \aoper(t) (u + qw),
				v \bigr\rangle_{V',V}$ is continuous for all $u,v,w \in V$ and all $t \in \tint$
				(i.e., $\aoper (t)$ is hemicontinuous);
	\item[(A$_3$)] There exists $C_0 > 0$ such that $\bigl\langle \aoper(t) u,u \bigr\rangle_{V',V}
				\geqslant C_0 \| u \|_V^2$ for all $u \in V$ and all $t \in \tint$ (i.e.,
				$\aoper (t)$ is coercive);
	\item[(A$_4$)] There exist a non-negative function $\beta \in L^2(\ZeToT)$ and a constant $C_1 > 0$
				such that $\| \aoper (t) u \|_{V'} \leqslant \beta(t) + C_1 \| u \|_V$ for all $u \in V$
				and all $t \in \tint$ (i.e., $\aoper$ satisfies a certain growth condition);
	\item[(A$_5$)] The $\tint \rightarrow \real$ function $t \mapsto \bigl\langle \aoper(t) u,
				v \bigr\rangle_{V',V}$ is measurable on $\tint$ for all $u,v \in V$ (i.e., $t \mapsto
				\aoper (t)$ is weakly measurable on $\tint$).
\end{itemize}
We have the following theorem on existence and uniqueness.
\begin{theo}
\label{th:fromzeidler}
Suppose that $\aoper \, : \, L^2 \bigl( \ZeToT; V \bigr) \rightarrow L^2 \bigl( \ZeToT; V' \bigr)$
satisfies \emph{(A}$_1$\emph{)}--\emph{(A}$_5$\emph{)} above and assume that $(V,H,V')$ forms an
evolution triple. Then, for every $f \in L^2 (\ZeToT;V')$ and $u^0 \in H$, there exists a unique
solution $u$ to \eqref{eq:operatoreq}.
\end{theo}
\begin{proof}[\boldproof]
See Theorem~30.A in \cite{ZeidIIB}.
\end{proof}
Let $X = L^2(\ZeToT;V)$ and $X' = L^2(\ZeToT;V')$, and consider a sequence $\{ \aoper^{\eps} \}$ of
monotone operators. Equivalently to \eqref{eq:operatoreq} for this sequence of operators, the evolution
problem can be formulated as
\begin{equation}
	\label{eq:weakevoleq}
	\left\{
		\begin{aligned}
			\bigl\langle \tfrac{\rmd}{\rmd t} u_{\eps},v \bigr\rangle_{X',X}
						+ \bigl\langle \aoper^{\eps} u_{\eps},v \bigr\rangle_{X',X}
					&= \bigl\langle f,v \bigr\rangle_{X',X} \;, \\
			u(0)	&= u^0 \in H, \\
			u_{\eps}		&\in H^1 (\ZeToT;V,V')
		\end{aligned}
	\right.
\end{equation}
for all $v \in X = L^2(\ZeToT;V)$, where $u^0 \in H$, $f \in X' = L^2(\ZeToT;V')$ and $(V,H,V')$
is an evolution triple.

Fix $H = L^2 (\Omega)$ and $V = H^1_0 (\Omega)$ with dual $V' = H^{-1} (\Omega)$. Then
\begin{equation*}
	\bigl( H^1_0 (\Omega), L^2 (\Omega), H^{-1} (\Omega)\bigr)
\end{equation*}
is an evolution triple. We let the operators $\aoper^{\eps} \, : \, L^2 \bigl( \ZeToT; H^1_0 (\Omega) \bigr)
\rightarrow L^2 \bigl( \ZeToT; H^{-1} (\Omega) \bigr)$ be defined in terms of a flux $a^{\eps} \, : \,
\bar{\Omega}_T \times \real^N \rightarrow \real^N$ by
\begin{equation}
	\label{eq:operdeffromflux}
	\bigl\langle \aoper^{\eps} u,v \bigr\rangle_{X',X}
		= \smint_{\Omega_T} a^{\eps}(x,t; \grad u) \cdot \grad v (x,t) \, \rmd x \rmd t
\end{equation}
for $u,v \in X = L^2 \bigl( \ZeToT;H^1_0 (\Omega) \bigr)$, which---by the definition \eqref{eq:timedeptimeindep}
of the time dependent operator---is the same as
\begin{align*}
		\bigl\langle \aoper^{\eps}(t) u,v \bigr\rangle_{H^{-1} (\Omega),H^1_0 (\Omega)}
			= \smint_{\Omega} a^{\eps}(x,t; \grad u)
				\cdot \grad v (x) \, \rmd x
\end{align*}
for $u,v \in H^1_0 (\Omega)$. We recall that $a^{\eps}$ is given via $a \, : \, \bar{\Omega}_T \times \real^{nN+m}
\times \real^N \rightarrow \real^N$ according to
\begin{equation}
	\label{eq:aepsdef}
	a^{\eps} (x,t; k) = a(x,t,\vecx_n^{\eps},\vect_m^{\eps}; k) \qquad \qquad ((x,t) \in \Omega_T, \, k \in \real^N).
\end{equation}
The problem
\begin{equation}
	\label{eq:monoparaprobgen}
	\left\{
		\begin{aligned}
			\tfrac{\d}{\d t} u_{\eps} (x,t) - \grad \cdot a^{\eps}(x,t; \grad u_{\eps})
							&= f(x,t)	&	&\textrm{in } \Omega_T, \\
			u_{\eps}(x,0)	&= u^0(x)	&	&\textrm{in } \Omega, \\
			u_{\eps}(x,t)	&= 0		&	&\textrm{on }\d \Omega \times \tint,
		\end{aligned}
	\right.
\end{equation}
is the same as \eqref{eq:monoparaprob} but generalised to $n+1$ spatial scales.
Clearly, with the conventions above, \eqref{eq:weakevoleq} is the weakly formulated version of \eqref{eq:monoparaprobgen}.
To conclude, the weak formulation is that, given $f \in X' = L^2 \bigl( \ZeToT; H^{-1} (\Omega) \bigr)$ and $u^0 \in L^2 (\Omega)$,
we want to find $u_{\eps} \in H^1 \bigl( \ZeToT; H^1_0 (\Omega), H^{-1} (\Omega) \bigr)$ such that
\begin{equation}
	\label{eq:weakformofevprob}
	\bigl\langle \tfrac{\d}{\d t} u_{\eps} , v \bigr\rangle_{X',X}
		+ \smint_{\Omega_T} a(x,t,\vecx_n^{\eps},\vect_m^{\eps}; \grad u_{\eps}) \cdot \grad v (x,t) \, \rmd x \rmd t
			= \smint_{\Omega_T} f(x,t) \, v (x,t) \, \rmd x \rmd t
\end{equation}
for all $v \in X = L^2 \bigl( \ZeToT; H^1_0 (\Omega) \bigr)$. The function $a$ should satisfy
the following five structure conditions:
\begin{itemize}
	\item[(B$_1$)] $a(x,t,\vecy_n,\vecs_m;0) = 0$ for all $(x,t) \in \bar{\Omega}_T$ and all $(\vecy_n,\vecs_m) \in \real^{nN+m}$;
	\item[(B$_2$)] $a(x,t,\, \cdot \; ;k)$ is $\sety_{nm}$-periodic for all $(x,t) \in \bar{\Omega}_T$ and all $k \in \real^N$,
				and $a(\, \cdot \; ;k)$ is continuous for all $k \in \real^N$;
	\item[(B$_3$)] $a(x,t,\vecy_n,\vecs_m; \, \cdot \, )$ is continuous for all $(x,t) \in \bar{\Omega}_T$
				and all $(\vecy_n,\vecs_m) \in \real^{nN+m}$;
	\item[(B$_4$)] There exists $C_0 > 0$ such that
				\begin{equation*}
					\bigl( a(x,t,\vecy_n,\vecs_m;k) - a(x,t,\vecy_n,\vecs_m;k') \bigr) \cdot (k - k') \geqslant C_0 |k - k'|^2
				\end{equation*}
				for all $(x,t) \in \Omega_T$, all $(\vecy_n,\vecs_m) \in \real^{nN+m}$ and all $k,k' \in \real^N$;
	\item[(B$_5$)] There exist $C_1 > 0$ and $0 < \alpha \leqslant 1$ such that
				\begin{equation*}
					\bigl| a(x,t,\vecy_n,\vecs_m;k) - a(x,t,\vecy_n,\vecs_m;k') \bigr|
						\leqslant C_1 (1 + |k| + |k'|)^{1-\alpha} |k - k'|^{\alpha}
				\end{equation*}
				for all $(x,t) \in \Omega_T$, all $(\vecy_n,\vecs_m) \in \real^{nN+m}$ and all $k,k' \in \real^N$.
\end{itemize}
We have the following proposition linking the structural conditions (B$_1$)--(B$_5$) for $a$ to the conditions
(A$_1$)--(A$_5$) for $\aoper^{\eps}$.
\begin{prop}
\label{prop:connectBA}
Suppose that $a \, : \, \bar{\Omega}_T \times \real^{nN+m} \times \real^N \rightarrow \real^N$ fulfils the structure
conditions \emph{(B}$_1$\emph{)}--\emph{(B}$_5$\emph{)}. Then $\aoper^{\eps} \, : \, L^2 \bigl( \ZeToT; H^1_0 (\Omega) \bigr)
\rightarrow L^2 \bigl( \ZeToT; H^{-1} (\Omega) \bigr)$ defined through \eqref{eq:operdeffromflux} satisfies the conditions
\emph{(A}$_1$\emph{)}--\emph{(A}$_5$\emph{)}.
\end{prop}
\begin{proof}[\boldproof]
We first prove that the monotonicity condition (A$_1$) holds. Fix an arbitrary $t \in \tint$.
Then, for any $u,v = H_0^1 (\Omega)$,
\begin{align*}
	\bigl\langle \aoper^{\eps}(t) u - \aoper^{\eps}&(t) v \, , \, u - v \bigr\rangle_{H^{-1} (\Omega),H_0^1 (\Omega)} \\
		&= \smint_{\Omega} \bigl( a(x,t,\vecx_n^{\eps},\vect_m^{\eps};\grad u) - a(x,t,\vecx_n^{\eps},\vect_m^{\eps};\grad v) \bigr)
			\cdot \bigl( \grad u(x) - \grad v(x) \bigr) \, \rmd x \\
		&\geqslant C_0 \smint_{\Omega} \bigl| \grad u(x) - \grad v(x) \bigr|^2 \, \rmd x
			= C_0 \| u - v \|_{H_0^1 (\Omega)}^2 \\
		&\geqslant 0,
\end{align*}
where we have employed the structure condition (B$_4$) to obtain the first inequality.

Next we prove the hemicontinuity condition (A$_2$). Fix arbitrary $t \in \tint$ and $q_0 \in [0,1]$, and let $q \in [0,1]$.
Then, for any $u,v,w \in H^1_0 (\Omega)$,
\begin{align*}
	&\Bigl| \bigl\langle \aoper^{\eps}(t) (u+qw) , v \bigl\rangle_{H^{-1}(\Omega),H^1_0 (\Omega)}
		- \bigl\langle \aoper^{\eps}(t) (u+q_0 w) , v \bigl\rangle_{H^{-1}(\Omega),H^1_0 (\Omega)} \Bigr| \\
	&\phantom{ \Bigl| \bigl\langle \aoper^{\eps}(t) (u }
		= \Bigl| \smint_{\Omega} \bigl( a(x,t,\vecx_n^{\eps},\vect_m^{\eps};\grad u + q \grad w)
				- a(x,t,\vecx_n^{\eps},\vect_m^{\eps};\grad u + q_0 \grad w) \bigr) \cdot \grad v (x) \, \rmd x \Bigr| \\
	&\phantom{ \Bigl| \bigl\langle \aoper^{\eps}(t) (u }
		\leqslant \smint_{\Omega} \bigl| a(x,t,\vecx_n^{\eps},\vect_m^{\eps};\grad u + q \grad w)
				- a(x,t,\vecx_n^{\eps},\vect_m^{\eps};\grad u + q_0 \grad w) \bigr| \, \bigl| \grad v (x) \bigr| \, \rmd x \\
	&\phantom{ \Bigl| \bigl\langle \aoper^{\eps}(t) (u }
		\leqslant C_1 \smint_{\Omega} \Bigl( 1 + \bigl| \grad u(x) + q \grad w(x) \bigr| \\
	&\phantom{ \Bigl| \bigl\langle \aoper^{\eps}(t) (u }
		\phantom{ \leqslant C_1 \smint_{\Omega} \Bigl( 1 } \; + \bigl| \grad u(x)
				+ q_0 \grad w(x) \bigr| \Bigr)^{1-\alpha} \bigl| (q-q_0) \grad w(x) \bigr|^{\alpha}
					\, \bigl| \grad v (x) \bigr| \, \rmd x \\
	&\phantom{ \Bigl| \bigl\langle \aoper^{\eps}(t) (u }
		\leqslant C_1 |q-q_0|^{\alpha} \smint_{\Omega} \Bigl( 1 + 2\bigl| \grad u(x) \bigr|
				+ 2\bigl| \grad w(x) \bigr| \Bigr)^{1-\alpha} \bigl| \grad w(x) \bigr|^{\alpha}
					\, \bigl| \grad v (x) \bigr| \, \rmd x \\
	&\phantom{ \Bigl| \bigl\langle \aoper^{\eps}(t) (u }
		< 2^{-\alpha} C_1 |q-q_0|^{\alpha} \smint_{\Omega} \Bigl( 1 + 2\bigl| \grad u(x) \bigr|
				+ 2\bigl| \grad w(x) \bigr| \Bigr) \, \bigl| \grad v (x) \bigr| \, \rmd x \\
	&\phantom{ \Bigl| \bigl\langle \aoper^{\eps}(t) (u }
		\leqslant 2^{-\alpha} C_1 \bigl\| \bigl( 1 + 2 | \grad u |
				+ 2 | \grad w | \bigr) \, | \grad v |  \bigr\|_{L^1(\Omega)} |q-q_0|^{\alpha} \\
	&\phantom{ \Bigl| \bigl\langle \aoper^{\eps}(t) (u }
		\leqslant 2^{-\alpha} C_1 \bigl( \| 1 \|_{L^2(\Omega)} + 2 \bigl\| | \grad u | \bigr\|_{L^2(\Omega)}
				+ 2 \bigl\| | \grad w | \bigr\|_{L^2(\Omega)} \bigr)
					\bigl\| | \grad v | \bigr\|_{L^2(\Omega)} |q-q_0|^{\alpha} \\
	&\phantom{ \Bigl| \bigl\langle \aoper^{\eps}(t) (u }
		= 2^{-\alpha} C_1 \bigl( | \Omega |^{\frac{1}{2}} + 2 \| u \|_{H_0^1(\Omega)}
				+ 2 \| w \|_{H_0^1(\Omega)} \bigr) \| v \|_{H_0^1(\Omega)} |q-q_0|^{\alpha} \\
	&\phantom{ \Bigl| \bigl\langle \aoper^{\eps}(t) (u }
		\rightarrow 0
\end{align*}
as $q \rightarrow q_0$, where we have utilised (B$_5$) for the second inequality and H\"{o}lder's inequality together with
the triangle inequality to obtain the last inequality. Thus, the hemicontinuity holds.

Let us move on to proving the coercivity condition (A$_3$). Fix $t \in \tint$. Then, for any $u \in H^1_0 (\Omega)$,
\begin{align*}
	\bigl\langle \aoper^{\eps}(t) u,u \bigr\rangle_{H^{-1}(\Omega),H^1_0 (\Omega)}
			&= \smint_{\Omega} a^{\eps}(x,t; \grad u) \cdot \grad u (x) \, \rmd x \\
			&= \smint_{\Omega} \bigl( a^{\eps}(x,t; \grad u) - a^{\eps}(x,t; 0) \bigr)
				\cdot \bigl( \grad u (x) - 0 \bigr) \, \rmd x \\
			&\geqslant C_0 \smint_{\Omega} \bigl| \grad u (x) - 0 \bigr|^2 \, \rmd x \\
			&= C_0 \| u \|_{H_0^1(\Omega)}^2,
\end{align*}
where we have used structure condition (B$_1$) to obtain the second equality and (B$_4$) for the inequality.

The growth condition (A$_4$) is proven in the following manner. We first note that by (B$_1$) and (B$_5$),
\begin{align}
	\bigl| a(x,t,\vecy_n,\vecs_m;k) \bigr|	&= \bigl| a(x,t,\vecy_n,\vecs_m;k) - a(x,t,\vecy_n,\vecs_m;0) \bigr| \nonumber \\
		&\leqslant C_1 \bigl( 1 + |k| \bigr)^{1-\alpha} |k|^{\alpha} \nonumber \\
		&< C_1 \bigl( 1 + |k| \bigr)^{1-\alpha} \bigl( 1 + |k| \bigr)^{\alpha} \nonumber \\
		&= C_1 \bigl( 1 + |k| \bigr) \label{eq:anineqfora}
\end{align}
for all $(x,t) \in \Omega_T$, all $(\vecy_n,\vecs_m) \in \real^{nN+m}$ and all $k \in \real^N$.
Fix $t \in \tint$ and let $u \in H_0^1(\Omega)$ be arbitrary. Then
\begin{align*}
	\bigl\| \aoper^{\eps} (t) u \bigr\|_{H^{-1}(\Omega)}
		&= \sup_{\| v \|_{H_0^1(\Omega)} \leqslant 1} \Bigl| \bigl\langle \aoper^{\eps} (t) u ,
			v \bigr\rangle_{H^{-1}(\Omega),H_0^1(\Omega)} \Bigr| \\
		&= \sup_{\| v \|_{H_0^1(\Omega)} \leqslant 1} \Bigl| \smint_{\Omega}
			a(x,t,\vecx_n^{\eps},\vect_m^{\eps};\grad u) \cdot \grad v(x) \, \rmd x \Bigr| \\
		&\leqslant \sup_{\| v \|_{H_0^1(\Omega)} \leqslant 1} \smint_{\Omega}
			\bigl| a(x,t,\vecx_n^{\eps},\vect_m^{\eps};\grad u) \bigr| \, \bigl| \grad v(x) \bigr| \, \rmd x \\
		&< C_1 \!\!\! \sup_{\| v \|_{H_0^1(\Omega)} \leqslant 1}
			\smint_{\Omega} \bigl( 1 + \bigl| \grad u (x) \bigr| \bigr) \, \bigl| \grad v(x) \bigr| \, \rmd x \\
		&= C_1 \!\!\! \sup_{\| \, | \grad v | \, \|_{L^2(\Omega)} \leqslant 1}
			\bigl\| \bigl( 1 + | \grad u | \bigr) \, | \grad v | \bigr\|_{L^1(\Omega)} \; ,
\end{align*}
where in the second inequality we have employed \eqref{eq:anineqfora}. By H\"{o}lder's inequality we obtain
\begin{align*}
	\bigl\| \aoper^{\eps} (t) u \bigr\|_{H^{-1}(\Omega)}
		&\leqslant C_1 \!\!\! \sup_{\| \, | \grad v | \, \|_{L^2(\Omega)} \leqslant 1}
			\bigl\| 1 + | \grad u | \bigr\|_{L^2(\Omega)} \bigl\| | \grad v | \bigr\|_{L^2(\Omega)} \\
		&\leqslant C_1 \bigl\| 1 + | \grad u | \bigr\|_{L^2(\Omega)} \leqslant C_1
			\bigl( \| 1 \|_{L^2(\Omega)} + \bigl\| | \grad u | \bigr\|_{L^2(\Omega)} \bigr) \\
		&= \beta + C_1 \| u \|_{H_0^1 (\Omega)}.
\end{align*}
This growth constraint is even more regular than anticipated since $\beta = C_1 \sqrt{| \Omega |}$
is independent of $t \in \tint$.

Finally, the weak measurability condition (A$_5$) follows readily from the continuity assumptions on $a$
and the boundedness property \eqref{eq:anineqfora}.
\end{proof}
The following important theorem follows immediately from Proposition~\ref{prop:connectBA} above together
with Theorem~\ref{th:fromzeidler}.
\begin{theo}
\label{th:existuniquesol}
Suppose that $a \, : \, \bar{\Omega}_T \times \real^{nN+m} \times \real^N \rightarrow \real^N$ fulfils
the structure conditions \emph{(B}$_1$\emph{)}--\emph{(B}$_5$\emph{)}. Then, for every $f \in L^2 (\Omega_T)$
and $u^0 \in L^2 (\Omega)$, the evolution problem \eqref{eq:monoparaprobgen} has a unique weak solution
$u_{\eps} \in H^1 \bigl( \ZeToT; H_0^1(\Omega), H^{-1}(\Omega) \bigr)$.
\end{theo}

\section{H-Convergence of Monotone Parabolic Problems}

In 1967 Spagnolo introduced the notion of G-convergence for linear problems governed by symmetric matrices
(see \cite{Spag67}; see also \cite{Spag68,Spag76,ColSpa77}). The name ``G''-conver\-gence comes from the fact that
this convergence mode corresponds to the convergence of the Green functions associated to the sequence of problems.
The G-convergence of symmetric matrices is defined via the weak convergence of solutions to the sequence of problems.

The concept of H-convergence---``H'' as in ``homogenisation''---is a generalisation of Spagnolo's G-convergence
to cover also non-symmetric matrices. It was introduced in 1976 by Tartar (see \cite{Tart78}; see also \cite{Tart79})
and further developed by Murat in 1978 (see \cite{Mura78,Mura78b}; see also \cite{Mura79}), and in 1977 Tartar defined
H-convergence for non-linear monotone problems (see \cite{Tart77}; see also \cite{ChDaDe90,ChiDef90}). Early studies of
H-convergence for non-linear monotone parabolic problems were conducted by Kun'ch and Pankov in 1986 (see \cite{KunPan86}),
Kun'ch in 1988 (see \cite{Kunch88}), and Svanstedt in 1992 (see \cite{Svan92}; see also
\cite{Svan99} by Svanstedt and \cite{Pank97} by Pankov for further developments).

We introduce a convenient set of flux functions in the following definition.
\begin{defin}
Suppose $C_0,C_1 > 0$ and $0 < \alpha \leqslant 1$. A function $a \, : \, \Omega_T \times \real^N \rightarrow \real^N$
is said to belong to $\mathcal{M}_{C_0,C_1}^{\alpha}(\Omega_T)$ if the following four structure conditions are satisfied:
\begin{itemize}
	\item[\emph{(C}$_1$\emph{)}] $a (x,t;0) = 0$ a.e.~$(x,t) \in \Omega_T$;
	\item[\emph{(C}$_2$\emph{)}] $a(\, \cdot \; ;k)$ is (Lebesgue) measurable for every $k \in \real^N$;
	\item[\emph{(C}$_3$\emph{)}] $\bigl( a(x,t;k) - a(x,t;k') \bigr) \cdot (k - k') \geqslant C_0 |k - k'|^2$
									\; a.e.~$(x,t) \in \Omega_T$ \, and \, for \, all \, $k,k' \in \real^N$;
	\item[\emph{(C}$_4$\emph{)}] $\bigl| a(x,t;k) - a(x,t;k') \bigr| \leqslant
									C_1 \bigl( 1 + |k| + |k'| \bigr)^{1-\alpha} |k - k'|^{\alpha}$
									a.e.~$(x,t) \in \Omega_T$ and for all $k,k' \in \real^N$.
\end{itemize}
If no values on $C_0,C_1 > 0$ and $0 < \alpha \leqslant 1$ are fixed we simply say that
$a \in \mathcal{M}(\Omega_T)$, i.e.,
\begin{equation*}
	\mathcal{M}(\Omega_T) = \mathsmaller{\bigcup\limits}_{\substack{C_0,C_1 > 0 \\ 0 < \alpha \leqslant 1}}
		\mathcal{M}_{C_0,C_1}^{\alpha}(\Omega_T).
\end{equation*}
We may drop $\Omega_T$ as soon as there is no hazard of confusion, i.e., $\mathcal{M}_{C_0,C_1}^{\alpha}
= \mathcal{M}_{C_0,C_1}^{\alpha}(\Omega_T)$ and $\mathcal{M} = \mathcal{M}(\Omega_T)$.
\end{defin}
The important concept of H-convergence of monotone parabolic problems---coined $\mathrm{H}_{\mathrm{MP}}$-convergence
in this paper for brevity---is introduced in the definition below.
\begin{defin}
\label{def:Hmp-conv}
Suppose $\{ a^{\eps} \}$ is a sequence of fluxes in $\mathcal{M}$.
We say that $\{ a^{\eps} \}$ $\mathrm{H}_{\mathrm{MP}}$-converges to the flux $b \in \mathcal{M}$ if,
for any $f \in L^2 \bigl( \ZeToT; H^{-1}(\Omega) \bigr)$ and any $u^0 \in L^2 (\Omega)$, the weak
solutions $u_{\eps} \in H^1 \bigl( \ZeToT; H_0^1(\Omega), H^{-1}(\Omega) \bigr)$ to the sequence
\begin{equation}
	\label{eq:defofseq}
	\left\{
		\begin{aligned}
			\tfrac{\d}{\d t} u_{\eps}(x,t) - \grad \cdot a^{\eps} (x,t;\grad u_{\eps})
							&= f(x,t)	&	&\textrm{in } \Omega_T, \\
			u_{\eps}(x,0)	&= u^0(x)	&	&\textrm{in } \Omega, \\
			u_{\eps}(x,t)	&= 0	&	&\textrm{on } \d \Omega \times \tint
		\end{aligned}
	\right.
\end{equation}
of evolution problems satisfy
\begin{equation*}
	\left\{
		\begin{aligned}
			u_{\eps}						&\rightharpoonup  u					
				&	&\textrm{in } L^2 \bigl( \ZeToT; H_0^1(\Omega) \bigr), \\
			a^{\eps}(\, \cdot \; ;\grad u_{\eps})	&\rightharpoonup  b(\, \cdot \; ;\grad u)
				&	&\textrm{in } L^2 (\Omega_T)^N,
		\end{aligned}
	\right.
\end{equation*}
where $u \in H^1 \bigl( \ZeToT; H_0^1(\Omega), H^{-1}(\Omega) \bigr)$ is the weak unique solution
to the evolution problem
\begin{equation}
	\label{eq:defoflim}
	\left\{
		\begin{aligned}
			\tfrac{\d}{\d t} u(x,t) - \grad \cdot b (x,t;\grad u)
						&= f(x,t)	&	&\textrm{in } \Omega_T, \\
			u(x,0)		&= u^0(x)	&	&\textrm{in } \Omega, \\
			u(x,t)	&= 0	&	&\textrm{on } \d \Omega \times \tint.
		\end{aligned}
	\right.
\end{equation}
Moreover, for brevity, we write this convergence $a^{\eps} \overset{\mathrm{H}_{\mathrm{MP}}}{\longrightarrow} b$,
and $b$ is called the $\mathrm{H}_{\mathrm{MP}}$-limit of $\{ a^{\eps} \}$.
\end{defin}
The definition above leads to the compactness result below.
\begin{theo}
Let $\{ a^{\eps} \}$ be a sequence of fluxes in $\mathcal{M}_{C_0,C_1}^{\alpha}$.
Then, up to a subsequence, $a^{\eps} \overset{\mathrm{H}_{\mathrm{MP}}}{\longrightarrow} b$ for some
$b \in \mathcal{M}_{C_0',C_1'}^{\alpha/(2-\alpha)}$ where $C_0',C_1' > 0$ only depend on the constants
$C_0, C_1, \alpha$.
\end{theo}
\begin{proof}[\boldproof]
This is just a special case of Theorem~3.1 in \cite{Svan99}.
\end{proof}
In the case that $\{ a^{\eps} \}$ is given according to \eqref{eq:aepsdef} we have the following proposition
linking the structure conditions (B$_1$)--(B$_5$) for $a$ to the conditions (C$_1$)--(C$_4$) for the sequence
$\{ a^{\eps} \}$ to be in $\mathcal{M}_{C_0,C_1}^{\alpha}$.
\begin{prop}
\label{prop:connectBwithadmiss}
Suppose that $a \, : \, \bar{\Omega}_T \times \real^{nN+m} \times \real^N \rightarrow \real^N$ fulfils the
structure conditions \emph{(B}$_1$\emph{)}--\emph{(B}$_5$\emph{)}. Then $\{ a^{\eps} \}$ defined through
\eqref{eq:aepsdef} is a sequence in $\mathcal{M}_{C_0,C_1}^{\alpha}$ where
$C_0, C_1$ and $\alpha$ are the constants introduced in \emph{(B}$_1$\emph{)}--\emph{(B}$_5$\emph{)}. 
\end{prop}
\begin{proof}[\boldproof]
We begin by recalling \eqref{eq:aepsdef}, i.e., 
\begin{equation*}
	a^{\eps} (x,t; k) = a(x,t,\vecx_n^{\eps},\vect_m^{\eps}; k) \qquad \qquad ((x,t) \in \Omega_T, \, k \in \real^N).
\end{equation*}

For condition (C$_1$) we have that
\begin{equation*}
	a^{\eps} (x,t; 0)	= a(x,t,\vecx_n^{\eps},\vect_m^{\eps}; 0) = 0
\end{equation*}
for all $(x,t) \in \Omega_T$ by (B$_1$).

Secondly, the (Lebesgue) measurability condition (C$_2$) follows from the
continuity and periodicity properties in condition (B$_2$).

Next we wish to verify (C$_3$). For all $(x,t) \in \Omega_T$ and all $k,k' \in \real^N$,
\begin{align*}
	\bigl( a^{\eps}(x,t;k) - a^{\eps}(x,t;k') \bigr) \cdot (k-k')
		&= \bigl( a(x,t,\vecx_n^{\eps},\vect_m^{\eps}; k) - a(x,t,\vecx_n^{\eps},\vect_m^{\eps}; k') \bigr) \cdot (k-k') \\
		&\geqslant C_0 |k-k'|^2
\end{align*}
according to structure condition (B$_4$).

Finally, (C$_4$) is to be checked. For all $(x,t) \in \Omega_T$ and all $k,k' \in \real^N$,
\begin{align*}
	\bigl| a^{\eps}(x,t;k) - a^{\eps}(x,t;k') \bigr|
		&= \bigl| a(x,t,\vecx_n^{\eps},\vect_m^{\eps}; k) - a(x,t,\vecx_n^{\eps},\vect_m^{\eps}; k') \bigr| \\
		&\leqslant C_1 \bigl( 1 + |k| + |k'| \bigr)^{1-\alpha} |k - k'|^{\alpha}.
\end{align*}
We conclude that $\{ a^{\eps} \}$ is in $\mathcal{M}_{C_0,C_1}^{\alpha}$
where $C_0$, $C_1$ and $\alpha$ are precisely the constants introduced in (B$_1$)--(B$_5$), and we are done.
\end{proof}
We have the following proposition governing an a priori estimate on the solutions to the
sequence of evolution problems.
\begin{prop}
\label{prop:twounifbound}
Suppose that $a \, : \, \bar{\Omega}_T \times \real^{nN+m} \times \real^N \rightarrow \real^N$ fulfils
the structure conditions \emph{(B}$_1$\emph{)}--\emph{(B}$_5$\emph{)}. Then the sequence $\{ u_{\eps}\}$
of weak solutions to the evolution problem \eqref{eq:defofseq} with $\{ a^{\eps} \}$ defined through
\eqref{eq:aepsdef} satisfies
\begin{equation}
	\label{eq:firstunifbound}
	\| u_{\eps} \|_{H^1 \bigl( \ZeToT; H_0^1(\Omega), H^{-1}(\Omega) \bigr)} \leqslant C
\end{equation}
for some $C > 0$. In other words, $\{ u_{\eps}\}$ is uniformly bounded in $H^1 \bigl( \ZeToT; H_0^1(\Omega),
H^{-1}(\Omega) \bigr)$.
\end{prop}
\begin{proof}[\boldproof]
For every fixed $\eps > 0$ we know as a matter of fact that we have a unique weak solution $u_{\eps} \in
H^1 \bigl( \ZeToT; H_0^1(\Omega), H^{-1}(\Omega) \bigr)$ to \eqref{eq:monoparaprobgen} by Theorem~\ref{th:existuniquesol}.

Let us now verify the uniform boundedness in $H^1 \bigl( \ZeToT; H_0^1(\Omega), H^{-1}(\Omega) \bigr)$, i.e.,
\eqref{eq:firstunifbound}. By Proposition~\ref{prop:connectBwithadmiss} we know that $\{ a^{\eps} \}$ is
in $\mathcal{M}$. We can then apply Proposition~2.3 and Corollary~2.1 in \cite{Svan99} which in this context
say that $\{ u_{\eps} \}$ and $\bigl\{ \tfrac{\d}{\d t} u_{\eps} \bigr\}$ are uniformly bounded in
$L^2 \bigl( \ZeToT; H_0^1(\Omega) \bigr)$ and $L^2 \bigl( \ZeToT; H^{-1}(\Omega) \bigr)$, respectively.
Thus, we have uniform boundedness in $H^1 \bigl( \ZeToT; H_0^1(\Omega), H^{-1}(\Omega) \bigr)$, i.e.,
\eqref{eq:firstunifbound} holds.
The proof is complete.
\end{proof}

\section{Homogenisation}

The notion of homogenisation of problems with multiple microscales was introduced in 1978 by Bensoussan, Lions
and Papanicolaou (see \cite{BeLiPa78}) who homogenised problems with two microscales characterised by the list
$\{ \eps,\eps^2 \}$ of scale functions. In 1996, Allaire and Briane (see \cite{AllBri96}) succeeded to generalise
this to homogenisation of linear elliptic problems with an arbitrary number of microscales---even infinitely
many---without even assuming the scale functions to be power functions using the notion of (well-)separatedness
instead. This was achieved by introducing the multiscale convergence technique. In 2001, Lions, Lukkassen, Persson
and Wall performed homogenisation of non-linear monotone elliptic problems with scale functions $\{ \eps,
\eps^2 \}$ (see \cite{LLPW01}), and in 2005 Holmbom, Svanstedt and Wellander studied homogenisation of linear
parabolic problems with pairs $\bigl( \{ \eps, \eps^2 \}, \eps^k \bigr)$ of lists of scale functions (see
\cite{HoSvWe05}). In 2006, Flod\'{e}n and Olsson generalised to monotone parabolic problems (see \cite{FloOls06};
see also \cite{FHOS07} by Flod\'{e}n, Olsson, Holmbom and Svanstedt for a related study from 2007 where there are no
temporal microscales), and in 2007 Flod\'{e}n and Olsson achieved homogenisation results for linear parabolic problems
involving pairs $\bigl( \eps, \{ \eps, \eps^r \} \bigr)$ of lists of scale functions (see \cite{FloOls07}); this was actually
the first time homogenisation was performed for problems with more than one temporal microscale. In 2009, Woukeng
studied non-linear non-monotone degenerated parabolic problems with the pair $\bigl( \eps, \{ \eps, \eps^k \} \bigr)$
of lists of spatial and temporal scale functions (see \cite{Wouk10}).

This paper deals with monotone parabolic problems with an arbitrary number of temporal microscales not necessarily
characterised by scale functions in the form of power functions but instead using the concept of (well-)separatedness
in spirit of \cite{AllBri96}. Furthermore---for simplicity---we only consider two spatial scales of which one is
microscopical, i.e., henceforth we fix $n = 1$.

Let $k \in \scope{m}$. Define $\jointly_{\wsep}^{m \sim k}$ to be the set of all pairs $\bigl( \eps ,
\{ \eps_j' \}_{j=1}^m \bigr)$ in $\jointly_{\wsep}^{1m}$ such that $\eps_k' \sim \eps$. (There is no loss
of generality to assume mere asymptotic equality rather than the ostensibly more general asymptotic equality modulo a
positive constant, i.e., $\eps_k' \sim C \eps$, $C \in \real$.) In other words, $\jointly_{\wsep}^{m \sim k}$
consists of pairs $\bigl( \eps , \{ \eps_j' \}_{j=1}^m \bigr)$ for which the temporal scale functions are separated and
the $k$-th temporal scale function coincides asymptotically with the spatial scale function. (This clearly explains the
convenient notation ``$\sim k$'' which could be read ``the spatial scale is asymptotically equal to the $k$-th temporal
scale''.)

Define the collection $\bigl\{ \jointly_{\wsep,i}^{m \sim k} \bigr\}_{i=1}^{1+2(m-k)}$ of $1 + 2(m - k)$ subsets of
$\jointly_{\wsep}^{m \sim k}$ by 
\begin{itemize}
	\item{} $\jointly_{\wsep,1}^{m \sim k} = \Bigl\{ \bigl( \eps , \{ \eps_j' \}_{j=1}^m \bigr)
				\in \jointly_{\wsep}^{m \sim k} \, : \, \tfrac{\eps^2}{\eps_m'} \rightarrow 0 \Bigr\}$,
	\item{} $\jointly_{\wsep,2}^{m \sim k} = \Bigl\{ \bigl( \eps , \{ \eps_j' \}_{j=1}^m \bigr)
				\in \jointly_{\wsep}^{m \sim k} \, : \, \eps_m' \sim \eps^2 \Bigr\}$,
	\item{} $\jointly_{\wsep,2+i-k}^{m \sim k} = \Bigl\{ \bigl( \eps , \{ \eps_j' \}_{j=1}^m \bigr)
				\in \jointly_{\wsep}^{m \sim k} \, : \, \tfrac{\eps_i'}{\eps^2} \rightarrow 0
					\textrm{ but } \tfrac{\eps_{i-1}'}{\eps^2} \rightarrow \infty \Bigr\}$ \quad ($i \in \scope{k+1,m}$),
	\item{} $\jointly_{\wsep,1+m+\iring-2k}^{m \sim k} = \Bigl\{ \bigl( \eps ,
				\{ \eps_j' \}_{j=1}^m \bigr) \in \jointly_{\wsep}^{m \sim k}
				\, : \, \eps_{\iring-1}' \sim \eps^2 \Bigr\}$ \qquad \quad \; \, \, ($\iring \in \scope{k+2,m}$). 
\end{itemize}
(Note that if $k  = m$, the collection of subsets of
$\jointly_{\wsep}^{m \sim m}$ reduces to merely
$\bigl\{ \jointly_{\wsep,1}^{m \sim m} \bigr\}$.) The subsets
$\jointly_{\wsep,1}^{m \sim k}$, $\jointly_{\wsep,2}^{m \sim k}$
and the collections of subsets $\bigl\{ \jointly_{\wsep,2+i-k}^{m \sim k} \bigr\}_{i=k+1}^m$
and $\bigl\{ \jointly_{\wsep,1+m+\iring-2k}^{m \sim k} \bigr\}_{\iring=k+2}^m$ of
$\jointly_{\wsep}^{m \sim k}$ correspond to slow temporal oscillations, slow resonance (i.e., ``slow''
self-similar case), rapid temporal oscillations and rapid resonance (i.e., ``rapid'' self-similar case),
respectively. 
\begin{theo}
\label{th:mutuallydisjoint}
The collection $\bigl\{ \jointly_{\wsep,i}^{m \sim k} \bigr\}_{i=1}^{1+2(m-k)}$
of $1 + 2(m - k)$ subsets of $\jointly_{\wsep}^{m \sim k}$ is mutually disjoint
for every $k \in \scope{m}$.
\end{theo}
\begin{proof}[\boldproof]
We must prove
\begin{equation*}
	\jointly_{\wsep,i}^{m \sim k}
		\cap \jointly_{\wsep,j}^{m \sim k} = \emptyset
\end{equation*}
for all $i,j \in \scope{1+2(m-k)}$ with $i \neq j$. That
\begin{equation*}
	\jointly_{\wsep,i}^{m \sim k}
		\cap \jointly_{\wsep,j}^{m \sim k} = \emptyset
\end{equation*}
for all $i,j \in \scope{2} \cup \scope{3+m-k,1+2(m-k)}$ with $i \neq j$, and that
\begin{equation*}
	\jointly_{\wsep,i}^{m \sim k} 
		\cap \jointly_{\wsep,j}^{m \sim k} = \emptyset
\end{equation*}
for all $i,j \in \scope{3,2+m-k}$ with $i \neq j$, are simple observations. It thus remains to show that
\begin{equation*}
	\jointly_{\wsep,i}^{m \sim k}
		\cap \jointly_{\wsep,j}^{m \sim k} = \emptyset
\end{equation*}
for all $i \in \scope{2} \cup \scope{3+m-k,1+2(m-k)}$ and all $j \in \scope{3,2+m-k}$. This is trivial for $k = m$
so it is understood that $k \in \scope{m-1}$ in the remainder of the proof.

Let $e \in \jointly_{\wsep,1}^{m \sim k}$ be arbitrary.
For this pair $e$ we have
\begin{equation*}
	\tfrac{\varepsilon^2}{\varepsilon_m'} \rightarrow 0,
\end{equation*}
which can be written
on the equivalent form
\begin{equation*}
	\tfrac{\varepsilon_m'}{\varepsilon^2} \rightarrow \infty,
\end{equation*}
or
\begin{equation*}
	\tfrac{\varepsilon_i'}{\varepsilon^2} \,
		\tfrac{\varepsilon_m'}{\varepsilon_i'} \rightarrow \infty
\end{equation*}
for every $i \in \scope{m}$. Furthermore,
\begin{equation*}
	\tfrac{\varepsilon_i'}{\varepsilon^2} \rightarrow \infty
\end{equation*}
since $\tfrac{\varepsilon_m'}{\varepsilon_i'}$ either tends to $0$ (if
$i \in \scope{m-1}$) or equals $1$ (if $i = m$). In particular this holds
for all $i \in \scope{k+1,m}$, and it is clear that $e \not\in
\jointly_{\wsep,2+i-k}^{m \sim k}$ for all $i \in \scope{k+1,m}$.
We have shown that
\begin{equation*}
	\jointly_{\wsep,1}^{m \sim k}
		\cap \jointly_{\wsep,j}^{m \sim k} = \emptyset
\end{equation*}
for all $j \in \scope{3,2+m-k}$.

Let $e \in \jointly_{\wsep,2}^{m \sim k}$ be arbitrary. Then we have
$\varepsilon_m' \sim \varepsilon^2$ for the chosen pair $e$ which gives
\begin{equation*}
	\tfrac{\varepsilon_i'}{\varepsilon^2} \sim \tfrac{\varepsilon_i'}{\varepsilon_m'},
\end{equation*}
$i \in \scope{m}$, which either tends to infinity (if $i \in \scope{m-1}$) or equals $1$
(if $i = m$). In particular this holds for all $i \in \scope{k+1,m}$. Thus, for every
$i \in \scope{k+1,m}$, $e \not\in \jointly_{\wsep,2+i-k}^{m \sim k}$, and
we have proven that
\begin{equation*}
	\jointly_{\wsep,2}^{m \sim k}
		\cap \jointly_{\wsep,j}^{m \sim k} = \emptyset
\end{equation*}
for all $j \in \scope{3,2+m-k}$.

Let $e \in \jointly_{\wsep,i}^{m \sim k}$, $i \in \scope{3+m-k,1+2(m-k)}$, be arbitrary.
The introduced pair $e$ satisfies $\varepsilon_{\iring-1}' \sim \varepsilon^2$, $\iring \in \scope{k+2,m}$, giving
\begin{multline*}
	\tfrac{\varepsilon_k'}{\varepsilon^2} \rightarrow \infty , \;
	\tfrac{\varepsilon_{k+1}'}{\varepsilon^2} \rightarrow \left\{
															\begin{aligned}
																1		&\textrm{ if } \iring = k+2 \\
																\infty	&\textrm{ if } \iring \in \scope{k+3,m}
															\end{aligned}
														\right. , \;
	\tfrac{\varepsilon_{k+2}'}{\varepsilon^2} \rightarrow \left\{
															\begin{aligned}
																0		&\textrm{ if } \iring = k+2 \\
																1		&\textrm{ if } \iring = k+3 \\
																\infty	&\textrm{ if } \iring \in \scope{k+4,m}
															\end{aligned}
														\right. , \\
	\ldots , \;
	\tfrac{\varepsilon_{m-2}'}{\varepsilon^2} \rightarrow \left\{
															\begin{aligned}
																0		&\textrm{ if } \iring \in \scope{m-2} \\
																1		&\textrm{ if } \iring = m-1 \\
																\infty	&\textrm{ if } \iring = m
															\end{aligned}
														\right. , \;
	\tfrac{\varepsilon_{m-1}'}{\varepsilon^2} \rightarrow \left\{
															\begin{aligned}
																0	&\textrm{ if } \iring \in \scope{m-1} k+2 \\
																1	&\textrm{ if } \iring = m
															\end{aligned}
														\right. , \;
	\tfrac{\varepsilon_m'}{\varepsilon^2} \rightarrow 0.
\end{multline*}
We see that $e \not\in \jointly_{\wsep,3}^{m \sim k}$. Indeed, to be in the subset requires $\tfrac{\varepsilon_{k+1}'}{\varepsilon^2} \rightarrow 0$
but $\tfrac{\varepsilon_k'}{\varepsilon^2} \rightarrow \infty$, which is clearly impossible.
We also see that $e \not\in \jointly_{\wsep,4}^{m \sim k}$, since being in the subset requires $\tfrac{\varepsilon_{k+2}'}{\varepsilon^2} \rightarrow 0$
but $\tfrac{\varepsilon_{k+1}'}{\varepsilon^2} \rightarrow \infty$; the former limit needs $\iring = k+2$ while the latter needs $\iring \in \scope{k+3,m}$.
We realise that $e \not\in \jointly_{\wsep,j}^{m \sim k}$ for all $j \in \scope{3,2+m-k}$. Hence,
\begin{equation*}
	\jointly_{\wsep,i}^{m \sim k}
		\cap \jointly_{\wsep,j}^{m \sim k} = \emptyset
\end{equation*}
for all $i \in \scope{3+m-k,1+2(m-k)}$ and all $j \in \scope{3,2+m-k}$. The mutual disjointness property has been verified.
\end{proof}
In the proposition below we will experience that the introduced collection of mutually disjoint subsets actually
forms a partition in the special but very important ``classical'' case of temporal scale functions expressed as
power functions. For this purpose, define the subset
\begin{align*}
	&\mathcal{P}^{m \sim k} = \Bigl\{ \bigl( \eps , \{ \eps_j' \}_{j=1}^m \bigr) \in
		\jointly_{\wsep}^{m \sim k} \, : \, \textrm{for every } \ell \in \scope{m} \\
	&\phantom{ \mathcal{P}^{m \sim k} = \Bigl\{ \bigl( \eps , \{ \eps_j' \}_{j=1}^m \bigr) \in }
		\textrm{ there exists a } p_{\ell} > 0 \textrm{ such that } \eps_{\ell}' = \eps^{p_{\ell}} \Bigr\}
\end{align*}
of $\jointly_{\wsep}^{m \sim k}$. We note that in the definition above for $\mathcal{P}^{m \sim k}$,
$p_k = 1$. Moreover, $0 < p_{\ell} < 1$ if $\ell \in \scope{k-1}$ (provided $k \in \scope{2,m}$) and $p_{\ell} > 1$
if $\ell \in \scope{k+1,m}$ (provided $k \in \scope{m-1}$). Furthermore, for each $i \in \scope{1+2(m-k)}$, define the
subsets
\begin{equation*}
	\mathcal{P}^{m \sim k}_i = \mathcal{P}^{m \sim k} \cap \jointly_{\wsep,i}^{m \sim k}
\end{equation*}
of $\mathcal{P}^{m \sim k}$. By Theorem~\ref{th:mutuallydisjoint} we already know that the collection
$\bigl\{ \mathcal{P}^{m \sim k}_i \bigr\}_{i=1}^{1+2(m-k)}$ is mutually disjoint. Below we will see that it actually
also covers all of $\mathcal{P}^{m \sim k}$.
\begin{prop}
\label{prop:collispart}
The collection $\bigl\{ \mathcal{P}^{m \sim k}_i \bigr\}_{i=1}^{1+2(m-k)}$ forms a partition of $\mathcal{P}^{m \sim k}$.
\end{prop}
\begin{proof}[\boldproof]
As already mentioned, the mutual disjointness property follows immediately from
Theorem~\ref{th:mutuallydisjoint}. It remains to show that
\begin{equation}
\label{eq:powerpartition}
	\mathcal{P}^{m \sim k} = \mathsmaller{\bigcup\limits}_{i=1}^{1+2(m-k)} \mathcal{P}^{m \sim k}_i,
\end{equation}
i.e., that the collection $\bigl\{ \mathcal{P}^{m \sim k}_i \bigr\}_{i=1}^{1+2(m-k)}$ of subsets covers all
of $\mathcal{P}^{m \sim k}$.

Suppose that there exists a pair
\begin{equation}
\label{eq:pairassumption}
	e \in \mathcal{P}^{m \sim k} \, \setminus \, \bigcup_{i=1}^{1+2(m-k)} \mathcal{P}^{m \sim k}_i,
\end{equation}
which means that we assume that $\bigl\{ \mathcal{P}^{m \sim k}_i \bigr\}_{i=1}^{1+2(m-k)}$ does not cover
all of $\mathcal{P}^{m \sim k}$. The introduced pair $e = \bigl( \eps, \{ \eps^{p_j} \}_{j=1}^m \bigr)$ must by definition satisfy
\begin{align}
	\frac{\eps^2}{\eps^{p_m}}		&\not\rightarrow 0
		&	&\textrm{since } e \not\in \mathcal{P}^{m \sim k}_1, \label{eq:powcond1} \\
	\eps^{p_m}						&\not\sim \eps^2
		&	&\textrm{since } e \not\in \mathcal{P}^{m \sim k}_2, \label{eq:powcond2} \\
	\frac{\eps^{p_i}}{\eps^2}	&\not\rightarrow 0 \textrm{ or } \frac{\eps^{p_{i-1}}}{\eps^2} \not\rightarrow \infty
		\quad \forall i \in \scope{k+1,m}
		&	&\textrm{since } e \not\in \bigcup_{\ell = 3}^{2+m-k} \mathcal{P}^{m \sim k}_{\ell}, \label{eq:powcond3}
\inter{and}
	\frac{\eps^{p_{\iring-1}}}{\eps^2} 	&\not\rightarrow 1 \qquad \forall \iring \in \scope{k+2,m}
		&	&\textrm{since } e \not\in \bigcup_{\ellring = 3+m-k}^{1+2(m-k)} \mathcal{P}^{m \sim k}_{\ellring}. \label{eq:powcond4}
\end{align}
The conditions \eqref{eq:powcond2} and \eqref{eq:powcond4} may be written
\begin{equation}
\label{eq:equiv2and4}
	( \, p_{k+1} \neq 2 \, ) \wedge ( \, p_{k+2} \neq 2 \, ) \wedge \ldots \wedge
		( \, p_{m-1} \neq 2 \, ) \wedge ( \, p_m \neq 2 \ ),
\end{equation}
and \eqref{eq:powcond3} can be expressed as
\begin{multline}
\label{eq:equiv3}
	\bigl( \, ( \, p_k \geqslant 2 \, ) \vee ( \, p_{k+1} \leqslant 2 \, ) \, \bigr)
		\wedge
	\bigl( \, ( \, p_{k+1} \geqslant 2 \, ) \vee ( \, p_{k+2} \leqslant 2 \, ) \, \bigr) \\
		\wedge \ldots \wedge
	\bigl( \, ( \, p_{m-2} \geqslant 2 \, ) \vee ( \, p_{m-1} \leqslant 2 \, ) \, \bigr)
		\wedge
	\bigl( \, ( \, p_{m-1} \geqslant 2 \, ) \vee ( \, p_m \leqslant 2 \, ) \, \bigr),
\end{multline}
where we employ the logic symbols  $\wedge$ `and' (i.e., logical conjunction) and
$\vee$ `or' (i.e., logical disjunction) for clarity.

We begin by noticing that $p_k = 1$ by definition, so \eqref{eq:equiv3} implies that $p_{k+1} \leqslant 2$.
This together with $p_{k+1} \neq 2 $ from \eqref{eq:equiv2and4} yields $p_{k+1} < 2$. Hence, using \eqref{eq:equiv3}
again and we conclude that $p_{k+1} \leqslant 2$. Consequently, \eqref{eq:equiv2and4} implies $p_{k+1} < 2$.
Continuing, we end up with $p_m < 2$. But this contradicts \eqref{eq:powcond1} which states that $p_m \geqslant 2$.
Thus, no pair $e$ fulfilling \eqref{eq:pairassumption} can exist so we indeed have \eqref{eq:powerpartition},
and the proof is complete.
\end{proof}
\begin{exam}
In \cite{FloOls07} one considers pairs of the type $\bigl( \eps, \{ \eps, \eps^r \} \bigr)$,
$r \in \real_+ \setminus \{ 1 \}$, in the context of linear parabolic problems. Define the mutually disjoint sets
\begin{align*}
	\mathcal{R}^-	&= \bigl\{ \bigl( \eps, \{ \eps^r, \eps \} \bigr)
		\in \mathcal{P}^{2 \sim 2} \, : \, 0 < r < 1 \bigr\}, \\
\inter{and}
	\mathcal{R}^+	&= \bigl\{ \bigl( \eps, \{ \eps, \eps^r \} \bigr)
		\in \mathcal{P}^{2 \sim 1} \, : \, r > 1 \bigr\},
\end{align*}
and let $\mathcal{R} = \mathcal{R}^- \cup \mathcal{R}^+$. Introduce the subsets
\begin{align*}
	\mathcal{R}^-_1	&= \mathcal{P}^{2 \sim 2}_1 \cap \mathcal{R}^- = \bigl\{ \bigl( \eps, \{ \eps^r, \eps \} \bigr)
			\in \mathcal{P}^{2 \sim 2} \, : \, 0 < r < 1 \bigr\} = \mathcal{R}^- \\
\inter{and}
	\mathcal{R}^-_2	&= \mathcal{P}^{2 \sim 2}_2 \cap \mathcal{R}^- = \emptyset
\end{align*}
of $\mathcal{R}^-$, and the subsets
\begin{align*}
	\mathcal{R}^+_1	&= \mathcal{P}^{2 \sim 1}_1 \cap \mathcal{R}^+
		= \bigl\{ \bigl( \eps, \{ \eps, \eps^r \} \bigr) \in \mathcal{P}^{2 \sim 1} \, : \, 1 < r < 2 \bigr\}, \\
	\mathcal{R}^+_2	&= \mathcal{P}^{2 \sim 1}_2 \cap \mathcal{R}^+
		= \bigl\{ \bigl( \eps, \{ \eps, \eps^r \} \bigr) \in \mathcal{P}^{2 \sim 1} \, : \, r = 2 \bigr\}
\inter{and}
	\mathcal{R}^+_3	&= \mathcal{P}^{2 \sim 1}_3 \cap \mathcal{R}^+
		= \bigl\{ \bigl( \eps, \{ \eps, \eps^r \} \bigr) \in \mathcal{P}^{2 \sim 1} \, : \, r > 2 \bigr\}
\end{align*}
of $\mathcal{R}^+$. By Proposition~\ref{prop:collispart}, $\mathcal{R}^-$ and $\mathcal{R}^+$
are partitioned by the collections $\{ \mathcal{R}^-_1,\mathcal{R}^-_2 \}$ and $\{ \mathcal{R}^+_1,\mathcal{R}^+_2,
\mathcal{R}^+_3 \}$, respectively. Thus, according to the developed theory, $\mathcal{R}$ should be partitioned by
the collection
\begin{equation*}
	\{ \mathcal{R}^-_1,\mathcal{R}^-_2 \, , \, \mathcal{R}^+_1,\mathcal{R}^+_2,\mathcal{R}^+_3 \},
\end{equation*}
which is verified by looking at the explicit expressions for the subsets derived above. Defining $\mathcal{R}_1 =
\mathcal{R}^-_1 \cup \mathcal{R}^+_1$, the partitioning collection of subsets
\begin{equation*}
	\{ \mathcal{R}_1, \mathcal{R}^+_2, \mathcal{R}^+_3 \}
\end{equation*}
of $\mathcal{R}$ is seen to correspond to the cases $0 < r < 2$ with $r \neq 1$, $r = 2$ and $r > 2$, respectively.
This is exactly the partition obtained in the homogenisation result of Theorem~10 in \cite{FloOls07} leading to three
distinct systems of local problems for $u_1$ corresponding to the mentioned distinct cases for
$r \in \real_+ \setminus \{ 1 \}$.
\end{exam}
Let $S = (0,1)$ and define $H_{\#}^1 (S; V,V')$, $V$ being any Banach space with topological dual $V'$,
as the space of functions $u$ satisfying $u \in L_{\#}^2(S;V)$ and $\tfrac{\rmd}{\rmd s} u
\in L_{\#}^2(S;V')$. In order to prove Theorem~\ref{th:maintheorem}---our first homogenisation
result---we first need the lemmas below.
\begin{lemma}
\label{lem:tensprodspdenseinsobsp}
The tensor product space $\bigl( \cont_{\#}^{\infty}(Y)/\real \bigr) \otimes \cont_{\#}^{\infty}(S)$
is dense in $H_{\#}^1 (S; \wrum,\wrum')$.
\end{lemma}
\begin{proof}[\boldproof]
This is just Proposition~4.6 in \cite{NguWou07} in which $\mathcal{E}$ and $\mathcal{V}$ correspond to
$\bigl( \cont_{\#}^{\infty}(Y)/\real \bigr) \otimes \cont_{\#}^{\infty}(S)$ and $H_{\#}^1 (S; \wrum,\wrum')$,
respectively, of the present paper.
\end{proof}
\begin{lemma}
\label{lem:vanishintegrsum}
Suppose that $u, v \in H_{\#}^1 (S; \wrum,\wrum')$. Then
\begin{equation*}
	\langle \d_s u, v \rangle_{L_{\#}^2(S;\wrum'),L_{\#}^2(S;\wrum)}
		+ \langle \d_s v, u \rangle_{L_{\#}^2(S;\wrum'),L_{\#}^2(S;\wrum)} = 0
\end{equation*}
holds. In particular,
\begin{equation*}
	\langle \d_s u, u \rangle_{L_{\#}^2(S;\wrum'),L_{\#}^2(S;\wrum)} = 0.
\end{equation*}
\end{lemma}
\begin{proof}[\boldproof]
This follows immediately from Corollary~4.1 in \cite{NguWou07}.
\end{proof}
Theorem~\ref{th:maintheorem} below is our first homogenisation result.
\begin{theo}
\label{th:maintheorem}
Let $k \in \scope{m}$. Suppose that the pair $e = \bigl( \eps , \{ \eps_j' \}_{j=1}^m \bigr)$ of lists of
spatial and temporal scale functions belongs to $\bigcup_{i=1}^{1+2(m-k)} \jointly_{\wsep,i}^{m \sim k}$. Let $\{ u_{\eps} \}$
be the sequence of weak solutions in $H^1 \bigl( \ZeToT; H_0^1(\Omega), H^{-1}(\Omega) \bigr)$ to the evolution
problem \eqref{eq:monoparaprob} with $a \, : \, \bar{\Omega}_T \times \real^{N+m} \times \real^N
\rightarrow \real^N$ satisfying the structure conditions \emph{(B}$_1$\emph{)}--\emph{(B}$_5$\emph{)}. Then
\begin{align*}
	u_{\eps}		&\rightarrow	u	&	&\textrm{in } L^2(\Omega_T), \\
	u_{\eps}		&\rightharpoonup  u	&	&\textrm{in } L^2 \bigl( \ZeToT; H_0^1(\Omega) \bigr), \\
\inter{and}
 	\grad u_{\eps}	&\scaleconv{(2,m+1)} \grad u + \grad_{y} u_1, & &
\end{align*}
where $u \in H^1 \bigl( \ZeToT; H_0^1(\Omega), H^{-1}(\Omega) \bigr)$ and $u_1 \in
L^2 \bigl( \Omega_T \times S^m; \wrum \bigr)$. Here $u$ is the unique weak solution to the homogenised
problem \eqref{eq:defoflim} with the homogenised flux $b \, : \, \bar{\Omega}_T \times \real^N
\rightarrow \real^N$ given by
\begin{equation}
	\label{eq:hmpliminth}
	b(x,t;\grad u) = \smint_{\sety_{1m}} a(x,t, y,\vecs_m; \grad u + \grad_{y} u_1) \, \rmd \vecs_m \rmd y.
\end{equation}
Moreover, we have the following characterisation of $u_1$:

$\bullet$ If $e \in \jointly_{\wsep,1}^{m \sim k}$ then the function $u_1$ is the unique weak solution to the local problem
\begin{equation*}
	- \grad_y \cdot a(x,t,y,\vecs_m; \grad u + \grad_y u_1) = 0.
\end{equation*}

$\bullet$ If $e \in \jointly_{\wsep,2}^{m \sim k}$, assuming
$u_1 \in L^2 \bigl( \Omega_T \times S^{m-1}; H_{\#}^1(S_m;\wrum,\wrum') \bigr)$,
then the function $u_1$ is the unique weak solution to the system of local problems
\begin{equation*}
	\d_{s_m} u_1(x,t,y,\vecs_m) - \grad_y \cdot a(x,t,y,\vecs_m; \grad u + \grad_y u_1) = 0.
\end{equation*}

$\bullet$ If $e \in \jointly_{\wsep,2+\bar{\ell}-k}^{m \sim k}$ for some $\bar{\ell} \in \scope{k+1,m}$,
provided $k \in \scope{m-1}$, then the function $u_1$ is the unique weak solution to the system of local problems
\begin{equation*}
 	\left\{
 		\begin{aligned}
 			- \grad_y \cdot \smint_{S^{\scope{\bar{\ell},m}}}
 				a(x,t,y,\vecs_m; \grad u + \grad_y u_1) \, \rmd \vecs_{\scope{\bar{\ell},m}}		&= 0, \\
				\d_{s_i} u_1(x,t,y,\vecs_m)	&= 0 \qquad (i \in \scope{\bar{\ell},m}).
 		\end{aligned}
 	\right.
\end{equation*}

$\bullet$ If $e \in \jointly_{\wsep,1+m+\ellring-2k}^{m \sim k}$ for
some $\ellring \in \scope{k+2,m}$, provided $k \in \scope{m-2}$ and assuming $u_1 \in
L^2 \bigl( \Omega_T \times S^{\ellring-2} \times S^{\scope{\ellring,m}};H_{\#}^1(S_{\ellring-1};\mathcal{W},\mathcal{W}') \bigr)$,
then the function $u_1$ is the unique weak solution to the system of local problems
\begin{equation*}
 	\left\{
 		\begin{aligned}
			\d_{s_{\ellring-1}} u_1(x,t,y,\vecs_m)
				- \grad_y \! \cdot \! \smint_{S^{\scope{\ellring,m}}}
					a(x,t,y,\vecs_m; \grad u + \grad_y u_1)
						\, \rmd \vecs_{\scope{\ellring,m}}	&= 0, \\
			\d_{s_i} u_1(x,t,y,\vecs_m)								&= 0
			\qquad (i \in \scope{\ellring,m}).
 		\end{aligned}
 	\right.
\end{equation*}
\end{theo}
\begin{proof}[\boldproof]
Since $a$ fulfils (B$_1$)--(B$_5$) we can use Proposition~\ref{prop:twounifbound} for the sequence
$\{ u_{\eps} \}$ of weak solutions; we have ensured uniform boundedness in $H^1 \bigl( \ZeToT; H_0^1(\Omega),
H^{-1}(\Omega) \bigr)$, i.e., \eqref{eq:firstunifbound} holds. We can then employ Theorem~\ref{th:gradchar}
(with $n = 1$) obtaining, up to a subsequence,
\begin{align}
	u_{\eps}		&\rightarrow	u
		& &\textrm{in } L^2(\Omega_T), \label{eq:strongconvofueps}  \\
	u_{\eps}		&\rightharpoonup  u
		& &\textrm{in } L^2 \bigl( \ZeToT; H_0^1(\Omega) \bigr), \label{eq:weakconvofueps} \\
\inter{and}
 	\grad u_{\eps}	&\scaleconv{(2,m+1)} \grad u + \grad_{y} u_1, & & \label{eq:1mconvofueps}
\end{align}
where $u \in H^1 \bigl( \ZeToT; H_0^1(\Omega), H^{-1}(\Omega) \bigr)$ and $u_1 \in
L^2 \bigl( \Omega_T \times S^m; \wrum \bigr)$. Consider the sequence $\{ a_{\eps} \}$ defined according to
\begin{align*}
	a_{\eps} (x,t)	&= a^{\eps} (x,t; \grad u_{\eps}) \\
							&= a(x,t, \tfrac{x}{\eps},\vect_m^{\eps} ;\grad u_{\eps})
								\qquad \qquad ((x,t) \in \Omega_T).
\end{align*}
We have that $\{ a_{\eps} \}$ is uniformly bounded in $L^2(\Omega_T)^N$. Indeed, using \eqref{eq:anineqfora},
the triangle inequality and \eqref{eq:firstunifbound} we get
\begin{align*}
	\| a_{\eps} \|_{L^2(\Omega_T)^N}^2
			&= \smint_{\Omega_T} \bigl| a(x,t, \tfrac{x}{\eps},\vect_m^{\eps} ;\grad u_{\eps}) \bigr|^2 \, \rmd x \rmd t \\
			&< C_1^2 \smint_{\Omega_T} \bigl( 1 + \bigl| \grad u_{\eps}(x,t) \bigr| \bigr)^2 \, \rmd x \rmd t \\
			&= C_1^2 \bigl\| 1 + | \grad u_{\eps} | \bigr\|_{L^2(\Omega_T)}^2 \\
			&\leqslant C_1^2 \Bigl( \| 1 \|_{L^2(\Omega_T)}
				+ \| u_{\eps} \|_{L^2 \bigl( \ZeToT; H_0^1(\Omega) \bigr)} \Bigr)^2 \\
			&\leqslant C_1^2 \Bigl( \bigl( T | \Omega | \bigr)^{\frac{1}{2}}
				+  \| u_{\eps} \|_{H^1 \bigl( \ZeToT; H_0^1(\Omega), H^{-1}(\Omega) \bigr)} \Bigr)^2 \\
			&\leqslant C_1^2 \bigl( ( T | \Omega | )^{\frac{1}{2}} + C \bigr)^2.
\end{align*}
By Theorem~\ref{th:nm-scaleconvres} (with $n = 1$) we then know that, up to a subsequence,
\begin{equation}
	\label{eq:1m-limofsigmaeps}
	a_{\eps} \scaleconv{(2,m+1)} a_0
\end{equation}
for some $a_0 \in L^2 (\Omega_T \times \sety_{1m})^N$.

Recall the weak form \eqref{eq:weakformofevprob} (with $n = 1$) of the evolution problem, i.e.,
\begin{equation}
		\label{eq:weakevolprob1}
		\bigl\langle \tfrac{\d}{\d t} u_{\eps} , \psi \bigr\rangle_{X',X}
			+ \smint_{\Omega_T} a_{\eps}(x,t) \cdot \grad \psi (x,t) \, \rmd x \rmd t
				= \smint_{\Omega_T} f(x,t) \, \psi (x,t) \, \rmd x \rmd t
\end{equation}
for every $\psi \in L^2 \bigl( \ZeToT; H_0^1(\Omega) \bigr)$.

Choose an arbitrary $\psi \in H_0^1(\Omega) \odot \distr(\ZeToT)$. Then we can shift the weak temporal derivative
$\tfrac{\d}{\d t}$ in \eqref{eq:weakevolprob1} from acting on $u_{\eps}$ to acting on $\psi$ instead, i.e.,
\begin{equation}
		\label{eq:weakevolprob2}
		\smint_{\Omega_T} \bigl( - u_{\eps} (x,t) \, \tfrac{\d}{\d t} \psi (x,t)
			+ a_{\eps}(x,t) \cdot \grad \psi (x,t) \bigr) \, \rmd x \rmd t
				= \smint_{\Omega_T} f(x,t) \, \psi (x,t) \, \rmd x \rmd t.
\end{equation}
Passing to the limit---using \eqref{eq:weakconvofueps} and \eqref{eq:1m-limofsigmaeps} on the first and second
terms on the left-hand side, respectively---we obtain, up to a subsequence,
\begin{multline*}
		\smint_{\Omega_T} \smint_{\sety_{1m}} \bigl( - u (x,t) \, \tfrac{\d}{\d t} \psi (x,t)
			+ a_0 (x,t,y,\vecs_m) \cdot \grad \psi (x,t) \bigr) \, \rmd \vecs_m \rmd y \rmd x \rmd t \\
				= \smint_{\Omega_T} f(x,t) \, \psi (x,t) \, \rmd x \rmd t,
\end{multline*}
or, in other words,
\begin{multline}
	\label{eq:weakevolproblim}
		\smint_{\Omega_T} \Bigl( - u (x,t) \, \tfrac{\d}{\d t} \psi (x,t)
			+ \smint_{\sety_{1m}} a_0 (x,t,y,\vecs_m) \, \rmd \vecs_m \rmd y \cdot \grad \psi (x,t) \Bigr) \, \rmd x \rmd t \\
				= \smint_{\Omega_T} f(x,t)  \, \psi (x,t) \, \rmd x \rmd t.
\end{multline}
Let again $\tfrac{\d}{\d t}$ act on $u$. By density, the obtained equality
\begin{align}
	&\bigl\langle \tfrac{\d}{\d t} u , \psi \bigr\rangle_{X',X}
		+ \smint_{\Omega_T} \smint_{\sety_{1m}} a_0 (x,t,y,\vecs_m)
			\, \rmd \vecs_m \rmd y \cdot \grad \psi (x,t) \, \rmd x \rmd t \nonumber \\
	&\phantom{ \bigl\langle \tfrac{\d}{\d t} u , \psi \bigr\rangle_{X',X}
				+ \smint_{\Omega_T} \smint_{\sety_{1m}} a_0 (x,t,y,\vecs_m) \, \rmd \vecs_m \rmd y \cdot \grad }
		= \smint_{\Omega_T} f(x,t)  \, \psi (x,t) \, \rmd x \rmd t \label{eq:weakevolproblim2}
\end{align}
holds for any $\psi \in L^2 \bigl( \ZeToT; H_0^1(\Omega) \bigr)$. We have obtained the weak form of the homogenised
evolution problem \eqref{eq:defoflim} with the limit flux given by
\begin{equation*}
	b(x,t;\grad u) = \smint_{\sety_{1m}} a_0 (x,t,y,\vecs_m) \, \rmd \vecs_m \rmd y.
\end{equation*}

What remains is to find the local problems for $u_1$ and to give the limit $a_0$ in terms of $a$. We will first extract
the pre-local-problems, i.e., the problems expressed in terms of $a_0$ which become the local problems once $a_0$
is given in terms of $a$. Introduce
$\omega_{\ell} \in \distr(\Omega) \! \odot \! \distr(\ZeToT) \, \odot \, \bigl( \cont_{\#}^{\infty}(Y) / \real \bigr)
\! \odot \! \bigl( \prod_{i=1}^{\ell} \cont_{\#}^{\infty}(S_i) \bigr)$, $\ell \in \scope{m}$. For each
$\ell \in \scope{m}$ we define the sequence $\{ \omega_{\ell}^{\eps} \}$ in the conventional manner. Let
$\{ r_{\eps} \}$ be a sequence of positive numbers such that $r_{\eps} \rightarrow 0$. We will now study sequences
of test functions $\{ \psi^{\eps} \}$ in \eqref{eq:weakevolprob2} such that
\begin{equation*}
	\psi_{\ell}^{\eps}(x,t) = r_{\eps} \omega_{\ell}^{\eps}(x,t) \qquad \qquad ((x,t) \in \Omega_T)
\end{equation*}
with appropriate choices of $\{ r_{\eps} \}$ and $\ell$ in order to extract the pre-local-problems. We note here that
\begin{align*}
	\grad \psi_{\ell}^{\eps}
		&= r_{\eps} \bigl( \grad_x + \tfrac{1}{\eps} \grad_y \bigr) \omega_{\ell}^{\eps}
\inter{and}
	\tfrac{\d}{\d t} \psi_{\ell}^{\eps}
		&= r_{\eps} \Bigl( \d_t + \ssum_{i=1}^{\ell} \tfrac{1}{\eps_i'} \d_{s_i} \Bigr) \omega_{\ell}^{\eps}.
\end{align*}
For the sequence $\{ \psi_{\ell}^{\eps} \}$, $\ell \in \scope{m}$, of test functions given above,
\eqref{eq:weakevolprob2} becomes
\begin{align*}
	\smint_{\Omega_T} \Bigl[ -u_{\eps}(x,t) r_{\eps} &\Bigl( \d_t
		+ \ssum_{i=1}^{\ell} \tfrac{1}{\eps_i'} \d_{s_i} \Bigr) \omega_{\ell}^{\eps}(x,t) \\
			&+ a_{\eps} (x,t) \cdot r_{\eps} \bigl( \grad_x
				+ \tfrac{1}{\eps} \grad_y \bigr) \omega_{\ell}^{\eps} (x,t) \Bigr] \, \rmd x \rmd t
					= \smint_{\Omega_T} f(x,t) r_{\eps} \omega_{\ell}^{\eps}(x,t) \, \rmd x \rmd t.
\end{align*}
The right-hand side and the $\d_t$ and $\grad_x$ terms in the left-hand side clearly vanish in the limit,
and what is left is
\begin{equation}
	\label{eq:eqwithsubtf}
	\lim_{\eps \rightarrow 0} \smint_{\Omega_T} \Bigl( -u_{\eps}(x,t) \ssum_{i=1}^{\ell}
		\tfrac{r_{\eps}}{\eps_i'} \d_{s_i} \omega_{\ell}^{\eps}(x,t) + a_{\eps} (x,t) \cdot
			\tfrac{r_{\eps}}{\eps_k'} \grad_y \omega_{\ell}^{\eps} (x,t) \Bigr) \, \rmd x \rmd t = 0
\end{equation}
recalling that $\eps_k' = \eps$.

Suppose that the real sequence $\{ \tfrac{r_{\eps}}{\eps_{\ell}'} \}$ is bounded,
then the limit equation becomes
\begin{equation}
	\label{eq:whatisleft}
	\lim_{\eps \rightarrow 0} \smint_{\Omega_T} \Bigl( -u_{\eps}(x,t) \tfrac{r_{\eps}}{\eps_{\ell}'}
		\d_{s_{\ell}} \omega_{\ell}^{\eps}(x,t) + a_{\eps} (x,t) \cdot
			\tfrac{r_{\eps}}{\eps_k'} \grad_y \omega_{\ell}^{\eps} (x,t) \Bigr) \, \rmd x \rmd t = 0.
\end{equation}
Choose $r_{\eps} = \eps_k'$, which implies that $\{ \tfrac{r_{\eps}}{\eps_{\ell}'} \}
= \{ \tfrac{\eps_k'}{\eps_{\ell}'} \}$ is bounded for $\ell \in \scope{k}$. Then \eqref{eq:whatisleft} becomes
\begin{equation}
	\label{eq:whatisleftnext}
	\lim_{\eps \rightarrow 0} \smint_{\Omega_T} \Bigl( -u_{\eps}(x,t) \tfrac{\eps_k'}{\eps_{\ell}'}
		\d_{s_{\ell}} \omega_{\ell}^{\eps}(x,t) + a_{\eps} (x,t) \cdot
			\grad_y \omega_{\ell}^{\eps} (x,t) \Bigr) \, \rmd x \rmd t = 0.
\end{equation}
If $\ell \in \scope{k-1}$ (provided $k \in \scope{2,m}$) the first term tends to $0$, and we get in this case
\begin{equation*}
	\lim_{\eps \rightarrow 0} \smint_{\Omega_T} a_{\eps} (x,t)
		\cdot \grad_y \omega_{\ell}^{\eps} (x,t) \, \rmd x \rmd t = 0,
\end{equation*}
which after taking the limit can be written
\begin{equation}
	\label{eq:firsttermgonelim}
	\smint_{\Omega_T} \smint_{\sety_{1m}} a_0 (x,t,y,\vecs_m)
		\cdot \grad_y \omega_{\ell} (x,t,y,\vecs_{\ell}) \, \rmd \vecs_m \rmd y \rmd x \rmd t = 0,
\end{equation}
i.e.,
\begin{equation*}
	\smint_{\Omega_T} \smint_{\sety_{1 \ell}} \smint_{S^{\scope{\ell+1,m}}}
		a_0 (x,t,y,\vecs_m) \rmd \vecs_{\scope{\ell+1,m}}
			\cdot \grad_y \omega_{\ell} (x,t,y,\vecs_{\ell}) \, \rmd \vecs_{\ell} \rmd y \rmd x \rmd t = 0.
\end{equation*}
Suppose $v_1 \in \cont_{\#}^{\infty}(Y) / \real$ is the factor of $\omega_{\ell}$ with respect to the $y$
variable. Then, by the Variational Lemma,
\begin{equation}
	\label{eq:case1astrong}
	\smint_Y \smint_{S^{\scope{\ell+1,m}}} a_0 (x,t,y,\vecs_m) \rmd \vecs_{\scope{\ell+1,m}}
		\, \cdot \, \grad_y v_1 (y) \, \rmd y = 0
\end{equation}
a.e.~on $\Omega_T \times S^{\ell}$. If $\ell = k$ the limit equation \eqref{eq:whatisleftnext} instead reduces to
\begin{equation*}
	\lim_{\eps \rightarrow 0} \smint_{\Omega_T} \bigl( -u_{\eps}(x,t) \d_{s_k} \omega_k^{\eps}(x,t)
		+ a_{\eps} (x,t) \cdot \grad_y \omega_k^{\eps} (x,t) \bigr) \, \rmd x \rmd t = 0,
\end{equation*}
which in the limit becomes
\begin{equation*}
	\smint_{\Omega_T} \smint_{\sety_{1m}} \bigl( -u(x,t) \d_{s_k} \omega_k (x,t,y,\vecs_k)
		+ a_0 (x,t,y,\vecs_m) \cdot \grad_y \omega_k (x,t,y,\vecs_{k}) \bigr) \, \rmd \vecs_m \rmd y \rmd x \rmd t = 0.
\end{equation*}
The first term gives no contribution since $\omega_k$ is $S_k$-periodic in the $s_k$ variable. Progressing like
in the case $\ell \in \scope{k-1}$ we finally arrive at \eqref{eq:case1astrong} which now also includes $\ell = k$,
i.e., \eqref{eq:case1astrong} holds for all $\ell \in \scope{k}$. But it is clear that \eqref{eq:case1astrong}
holding for $\ell = k$ implies that it holds also for any $\ell \in \scope{k-1}$ (provided $k \in \scope{2,m}$).
Thus, we only have to consider \eqref{eq:case1astrong} for $\ell = k$, i.e., we have so far obtained
\begin{equation}
	\label{eq:case1astrongonlyk}
	\smint_Y \smint_{S^{\scope{k+1,m}}} a_0 (x,t,y,\vecs_m) \, \rmd \vecs_{\scope{k+1,m}} \cdot \grad_y v_1(y) \, \rmd y = 0.
\end{equation}
It should be emphasised here that this equation is always true for $\jointly_{\wsep}^{m \sim k}$ and
is not confined to any particular subset $\jointly_{\wsep,j}^{m \sim k}$, $j \in \scope{1+2(m-k)}$.

If we study the limit equation \eqref{eq:eqwithsubtf} extracting a factor $\tfrac{1}{\eps}$ in the first term we obtain
\begin{equation*}
	\lim_{\eps \rightarrow 0} \smint_{\Omega_T} \Bigl( - \tfrac{1}{\eps} u_{\eps}(x,t) \ssum_{i=1}^{\ell}
		\tfrac{r_{\eps}\eps_k'}{\eps_i'} \d_{s_i} \omega_{\ell}^{\eps}(x,t) + a_{\eps} (x,t) \cdot
			\tfrac{r_{\eps}}{\eps_k'} \grad_y \omega_{\ell}^{\eps}(x,t) \Bigr) \, \rmd x \rmd t = 0,
\end{equation*}
where we have recalled $\eps_k' = \eps$. Suppose that $\{ \tfrac{r_{\eps} \eps_k'}{\eps_{\ell}'} \}$ is bounded
(in $\real$), it is then clear that the limit equation above reduces to
\begin{equation}
	\label{eq:whatisleftdiff}
	\lim_{\eps \rightarrow 0} \smint_{\Omega_T} \Bigl( -\tfrac{1}{\eps} u_{\eps}(x,t)
		\tfrac{r_{\eps} \eps_k'}{\eps_{\ell}'} \d_{s_{\ell}} \omega_{\ell}^{\eps}(x,t)
			+ a_{\eps} (x,t) \cdot \tfrac{r_{\eps}}{\eps_k'} \grad_y \omega_{\ell}^{\eps}(x,t) \Bigr) \, \rmd x \rmd t = 0.
\end{equation}

$\bullet$ Suppose $\bigl( \eps, \{ \eps_j' \}_{j=1}^m \bigr) \in \jointly_{\wsep,1}^{m \sim k}$.
By definition this means that $\bigl( \eps, \{ \eps_j' \} \bigr) \in \jointly_{\wsep}^{m \sim k}$
and $\tfrac{\eps_k'^2}{\eps_m'} \rightarrow 0$. Consider first $\eps_m' \sim \eps_k'$, i.e., $k = m$. We have
already extracted \eqref{eq:case1astrongonlyk} which in this case, $k = m$, is merely
\begin{equation}
	\label{eq:locprobfirst}
	\smint_Y a_0 (x,t,y,\vecs_m) \cdot \grad_y v_1(y) \, \rmd y = 0,
\end{equation}
which is the pre-local-problem.

Consider now the situation $\eps_m' \not\sim \eps_k'$, i.e., $k \in \scope{m-1}$ requiring $m > 1$.
We first note that we have already extracted \eqref{eq:case1astrongonlyk} which at this point carries
at least one integral (over $S_m$). We want to employ \eqref{eq:whatisleftdiff} for $\ell \in \scope{k+1,m}$.
Choose $r_{\eps} = \eps_k'$, and we get that
\begin{align*}
	\tfrac{r_{\eps} \eps_k'}{\eps_{\ell}'}	&= \tfrac{\eps_k'^2}{\eps_{\ell}'}
			=\tfrac{\eps_k'^2}{\eps_m'} \tfrac{\eps_m'}{\eps_{\ell}'} \\
			&\rightarrow 0,
\end{align*}
so $\bigl\{ \tfrac{r_{\eps} \eps_k'}{\eps_{\ell}'} \bigr\}$ is indeed bounded (we even have a vanishing limit).
We can now use \eqref{eq:whatisleftdiff} which yields
\begin{equation*}
	\lim_{\eps \rightarrow 0} \smint_{\Omega_T} \Bigl( -\tfrac{1}{\eps} u_{\eps}(x,t)
		\tfrac{\eps_k'^2}{\eps_{\ell}'} \d_{s_{\ell}} \omega_{\ell}^{\eps}(x,t)
			+ a_{\eps} (x,t) \cdot \grad_y \omega_{\ell}^{\eps}(x,t) \Bigr) \, \rmd x \rmd t = 0,
\end{equation*}
which in the limit becomes \eqref{eq:firsttermgonelim}; this can be realised by utilising Theorem~\ref{th:nmtypelimit},
considering the final remark in Remark~\ref{re:nmtypelimit} and using $\tfrac{\eps_k'^2}{\eps_{\ell}'} \rightarrow 0$ such
that the contribution from the first term vanishes in the limit. Hence, we have again \eqref{eq:case1astrong} but for
$\ell \in \scope{k+1,m}$. Apparently we end up at the pre-local-problem \eqref{eq:locprobfirst} again since \eqref{eq:case1astrong}
in the case $\ell = m$ implies that \eqref{eq:case1astrong} holds automatically for any $\ell \in \scope{m-1}$.

$\bullet$ Suppose $\bigl( \eps, \{ \eps_j' \}_{j=1}^m \bigr) \in \jointly_{\wsep,2}^{m \sim k}$. By definition
this means that $\bigl( \eps, \{ \eps_j' \} \bigr) \in \jointly_{\wsep}^{m \sim k}$ and $\eps_m' \sim \eps_k'^2$.
Let $\ell = m$ in \eqref{eq:whatisleftdiff}. Choose $r_{\eps} = \eps_k'$ again, giving
\begin{align*}
	\tfrac{r_{\eps} \eps_k'}{\eps_{\ell}'} = \tfrac{\eps_k'^2}{\eps_m'} \sim 1,
\end{align*}
so $\bigl\{ \tfrac{r_{\eps} \eps_k'}{\eps_{\ell}'} \bigr\}$ is bounded. The equation \eqref{eq:whatisleftdiff}
then becomes
\begin{equation*}
	\lim_{\eps \rightarrow 0} \smint_{\Omega_T} \Bigl( -\tfrac{1}{\eps} u_{\eps}(x,t)
		\tfrac{\eps_k'^2}{\eps_m'} \d_{s_m} \omega_m^{\eps}(x,t) + a_{\eps} (x,t)
			\cdot \grad_y \omega_m^{\eps}(x,t) \Bigr) \, \rmd x \rmd t = 0,
\end{equation*}
and by Theorem~\ref{th:nmtypelimit} the limit is
\begin{align*}
	&\smint_{\Omega_T} \smint_{\sety_{1m}} \bigl( -u_1(x,t,y,\vecs_m) \d_{s_m} \omega_m(x,t,\vecy_m) \\
	&\phantom{ \smint_{\Omega_T} \smint_{\sety_{1m}} \bigl( -u_1(x,t,y,\vecs_m) }
		+ a_0 (x,t,y,\vecs_m) \cdot \grad_y \omega_m(x,t,y,\vecs_{m}) \bigr) \, \rmd \vecs_m \rmd y \rmd x \rmd t = 0.
\end{align*}
Suppose $v_1 \in \cont_{\#}^{\infty}(Y)/\real$ and $c_m \in \cont_{\#}^{\infty}(S_m)$ are the factors of
$\omega_m$ with respect to the $y$ and $s_m$ variables. Utilising the Variational Lemma we then arrive at
\begin{align}
	&\smint_Y \smint_{S_m} \bigl( -u_1(x,t,y,\vecs_m) \, v_1(y) \, \d_{s_m} c_m(s_m) \nonumber \\
	&\phantom{ \smint_Y \smint_{S_m} \bigl( -u_1(x,t,y,\vecs_m) \, v_1(y) }
		+ a_0 (x,t,y,\vecs_m) \cdot \grad_y v_1(y) \, c_m(s_m) \bigr) \, \rmd s_m \rmd y = 0 \label{eq:miniparabolic}
\end{align}
a.e.~on $\Omega_T \times S^{m-1}$, which is our pre-local-problem.

$\bullet$ Suppose $\bigl( \eps, \{ \eps_j' \}_{j=1}^m \bigr) \in \jointly_{\wsep,2+\bar{\ell}-k}^{m \sim k}$
for some $\bar{\ell} \in \scope{k+1,m}$ where $k \in \scope{m-1}$ is required. By definition this means that
$\bigl( \eps, \{ \eps_j' \} \bigr) \in \jointly_{\wsep}^{m \sim k}$ and $\tfrac{\eps_{\bar{\ell}}'}{\eps_k'^2}
\rightarrow  0$ but $\tfrac{\eps_{\bar{\ell}-1}'}{\eps_k'^2} \rightarrow \infty$. We first note that we have already
extracted \eqref{eq:case1astrongonlyk} which at this point carries at least one integral and it happens to be independent
of $\bar{\ell}$. Choose $r_{\eps} = \tfrac{\eps_{i}'}{\eps_k'}$ where $i \in \scope{\bar{\ell},m}$. Apparently, $r_{\eps}
\rightarrow 0$ is guaranteed since $i \in \scope{k+1,m}$. Trivially, $\bigl\{ \tfrac{r_{\eps} \eps_k'}{\eps_{i}'} \bigr\}$
is bounded. Finally,
\begin{align*}
	\tfrac{r_{\eps}}{\eps_k'}	&= \tfrac{\eps_{i}'}{\eps_{\bar{\ell}}'} \tfrac{\eps_{\bar{\ell}}'}{\eps_k'^2} \\
								&\rightarrow 0
\end{align*}	
by assumption and separatedness. Hence, we can utilise \eqref{eq:whatisleftdiff} (with $\ell = i$) giving
\begin{equation*}
	\lim_{\eps \rightarrow 0} \smint_{\Omega_T} \Bigl( -\tfrac{1}{\eps} u_{\eps}(x,t) \d_{s_{i}} \omega_{i}^{\eps}(x,t)
		+ a_{\eps} (x,t) \cdot \tfrac{\eps_{i}'}{\eps_k'^2} \grad_y \omega_{i}^{\eps}(x,t) \Bigr) \, \rmd x \rmd t = 0,
\end{equation*}
and taking the limit by using Theorem~\ref{th:nmtypelimit},
\begin{equation*}
	-\smint_{\Omega_T} \smint_{\sety_{1m}} u_1(x,t,y,\vecs_m) \d_{s_{i}} \omega_{i}(x,t,y,\vecs_{i})
		\, \rmd \vecs_m \rmd y \rmd x \rmd t = 0 \qquad \qquad (i \in \scope{\bar{\ell},m}).
\end{equation*}
Proceeding like before, the equation above leads to the pre-local-problem
\begin{equation}
	\label{eq:indepoftemplocvar}
	-\smint_{S_i} u_1(x,t,y,\vecs_m) \d_{s_{i}} c_{i}(s_i)
		\, \rmd s_i \rmd y = 0 \quad \textrm{for all } c_i \in \cont^{\infty}_{\#}(S_i) \qquad (i \in \scope{\bar{\ell},m}).
\end{equation}
Note that this means that $u_1$ is essentially independent of the temporal local variables $\vecs_{\scope{\bar{\ell},m}}
\in S^{\scope{\bar{\ell},m}}$. Choose now $r_{\eps} = \eps_k'$ (which indeed tends to $0$) and let
$i \in \scope{k+1,\bar{\ell}-1}$ which requires $\bar{\ell} \in \scope{k+2,m}$ (which, of course, in turn requires
$k \in \scope{m-2}$). Then $\bigl\{ \tfrac{r_{\eps} \eps_k'}{\eps_{i}'} \bigr\}$ is bounded since, by assumption
and separatedness,
\begin{align*}
	\tfrac{r_{\eps} \eps_k'}{\eps_{i}'}	&= \tfrac{\eps_k'^2}{\eps_{\bar{\ell}-1}'} \tfrac{\eps_{\bar{\ell}-1}'}{\eps_{i}'} \\
										&\rightarrow 0.
\end{align*}
We have shown that we can employ \eqref{eq:whatisleftdiff} (with $\ell = i$), giving
\begin{equation*}
	\lim_{\eps \rightarrow 0} \smint_{\Omega_T} \Bigl( -\tfrac{1}{\eps} u_{\eps}(x,t)
		\tfrac{\eps_k'^2}{\eps_{\bar{\ell}-1}'} \tfrac{\eps_{\bar{\ell}-1}'}{\eps_{i}'} \d_{s_{i}} \omega_{i}^{\eps}(x,t)
			+ a_{\eps} (x,t) \cdot \grad_y \omega_{i}^{\eps}(x,t) \Bigr) \, \rmd x \rmd t = 0
\end{equation*}
for $i \in \scope{k+1,\bar{\ell}-1}$. Taking the limit by using Theorem~\ref{th:nmtypelimit}, remembering that
$\tfrac{\eps_k'^2}{\eps_{\bar{\ell}-1}'} \tfrac{\eps_{\bar{\ell}-1}'}{\eps_{i}'} \rightarrow 0$ and taking into
consideration the final remark of Remark~\ref{re:nmtypelimit}, we arrive at
\begin{equation*}
	\smint_{\Omega_T} \smint_{\sety_{1m}} a_0 (x,t,y,\vecs_m)
		\cdot \grad_y \omega_{i}(x,t,y,\vecs_{i}) \, \rmd \vecs_m \rmd y \rmd x \rmd t = 0.
\end{equation*}
Proceeding in the same way as in the derivation of \eqref{eq:case1astrong} we get
\begin{equation*}
	\smint_Y \smint_{S^{\scope{i+1,m}}} a_0 (x,t,y,\vecs_m) \, \rmd \vecs_{\scope{i+1,m}} \cdot \grad_y v_1(y) \, \rmd y = 0
		\qquad \qquad (i \in \scope{k+1,\bar{\ell}-1}).
\end{equation*}
We conclude that
\begin{equation}
	\label{eq:lovprobassomeints}
	\smint_Y \smint_{S^{\scope{\bar{\ell},m}}} a_0 (x,t,y,\vecs_m)
		\, \rmd \vecs_{\scope{\bar{\ell},m}} \cdot \grad_y v_1(y) \, \rmd y = 0
\end{equation}
since the case $\bar{\ell} = k+1$ is taken care of by \eqref{eq:case1astrongonlyk}. The extracted pre-local-problems are
\eqref{eq:indepoftemplocvar} and \eqref{eq:lovprobassomeints} in this case.

$\bullet$ Suppose $\bigl( \eps, \{ \eps_j' \}_{j=1}^m \bigr) \in
\jointly_{\wsep,1+m+\ellring-2k}^{m \sim k}$ for some $\ellring
\in \scope{k+2,m}$ where it is required that $k \in \scope{m-2}$. By definition this means
that $\bigl( \eps, \{ \eps_j' \} \bigr) \in \jointly_{\wsep}^{m \sim k}$ and that $\eps_{\ellring-1}' \sim \eps_k'^2$.
Choose $r_{\eps} = \tfrac{\eps_{i}'}{\eps_k'}$ and let $i \in \scope{\ellring,m}$.
It is clearly guaranteed that $r_{\eps} \rightarrow 0$ since $i \in \scope{k+2,m}$. Moreover, it is trivial that
$\bigl\{ \tfrac{r_{\eps} \eps_k'}{\eps_{i}'} \bigr\}$ is bounded. Finally,
\begin{equation*}
	\tfrac{r_{\eps}}{\eps_k'} = \tfrac{\eps_{i}'}{\eps_k'^2}
		= \tfrac{\eps_{i}'}{\eps_{\ellring-1}'} \tfrac{\eps_{\ellring-1}'}{\eps_k'^2} \rightarrow 0
\end{equation*}	
by assumption and separatedness. Hence, we can utilise \eqref{eq:whatisleftdiff}
with $\ell = i$ giving
\begin{equation*}
	\lim_{\eps \rightarrow 0} \smint_{\Omega_T} \Bigl( -\tfrac{1}{\eps}
		u_{\eps}(x,t) \d_{s_{i}} \omega_{i}^{\eps}(x,t) + a_{\eps} (x,t)
			\cdot \tfrac{\eps_{i}'}{\eps_k'^2} \grad_y \omega_{i}^{\eps}(x,t) \Bigr)
				\, \rmd x \rmd t = 0,
\end{equation*}
and taking the limit by using Theorem~\ref{th:nmtypelimit},
\begin{equation*}
	-\smint_{\Omega_T} \smint_{\mathcal{Y}_{1m}} u_1(x,t,y,\vecs_m) \d_{s_{i}}
		\omega_{i}(x,t,y,\vecs_{i}) \, \rmd \vecs_m \rmd y \rmd x \rmd t
			= 0 \qquad \qquad (i \in \scope{\ellring,m}).
\end{equation*}
Proceeding like before, the equation above leads to the pre-local-problem
\begin{equation}
	\label{eq:indepoftemplocvarres}
	-\smint_{S_i} u_1(x,t,y,\vecs_m) \d_{s_{i}} c_i(s_i) \, \rmd s_i = 0
		\quad \textrm{for all }c_i \in \cont^{\infty}_{\#}(S_i) \qquad (i \in \scope{\ellring,m}).
\end{equation}
Note that this means that $u_1$ is essentially independent of the temporal local variables
$\vecs_{\scope{\ellring,m}} \in S^{\scope{\ellring,m}}$. In particular, \eqref{eq:indepoftemplocvarres}
implies that
\begin{equation}
	\label{eq:inparticularresult}
	\smint_{S^{\scope{\ellring,m}}} u_1(x,t,y,\vecs_m) \, \rmd \vecs_{\scope{\ellring,m}} = u_1(x,t,y,\vecs_m)
\end{equation}
holds a.e.~on $\Omega_T \times Y \times S^m$. For the second pre-local-problem, choose
$r_{\eps} = \eps_k'$ and let $i = \ellring-1$. Then
$\bigl\{ \tfrac{r_{\eps} \eps_k'}{\eps_{i}'} \bigr\}$ is bounded
since, by assumption,
\begin{equation}
	\label{eq:boundedbyassumption}
	\tfrac{r_{\eps} \eps_k'}{\eps_{i}'} = \tfrac{\eps_k'^2}{\eps_{\ellring-1}'} \rightarrow 1.
\end{equation}
We have shown that we can employ \eqref{eq:whatisleftdiff}, giving
\begin{equation*}
	\lim_{\eps \rightarrow 0} \smint_{\Omega_T} \Bigl( -\tfrac{1}{\eps}
		u_{\eps}(x,t) \tfrac{\eps_k'^2}{\eps_{\ellring-1}'}
			\d_{s_{\ellring-1}} \omega_{\ellring-1}^{\eps}(x,t)
				+ a_{\eps} (x,t) \cdot \grad_y \omega_{\ellring-1}^{\eps}(x,t) \Bigr) \, \rmd x \rmd t = 0.
\end{equation*}
Taking the limit by using Theorem~\ref{th:nmtypelimit} and \eqref{eq:boundedbyassumption}, we arrive at
\begin{multline*}
	\smint_{\Omega_T} \smint_{\mathcal{Y}_{1m}} \bigl( -u_1(x,t,y,\vecs_m)
		\d_{s_{\ellring-1}} \omega_{\ellring-1}(x,t,y,\vecs_{\ellring-1}) \\
			+ a_0 (x,t,y,\vecs_m) \cdot \grad_y \omega_{\ellring-1}(x,t,y,\vecs_{\ellring-1}) \bigr)
				\, \rmd \vecs_m \rmd y \rmd x \rmd t = 0.
\end{multline*}
Utilising property \eqref{eq:inparticularresult}, this becomes
\begin{multline*}
	\smint_{\Omega_T} \smint_Y \smint_{S^{\ellring-1}} \Bigl( -u_1(x,t,y,\vecs_m)
		\d_{s_{\ellring-1}} \omega_{\ellring-1}(x,t,y,\vecs_{\ellring-1}) \\
			+ \smint_{S^{\scope{\ellring,m}}} a_0 (x,t,y,\vecs_m) \, \rmd \vecs_{\scope{\ellring,m}}
				\cdot \grad_y \omega_{\ellring-1}(x,t,y,\vecs_{\ellring-1}) \Bigr)
					\, \rmd \vecs_{\ellring-1} \rmd y \rmd x \rmd t = 0.
\end{multline*}
Suppose $v_1 \in \cont_{\#}^{\infty}(Y)/\mathbb{R}$ and $c_{\ellring-1} \in \cont_{\#}^{\infty}(S_{\ellring-1})$
are the factors of $ \omega_{\ellring-1}$ with respect to the $y$ and $s_{\ellring-1}$ local variables, respectively.
Employing the Variational Lemma we then get
\begin{multline}
	\label{eq:seclocprobrapres}
	\smint_Y \smint_{S_{\ellring-1}} \Bigl( -u_1(x,t,y,\vecs_m)
		v_1(y) \d_{s_{\ellring-1}} c_{\ellring-1}(s_{\ellring-1}) \\
			+ \smint_{S^{\scope{\ellring,m}}} a_0 (x,t,y,\vecs_m) \, \rmd \vecs_{\scope{\ellring,m}}
				\cdot \grad_y v_1(y)  c_{\ellring-1}(s_{\ellring-1}) \Bigr)
					\, \rmd s_{\ellring-1} \rmd y = 0.
\end{multline}
a.e.~on $\Omega_T \times S^{\ellring-2} \times S^{\scope{\ellring,m}}$, which is our second pre-local-problem.
Concluding the present case, the extracted pre-local-problems are \eqref{eq:indepoftemplocvarres} and \eqref{eq:seclocprobrapres}. \\

What is left to do is to characterise $a_0$ in terms of $a$ such that the pre-local-problems become true local problems,
and for this we introduce a sequence $\{ p_{\mu} \}_{\mu=1}^{\infty}$ in $\distr \bigl( \Omega_T; \cont_{\#}^{\infty}(\sety_{1m})^N \bigr)$
of Evans's perturbed test functions (see \cite{Evan89,Evan92}) defined according to
\begin{equation*}
	p_{\mu} = \pi_{\mu} + {\pi_1}_{\mu} + \delta c \qquad \qquad (\mu \in \integ_+),
\end{equation*}
where $\delta > 0$, $\pi_{\mu} \in \distr (\Omega_T)^N$ and ${\pi_1}_{\mu},c \in
\distr \bigl( \Omega_T; \cont_{\#}^{\infty}(\sety_{1m})^N \bigr)$ for all $\mu \in \integ_+$. Let
$\{ \pi_{\mu} \}_{\mu=1}^{\infty}$ and $\{ {\pi_1}_{\mu} \}_{\mu=1}^{\infty}$ be such that
\begin{align*}
	&\left\{
		\begin{aligned}
			\pi_{\mu}			&\rightarrow \grad u			&	&\qquad \textrm{in } L^2(\Omega_T)^N, \\
			\pi_{\mu} (x,t)	&\rightarrow \grad u (x,t)	&	&\qquad \textrm{a.e.~on } \Omega_T,
		\end{aligned}
	\right. \\
\inter{and}
	&\left\{
		\begin{aligned}
			{\pi_1}_{\mu}						&\rightarrow \grad_y u_1
				&	&\qquad \textrm{in } L^2(\Omega_T \times \sety_{1m})^N, \\
			{\pi_1}_{\mu} (x,t,y,\vecs_m)	&\rightarrow \grad_y u_1 (x,t,y,\vecs_m)
				&	&\qquad \textrm{a.e.~on } \Omega_T \times \sety_{1m}
		\end{aligned}
	\right.
\end{align*}
as $\mu \rightarrow \infty$. Strictly speaking, the last convergence should hold a.e.~on $\Omega_T \times \real^{n+m}$.
By periodicity, this is implied from the given assumption, though. For each fixed $\mu \in \integ_+$, introduce the sequence
$\{ p_{\mu}^{\eps} \}$ defined by
\begin{equation*}
	p_{\mu}^{\eps}(x,t) = p_{\mu} (x,t,\tfrac{x}{\eps},\vect_m^{\eps}) \qquad \qquad ((x,t) \in \Omega_T).
\end{equation*}

A crucial result for the remainder of the proof is
\begin{align}
	&\smint_{\Omega_T} \smint_{\sety_{1m}} \bigl( - a_0 (x,t,y,\vecs_m)
		+ a(x,t,y,\vecs_m ; \grad u + \grad_y u_1 + \delta c) \bigr) \nonumber \\
	&\phantom{ \smint_{\Omega_T} \smint_{\sety_{1m}} \bigl( - a_0 (x,t,y,\vecs_m)
				+ a(x,t,y,\vecs_m ; \grad u + }
		\cdot \delta c (x,t,y,\vecs_m) \, \rmd \vecs_m \rmd y \rmd x \rmd t \geqslant 0 \label{eq:desiredinequality}
\end{align}
for every $\delta > 0$ and every $c \in \distr \bigl( \Omega_T; \cont_{\#}^{\infty}(\sety_{1m})^N \bigr)$.
Hence, let us prove \eqref{eq:desiredinequality}. The point of departure is property (B$_4$)
which implies the inequality
\begin{equation*}
	\bigl( a(x,t,\tfrac{x}{\eps},\vect_m^{\eps}; \grad u_{\eps}) - a(x,t,\tfrac{x}{\eps},\vect_m^{\eps}; p_{\mu}^{\eps}) \bigr)
		\cdot \bigl( \grad u_{\eps}(x,t) - p_{\mu}^{\eps}(x,t) \bigr) \geqslant 0 \qquad \quad ((x,t) \in \Omega_T),
\end{equation*}
which after integration over $\Omega_T$ and expansion of the scalar product becomes
\begin{align*}
	\smint_{\Omega_T} \bigl( &a^{\eps}(x,t; \grad u_{\eps}) \cdot \grad u_{\eps}(x,t)
		- a^{\eps}(x,t; \grad u_{\eps}) \cdot p_{\mu}^{\eps}(x,t) \\
			&\qquad \quad - a^{\eps}(x,t; p_{\mu}^{\eps}) \cdot \grad u_{\eps}(x,t)
				+ a^{\eps}(x,t; p_{\mu}^{\eps}) \cdot p_{\mu}^{\eps}(x,t) \bigr) \, \rmd x \rmd t \geqslant 0.
\end{align*}
We can rewrite the first term by \eqref{eq:weakevolprob1} to obtain
\begin{multline*}
	- \bigl\langle \tfrac{\d}{\d t} u_{\eps} , u_{\eps} \bigr\rangle_{X',X}
		+ \smint_{\Omega_T} f(x,t) \, u_{\eps}(x,t) \,  \rmd x \rmd t \\
			+ \smint_{\Omega_T} \bigl( - a^{\eps}(x,t; \grad u_{\eps}) \cdot p_{\mu}^{\eps}(x,t)
				- a^{\eps}(x,t; p_{\mu}^{\eps}) \cdot \grad u_{\eps}(x,t) \\
					+ a^{\eps}(x,t; p_{\mu}^{\eps}) \cdot p_{\mu}^{\eps}(x,t) \bigr)
						\, \rmd x \rmd t \geqslant 0,
\end{multline*}
which is realised to tend to, as $\eps \rightarrow 0$ and up to a subsequence, the inequality
\begin{align}
	&- \bigl\langle \tfrac{\d}{\d t} u , u \bigr\rangle_{X',X}
		+ \smint_{\Omega_T} f(x,t) \, u(x,t) \, \rmd x \rmd t \nonumber \\
	&\qquad + \smint_{\Omega_T} \smint_{\mathcal{Y}_{1m}} \Bigl( - a_0 (x,t,y,\vecs_m)
		\cdot p_{\mu}(x,t,y,\vecs_m) \nonumber \\
	&\phantom{ \qquad + \smint_{\Omega_T} \smint_{\mathcal{Y}_{1m}} \Bigl( }
		- a(x,t,y,\vecs_m; p_{\mu}) \cdot \bigl( \grad u(x,t) + \grad_y u_1(x,t,y,\vecs_m) \bigr) \nonumber \\
	&\phantom{ \qquad + \smint_{\Omega_T} \smint_{\mathcal{Y}_{1m}} \Bigl( }
		+ a(x,t,y,\vecs_m; p_{\mu}) \cdot p_{\mu}(x,t,y,\vecs_m) \Bigr)
			\, \rmd \vecs_m \rmd y \rmd x \rmd t \geqslant 0 \label{eq:afterepstozero}
\end{align}
since
\begin{equation*}
	\bigl\langle \tfrac{\d}{\d t} u , u \bigr\rangle_{X',X}
		\leqslant \liminf_{\eps \rightarrow 0}
			\, \bigl\langle \tfrac{\d}{\d t} u_{\eps} , u_{\eps} \bigr\rangle_{X',X}
\end{equation*}
(see, e.g., the end of the proof of Theorem~3.1 in \cite{NguWou04}).
We will now investigate what happens when we let $\mu \rightarrow \infty$ in \eqref{eq:afterepstozero}.
Immediately from the assumptions on $\{ p_{\mu} \}_{\mu=1}^{\infty}$ we have, as $\mu \rightarrow \infty$,
\begin{equation}
	\label{eq:perturbedlimit}
	p_{\mu} \rightarrow \grad u + \grad_y u_1 + \delta c \qquad \textrm{in }
		L^2(\Omega_T \times \sety_{1m})^N \textrm{ and a.e.~on } \Omega_T \times \sety_{1m},
\end{equation}
which takes care of the first term of the second integral in \eqref{eq:afterepstozero}.
Moreover, we clearly have
\begin{equation*}
	a(x,t,y,\vecs_m;p_{\mu})
		\rightarrow a(x,t,y,\vecs_m; \grad u + \grad_y u_1 + \delta c) \qquad (\textrm{a.e.~on } \Omega_T \times \sety_{1m}),
\end{equation*}
which takes care of the mid term of the second integral in \eqref{eq:afterepstozero},
and for the last term of the second integral in \eqref{eq:afterepstozero},
\begin{align*}
	a(x,t,y,\vecs_m;p_{\mu}) \cdot p_{\mu}(x,t,y,\vecs_m)
		&\rightarrow a(x,t,y,\vecs_m; \grad u + \grad_y u_1 + \delta c) \\
			&\phantom{ \rightarrow a(x, } \cdot \bigl( \grad u(x,t)
				+ \grad_y u_1(x,t,y,\vecs_m) + \delta c(x,t,y,\vecs_m) \bigr)
\end{align*}
a.e.~on $\Omega_T \times \sety_{1m}$.
The key to come any further is to use Lebesgue's Generalised Dominated Convergence Theorem (LGDCT) on this
last integral term. (See, e.g., Theorem~(19a) on p.~1015 in \cite{ZeidIIB} for the formulation of LGDCT.)
What remains in order to employ LGDCT is to establish majorising, non-negative sequences of functions.
By \eqref{eq:anineqfora} (with $n = 1$), we have
\begin{equation*}
	\bigl| a(x,t,y,\vecs_m; p_{\mu}) \bigr| < C_1 \bigl( 1 + \bigl| p_{\mu}(x,t,y,\vecs_m) \bigr| \bigr)
		\qquad \qquad ((x,t) \in \Omega_T, (y,\vecs_m) \in \sety_{1m}).
\end{equation*}
Hence, by applying this observation and the Cauchy--Schwarz inequality, we have for the last term of
the second integral in \eqref{eq:afterepstozero} the majorisation
\begin{align*}
	&\bigl| a(x,t,y,\vecs_m; p_{\mu}) \cdot p_{\mu}(x,t,y,\vecs_m) \bigr|
		\leqslant \bigl| a(x,t,y,\vecs_m; p_{\mu}) \bigr| \, \bigl| p_{\mu}(x,t,y,\vecs_m) \bigr| \\
	&\phantom{ \bigl| a(x,t,y,\vecs_m; p_{\mu}) \cdot p_{\mu}(x,t,y,\vecs_m) \bigr| }
		< C_1 \bigl( 1 + \bigl| p_{\mu}(x,t,y,\vecs_m) \bigr| \bigr) \bigl| p_{\mu}(x,t,y,\vecs_m) \bigr| \\
	&\phantom{ \bigl| a(x,t,y,\vecs_m; p_{\mu}) \cdot p_{\mu}(x,t,y,\vecs_m) \bigr| }
		= C_1  \Bigl( \bigl| p_{\mu}(x,t,y,\vecs_m) \bigr| + \bigl| p_{\mu}(x,t,y,\vecs_m) \bigr|^2 \Bigr)
\end{align*}
a.e.~on $\Omega_T \times \sety_{1m}$. Due to \eqref{eq:perturbedlimit}, the majorising right-hand
side fulfils, as $\mu \rightarrow \infty$, both
\begin{align*}
	&C_1  \Bigl( \bigl| p_{\mu}(x,t,y,\vecs_m) \bigr| + \bigl| p_{\mu}(x,t,y,\vecs_m) \bigr|^2 \Bigr) \\
	&\phantom{ C_1  \Bigl( \bigl| p_{\mu}(x,t,y,\vecs_m) \bigr| }
		\rightarrow C_1  \Bigl( \bigl| \grad u(x,t) + \grad_y u_1(x,t,y,\vecs_m) + \delta c(x,t,y,\vecs_m) \bigr| \\
	&\phantom{ C_1  \Bigl( \bigl| p_{\mu}(x,t,y,\vecs_m) \bigr|
				\rightarrow C_1  \Bigl( }	
		+ \bigl| \grad u(x,t) + \grad_y u_1(x,t,y,\vecs_m) + \delta c(x,t,y,\vecs_m) \bigr|^2 \Bigr),
\end{align*}
a.e.~on $\Omega_T \times \sety_{1m}$, and
\begin{align*}
	&\smint_{\Omega_T} \smint_{\mathcal{Y}_{1m}} C_1  \Bigl( \bigl| p_{\mu}(x,t,y,\vecs_m) \bigr|
		+ \bigl| p_{\mu}(x,t,y,\vecs_m) \bigr|^2 \Bigr) \, \rmd \vecs_m \rmd y \rmd x \rmd t \\
	&\qquad \rightarrow \smint_{\Omega_T} \smint_{\mathcal{Y}_{1m}} C_1  \Bigl( \bigl| \grad u(x,t)
		+ \grad_y u_1(x,t,y,\vecs_m) + \delta c(x,t,y,\vecs_m) \bigr| \\
	&\phantom{ \qquad \rightarrow \smint_{\Omega_T} \smint_{\mathcal{Y}_{1m}} C_1  \Bigl( }
		+ \bigl| \grad u(x,t) + \grad_y u_1(x,t,y,\vecs_m) + \delta c(x,t,y,\vecs_m) \bigr|^2 \Bigr)
			\, \rmd \vecs_m \rmd y \rmd x \rmd t;
\end{align*}
thus, LGDCT is applicable. Hence, by finally utilising LGDCT, \eqref{eq:afterepstozero} converges to the inequality
\begin{align*}
	&- \bigl\langle  \tfrac{\d}{\d t} u , u \bigr\rangle_{X',X}
		+ \smint_{\Omega_T} f(x,t) \, u(x,t) \, \rmd x \rmd t \\
	&\quad + \smint_{\Omega_T} \smint_{\mathcal{Y}_{1m}} \Bigl( - a_0(x,t,y,\vecs_m)
		\cdot \bigl( \grad u(x,t) + \grad_y u_1(x,t,y,\vecs_m) + \delta c(x,t,y,\vecs_m) \bigr) \\
	&\phantom{ \quad + \smint_{\Omega_T} \smint_{\mathcal{Y}_{1m}} \Bigl( }
		- a(x,t,y,\vecs_m ; \grad u + \grad_y u_1 + \delta c)
			\cdot \bigl( \grad u (x,t) + \grad_y u_1 (x,t,y,\vecs_m) \bigr) \\
	&\phantom{ \quad + \smint_{\Omega_T} \smint_{\mathcal{Y}_{1m}} \Bigl( }
		+ a(x,t,y,\vecs_m ; \grad u + \grad_y u_1 + \delta c) \\
	&\phantom{ \quad + \smint_{\Omega_T} \smint_{\mathcal{Y}_{1m}} \Bigl( } \quad
		\cdot \bigl( \grad u(x,t) + \grad_y u_1(x,t,y,\vecs_m) + \delta c(x,t,y,\vecs_m) \bigr) \Bigr)
			\, \rmd \vecs_m \rmd y \rmd x \rmd t \geqslant 0.
\end{align*}
The inequality above can be written
\begin{multline}
	\label{eq:firsttermshouldvanish}
	\smint_{\Omega_T} \smint_{\sety_{1m}} \bigl( - a_0(x,t,y,\vecs_m) \cdot \grad_y u_1(x,t,y,\vecs_m)
		- a_0(x,t,y,\vecs_m) \cdot \delta c(x,t,y,\vecs_m) \\
			+ a(x,t,y,\vecs_m ; \grad u + \grad_y u_1 + \delta c) \cdot \delta c(x,t,y,\vecs_m) \bigr)
				\, \rmd \vecs_m \rmd y \rmd x \rmd t \geqslant 0,
\end{multline}
where we have used \eqref{eq:weakevolproblim2} to lose the $\langle \tfrac{\d}{\d t} u , u \rangle$
and the $\int f u$ terms. We want to lose the first term in
the integrand, and in order to achieve this we must utilise the pre-local-problems.

$\bullet$ Suppose $\bigl( \eps, \{ \eps_j' \}_{j=1}^m \bigr) \in \jointly_{\wsep,1}^{m \sim k}$.
By density, the pre-local-problem \eqref{eq:locprobfirst} holds for all $v_1 \in \wrum = H_{\#}^1 (Y) / \real$.
(The density property follows from the fact that $H_{\#}^1 (Y)$ is defined to be the closure
of $\cont_{\#}^{\infty}(Y)$ in the $H^1(Y)$-norm; see, e.g., Definition~3.48 in \cite{CioDon99}.) Hence, since
$u_1(x,t, \vecs_m) \in \wrum = H_{\#}^1 (Y) / \real$ a.e.~on $\Omega_T \times S^m$,
\begin{equation*}
	- \smint_{\Omega_T}\smint_{\sety_{1m}} a_0 (x,t,y,\vecs_m)
		\cdot \grad_y u_1 (x,t,y,\vecs_m) \, \rmd \vecs_m \rmd y \rmd x \rmd t = 0,
\end{equation*}
i.e., the first term in the integrand of \eqref{eq:firsttermshouldvanish} gives no contribution in this case.

$\bullet$ Suppose $\bigl( \eps, \{ \eps_j' \}_{j=1}^m \bigr) \in \jointly_{\wsep,2}^{m \sim k}$.
The pre-local-problem \eqref{eq:miniparabolic} can be written
\begin{multline*}
	\smint_Y \smint_{S_m} \bigl( u_1(x,t,y,\vecs_m) \, \d_{s_m} \omega(y,s_m) \\
		- a_0 (x,t,y,\vecs_m) \cdot \grad_y \omega(y,s_m) \bigr) \, \rmd s_m \rmd y
			= 0 \qquad \quad (\textrm{a.e.~on } \Omega_T \times S^{m-1}),
\end{multline*}
i.e.,
\begin{multline*}
		- \smint_Y \smint_{S_m} a_0 (x,t,y,\vecs_m) \cdot \grad_y \omega(y,s_m) \, \rmd s_m \rmd y \\
			= - \smint_Y \smint_{S_m} u_1(x,t,y,\vecs_m) \, \d_{s_m} \omega(y,s_m) \, \rmd s_m \rmd y
				\qquad \quad (\textrm{a.e.~on } \Omega_T \times S^{m-1})
\end{multline*}
for all $\omega \in \bigl( \cont_{\#}^{\infty}(Y)/\real \bigr) \odot \cont_{\#}^{\infty}(S_m)$ and hence,
by the density result of Lemma~\ref{lem:tensprodspdenseinsobsp} and the fact that the tensor product set
spans the corresponding tensor product space, for all $\omega \in H_{\#}^1 (S_m;\wrum,\wrum')$. In this case
we have by assumption that $u_1 \in L^2 \bigl( \Omega_T \times S^{m-1};  H_{\#}^1 (S_m;\wrum,\wrum') \bigr)$,
which implies $u_1(x,t,\vecs_{m-1}) \in H_{\#}^1 (S_m;\wrum,\wrum')$ a.e.~on $\Omega_T \times S^{m-1}$. Thus,
\begin{align*}
	&- \smint_{\Omega_T} \smint_{\sety_{1m}} a_0 (x,t,y,\vecs_m)
		\cdot \grad_y u_1(x,t,y,\vecs_{m}) \, \rmd \vecs_m \rmd y \rmd x \rmd t \\
	&\phantom{ - \smint_{\Omega_T} \smint_{\sety_{1m}} }
		= \smint_{\Omega_T} \smint_{S^{m-1}} \Bigl( - \smint_Y \smint_{S_m} a_0 (x,t,y,\vecs_m)
			\cdot \grad_y u_1(x,t,y,\vecs_{m}) \, \rmd s_m \rmd y \Bigr) \, \rmd \vecs_{m-1} \rmd x \rmd t  \\
	&\phantom{ - \smint_{\Omega_T} \smint_{\sety_{1m}} }
		= - \smint_{\Omega_T} \smint_{S^{m-1}} \bigl\langle \d_{s_m} u_1 (x,t,\vecs_{m-1}),
			u_1 (x,t,\vecs_{m-1}) \bigr\rangle_{L_{\#}^2(S_m; \wrum'),L_{\#}^2(S_m; \wrum)}
				\, \rmd \vecs_{m-1} \rmd x \rmd t .
\end{align*}
By Lemma~\ref{lem:vanishintegrsum}, the duality pairing in the right-hand side
vanishes, so also in this case the first term in the integrand of
\eqref{eq:firsttermshouldvanish} gives no contribution.

$\bullet$ Suppose $\bigl( \eps , \{ \eps_j' \}_{j=1}^m \bigr) \in
\jointly_{\wsep,2+\bar{\ell}-k}^{m \sim k}$ for some $\bar{\ell} \in \scope{k+1,m}$
where $k \in \scope{m-1}$ is required. By density, the pre-local-problem
\eqref{eq:lovprobassomeints} becomes
\begin{equation}
	\label{eq:givesvanishing}
	- \smint_Y \smint_{S^{\scope{\bar{\ell},m}}} a_0 (x,t,y,\vecs_m)
		\, \rmd \vecs_{\scope{\bar{\ell},m}} \cdot \grad_y v_1 (y) \, \rmd y
			= 0 \qquad (\textrm{a.e.~on } \Omega_T \times S^{\bar{\ell}-1})
\end{equation}
for all $v_1 \in \wrum = H_{\#}^1 (Y) / \real$. Since $u_1$ is almost everywhere constant with
respect to $\vecs_{\scope{\bar{\ell},m}} \in S^{\scope{\bar{\ell},m}}$ due to the pre-local-problem
\eqref{eq:indepoftemplocvar}, and $u_1(x,t, \vecs_m) \in \wrum = H_{\#}^1 (Y) / \real$ a.e.~on
$\Omega_T \times S^m$, we have
\begin{multline*}
	- \smint_{\Omega_T} \smint_{\sety_{1m}} a_0 (x,t,y,\vecs_m)
		\cdot \grad_y u_1(x,t,y,\vecs_{m}) \, \rmd \vecs_m \rmd y \rmd x \rmd t \\
			= \smint_{\Omega_T} \smint_{S^{\bar{\ell}-1}} \Bigl( - \smint_Y
				\smint_{S^{\scope{\bar{\ell},m}}} a_0 (x,t,y,\vecs_m) \, \rmd \vecs_{\scope{\bar{\ell},m}}
					\cdot \grad_y u_1(x,t,y,\vecs_{m}) \, \rmd y \Bigr) \, \rmd \vecs_{\bar{\ell}-1} \rmd x \rmd t,
\end{multline*}
which clearly vanishes due to \eqref{eq:givesvanishing}. Again, the first term in the integrand of
\eqref{eq:firsttermshouldvanish} gives no contribution.

$\bullet$ Suppose $\bigl( \eps , \{ \eps_j' \}_{j=1}^m \bigr) \in
\jointly_{\wsep,1+m+\ellring-2k}^{m \sim k}$ for some $\ellring
\in \scope{k+2,m}$ where $k \in \scope{m-2}$ is required. The pre-local-problem
\eqref{eq:seclocprobrapres} may be written as
\begin{multline*}
	\smint_Y \smint_{S_{\ellring-1}} \Bigl( u_1(x,t,y,\vecs_m)
		\d_{s_{\ellring-1}} \omega (y,s_{\ellring-1}) \\
			- \smint_{S^{\scope{\ellring,m}}} a_0 (x,t,y,\vecs_m) \, \rmd \vecs_{\scope{\ellring,m}}
				\cdot \grad_y \omega (y,s_{\ellring-1}) \Bigr)
					\, \rmd s_{\ellring-1} \rmd y = 0
\end{multline*}
a.e.~on $\Omega_T \times S^{\ellring-2} \times S^{\scope{\ellring,m}}$ for all $\omega \in
\bigl( \cont_{\#}^{\infty}(Y)/\real \bigr) \odot \cont_{\#}^{\infty}(S_{\ellring-1})$, i.e.,
\begin{multline*}
	- \smint_Y \smint_{S^{\scope{\ellring-1,m}}} a_0 (x,t,y,\vecs_m) \cdot
		\grad_y \omega (y,s_{\ellring-1}) \, \rmd \vecs_{\scope{\ellring-1,m}} \rmd y \\
			= - \smint_Y \smint_{S_{\ellring-1}} u_1(x,t,y,\vecs_m) \d_{s_{\ellring-1}}
				\omega (y,s_{\ellring-1}) \, \rmd s_{\ellring-1} \rmd y
\end{multline*}
a.e.~on $\Omega_T \times S^{\ellring-2} \times S^{\scope{\ellring,m}}$ for all $\omega \in
\bigl( \cont_{\#}^{\infty}(Y)/\real \bigr) \odot \cont_{\#}^{\infty}(S_{\ellring-1})$
and hence, by the density result of Lemma~\ref{lem:tensprodspdenseinsobsp} and the fact that the tensor
product set spans the corresponding tensor product space, for all $\omega \in H_{\#}^1 (S_{\ellring-1};\wrum,\wrum')$.
By assumption, $u_1 \in L^2 \bigl( \Omega_T \times S^{\ellring-2} \times S^{\scope{\ellring,m}};
H_{\#}^1 (S_{\ellring-1};\wrum,\wrum') \bigr)$, implying $u_1 \in H_{\#}^1 (S_{\ellring-1};\wrum,\wrum')$
a.e.~on $\Omega_T \times S^{\ellring-2} \times S^{\scope{\ellring,m}}$. Thus,
\begin{align*}
	&- \smint_{\Omega_T} \smint_{\mathcal{Y}_{1m}} a_0 (x,t,y,\vecs_m) \cdot
		\grad_y u_1(x,t,y,\vecs_{m}) \, \rmd \vecs_m
			\rmd y \rmd x \rmd t \\
	&\phantom{ - \smint_{\Omega_T} \smint_{\mathcal{Y}_{1m}} }
		= \smint_{\Omega_T} \smint_{S^{\ellring-2}}
			\Bigl( - \smint_Y \smint_{S^{\scope{\ellring-1,m}}} a_0 (x,t,y,\vecs_m) \\
	&\phantom{	- \smint_{\Omega_T} \smint_{\mathcal{Y}_{1m}}
				= \smint_{\Omega_T} \smint_{S^{\ellring-2}} \Bigl( - \smint_Y \smint_{S^{\scope{\ellring-1,m}}} }
		\cdot \grad_y u_1(x,t,y,\vecs_{m}) \, \rmd \vecs_{\scope{\ellring-1,m}}
			\rmd y \Bigr) \, \rmd \vecs_{\ellring-2} \rmd x \rmd t \\
	&\phantom{ - \smint_{\Omega_T} \smint_{\mathcal{Y}_{1m}} }
		= - \smint_{\Omega_T} \smint_{S^{\ellring-2}}
		\bigl\langle \d_{s_{\ellring-1}} u_1, u_1 \bigr\rangle_{L_{\#}^2(S_{\ellring-1};\wrum'),L_{\#}^2(S_{\ellring-1};\wrum)}
		\, \rmd \vecs_{\ellring-2} \rmd x \rmd t.
\end{align*}
By Lemma~\ref{lem:vanishintegrsum}, the duality pairing in the right-hand
side vanishes implying that the first term in the integrand of
\eqref{eq:firsttermshouldvanish} gives no contribution.

To conclude, we have proven the inequality \eqref{eq:desiredinequality} for all considered cases. 
Divide \eqref{eq:desiredinequality} by $\delta$, let $\delta \rightarrow 0$ and finally use the
Variational Lemma. Then we clearly have
\begin{equation*}
	a_0 (x,t,y,\vecs_m) = a(x,t,y,\vecs_m; \grad u + \grad_y u_1)
		\qquad \qquad (\textrm{a.e.~on } \Omega_T \times \sety_{1m})
\end{equation*}
as desired. This establishes an H$_\textrm{MP}$-limit $b$ on the form \eqref{eq:hmpliminth}. Since
$u$ is the unique solution to the homogenised equation and $u_1$ is the unique solution to the local
problems, the convergences \eqref{eq:strongconvofueps}--\eqref{eq:1mconvofueps} hold not only for the
extracted subsequence but for the whole sequence as well. The proof is complete.
\end{proof}
\begin{rema}
The assumption $u_1 \in L^2 \bigl( \Omega_T \times S^{m-1}; H_{\#}^1 (S_m;\wrum,\wrum') \bigr)$ in
the slow resonant case $\jointly_{\wsep,2}^{m \sim k}$ merely amounts to the supposition
$\d_{s_m} u_1 \in L^2 \bigl( \Omega_T \times S^{m-1}; L_{\#}^2(S_m; \wrum') \bigr)$ since we
already know $u_1 \in L^2 \bigl( \Omega_T \times S^{m-1}; L_{\#}^2(S_m; \wrum) \bigr)$ as a
fact due to Theorem~\ref{th:gradchar} (with $n=1$).
Similarly, in the rapid resonant case $\mathcal{J}_{\mathrm{wsep},1+m+\ellring-2k}^{m \sim k}$,
$\ellring \in \scope{k+2,m}$, the assumption $u_1 \in L^2 \bigl( \Omega_T \times S^{\ellring-2}
\times S^{\scope{\ellring,m}}; H_{\#}^1 (S_{\ellring-1};\wrum,\wrum') \bigr)$ boils down to requiring
$\d_{s_{\ellring-1}} u_1 \in L^2 \bigl( \Omega_T \times S^{\ellring-2} \times S^{\scope{\ellring,m}};
L_{\#}^2(S_{\ellring-1}; \wrum') \bigr)$.
\end{rema}
Define $\scope{\ell}_0 = \scope{\ell} \cup \{ 0 \} = \{ 0,1,\ldots,\ell\}$ for any non-negative
integer $\ell$. Let $k \in \scope{m}_0$. Define $\jointly_{\wsep}^{m \prec k}$ to be
the set of all pairs $\bigl( \eps , \{ \eps_j' \}_{j=1}^m \bigr)$ of lists in
$\jointly_{\wsep}^{1m}$ such that
\begin{equation*}
	\left\{
		\begin{aligned}
			&\{ \eps \, , \eps_1',\ldots,\eps_m' \}
				& &\qquad \textrm{if } k = 0, \\
			&\{ \eps_1',\ldots,\eps_k', \, \eps \, , \eps_{k+1}',\ldots,\eps_m' \}
				& &\qquad \textrm{if } k \in \scope{m-1}, \\
			&\{ \eps_1',\ldots,\eps_m' \, , \, \eps \}
				& &\qquad \textrm{if } k = m
		\end{aligned}
	\right.
\end{equation*}
is a well-separated list of scale functions. (Hence, for small enough $\eps$, $\eps < \eps_k'$,
explaining the notation ``$\prec k$''. This could be read as ``the spatial scale is asymptotically
less than the $k$-th temporal scale''.) Define the collection $\bigl\{ \jointly_{\wsep,i}^{m \prec k}
\bigr\}_{i=1}^{1+2(m-k)}$ of $1 + 2(m - k)$ subsets of $\jointly_{\wsep}^{m \prec k}$ according to
\begin{itemize}
	\item{} $\jointly_{\wsep,1}^{m \prec k} = \Bigl\{ \bigl( \eps , \{ \eps_j' \}_{j=1}^m \bigr)
				\in \jointly_{\wsep}^{m \prec k} \, : \, \tfrac{\eps^2}{\eps_m'} \rightarrow 0 \Bigr\}$,
	\item{} $\jointly_{\wsep,2}^{m \prec k} = \Bigl\{ \bigl( \eps , \{ \eps_j' \}_{j=1}^m \bigr)
				\in \jointly_{\wsep}^{m \prec k} \, : \, \eps_m' \sim \eps^2 \Bigr\}$,
	\item{} $\jointly_{\wsep,2+i-k}^{m \prec k} = \Bigl\{ \bigl( \eps , \{ \eps_j' \}_{j=1}^m \bigr)
				\in \jointly_{\wsep}^{m \prec k} : \tfrac{\eps_i'}{\eps^2} \rightarrow 0
					\textrm{ but } \tfrac{\eps_{i-1}'}{\eps^2} \rightarrow \infty \Bigr\}$
						\quad $\matris{\textrm{(} \, i \in \scope{k+1,m}, \\ \phantom{ \textrm{(} \, } (k,i) \neq (0,1) \, \textrm{)}, }$
	\item{} $\jointly_{\wsep,1+m+\iring-2k}^{m \prec k} = \Bigl\{ \bigl( \eps ,
				\{ \eps_j' \}_{j=1}^m \bigr) \in \jointly_{\wsep}^{m \sim k}
				\, : \, \eps_{\iring-1}' \sim \eps^2 \Bigr\}$ \qquad \quad \; \, \, ($\iring \in \scope{k+2,m}$),
\end{itemize}
and
\begin{equation}
	\label{eq:kayeyezeroone}
	\jointly_{\wsep,3}^{m \prec 0}
	= \Bigl\{ \bigl( \varepsilon , \{ \varepsilon_j' \}_{j=1}^m \bigr)
			\in \jointly_{\wsep}^{m \prec 0} \, : \,
				\tfrac{\varepsilon_1'}{\varepsilon^2} \rightarrow 0 \Bigr\}
\end{equation}
for $(k,i) = (0,1)$. Actually, $\jointly_{\wsep,3}^{m \prec k}$ does not really need the second condition---i.e., the
non-convergence to $0$---since it is already implied by the fact that we are in $\jointly_{\wsep}^{m \prec k}$.
Since there does not exist any ``$\varepsilon_0'$'', we note that we need to impose a special definition \eqref{eq:kayeyezeroone}
for $\jointly_{\wsep,3}^{m \prec 0}$ without the extra condition. The collection  $\bigl\{ \jointly_{\wsep,i}^{m \prec k}
\bigr\}_{i=1}^{1+2(m-k)}$ of subsets of $\jointly_{\wsep}^{m \prec k}$ is clearly mutually disjoint. (Note that
if $k = m$, the introduced collection of subsets of $\jointly_{\wsep}^{m \prec m}$ reduces to merely
$\bigl\{ \jointly_{\wsep,1}^{m \prec m} \bigr\}$.)

The theorem below is a modification of Theorem~\ref{th:maintheorem} where the spatial scale function is not
allowed to coincide with any temporal scale function.
\begin{theo}
\label{th:maintheorem2}
Let $k \in \scope{m}_0$. Suppose that the pair $e = \bigl( \eps , \{ \eps_j' \}_{j=1}^m \bigr)$ of lists of
spatial and temporal scale functions belongs to $\bigcup_{i=1}^{1+2(m-k)} \jointly_{\wsep,i}^{m \prec k}$. Let
$\{ u_{\eps} \}$ be the sequence of weak solutions in $H^1 \bigl( \ZeToT; H_0^1(\Omega), H^{-1}(\Omega) \bigr)$
to the evolution problem \eqref{eq:monoparaprob} with $a \, : \, \bar{\Omega}_T \times \real^{N+m} \times \real^N
\rightarrow \real^N$ satisfying the structure conditions \emph{(B}$_1$\emph{)}--\emph{(B}$_5$\emph{)}. Then
\begin{align*}
	u_{\eps}		&\rightarrow	u	&	&\textrm{in } L^2(\Omega_T), \\
	u_{\eps}		&\rightharpoonup  u	&	&\textrm{in } L^2 \bigl( \ZeToT; H_0^1(\Omega) \bigr), \\
\inter{and}
 	\grad u_{\eps}	&\scaleconv{(2,m+1)} \grad u + \grad_{y} u_1, & &
\end{align*}
where $u \in H^1 \bigl( \ZeToT; H_0^1(\Omega), H^{-1}(\Omega) \bigr)$ and $u_1 \in
L^2 ( \Omega_T \times S^m; \wrum )$. Here $u$ is the unique weak solution to
the homogenised problem \eqref{eq:defoflim} with the homogenised flux $b \, : \,
\bar{\Omega}_T \times \real^N \rightarrow \real^N$ given by
\begin{equation*}
	b(x,t;\grad u) = \smint_{\sety_{1m}} a(x,t, y,\vecs_m; \grad u + \grad_{y} u_1) \, \rmd \vecs_m \rmd y.
\end{equation*}
Moreover, we have the following characterisation of $u_1$:

$\bullet$ If $e \in \jointly_{\wsep,1}^{m \prec k}$ then the function $u_1$ is
the unique weak solution to the local problem
\begin{equation*}
	- \grad_y \cdot a(x,t,y,\vecs_m; \grad u + \grad_y u_1) = 0.
\end{equation*}

$\bullet$ If $e \in \jointly_{\wsep,2}^{m \prec k}$, assuming $u_1 \in L^2 \bigl( \Omega_T \times S^{m-1};
H_{\#}^1(S_m;\wrum,\wrum') \bigr)$, then the function $u_1$ is the unique weak solution to the local problem
\begin{equation*}
	\d_{s_m} u_1(x,t,y,\vecs_m) - \grad_y \cdot a(x,t,y,\vecs_m; \grad u + \grad_y u_1)	= 0.
\end{equation*}

$\bullet$ If $e \in \jointly_{\wsep,2+\bar{\ell}-k}^{m \prec k}$ for some $\bar{\ell} \in \scope{k+1,m}$,
provided $k \in \scope{m-1}_0$, then the function $u_1$ is the unique weak solution to the system of local problems
\begin{equation*}
	\left\{
		\begin{aligned}
			- \grad_y \cdot \smint_{S^{\scope{\bar{\ell},m}}}
				a(x,t,y,\vecs_m; \grad u + \grad_y u_1) \, \rmd \vecs_{\scope{\bar{\ell},m}}	&= 0, \\
				\d_{s_i} u_1(x,t,y,\vecs_m)	&= 0 \qquad (i \in \scope{\bar{\ell},m}).
		\end{aligned}
	\right.
\end{equation*}

$\bullet$ If $e \in \jointly_{\wsep,1+m+\ellring-2k}^{m \prec k}$ for
some $\ellring \in \scope{k+2,m}$, provided $k \in \scope{m-2}_0$ and assuming $u_1 \in
L^2 \bigl( \Omega_T \times S^{\ellring-2} \times S^{\scope{\ellring,m}};H_{\#}^1(S_{\ellring-1};\mathcal{W},\mathcal{W}') \bigr)$,
then the function $u_1$ is the unique weak solution to the system of local problems
\begin{equation*}
 	\left\{
 		\begin{aligned}
			\d_{s_{\ellring-1}} u_1(x,t,y,\vecs_m)
				- \grad_y \! \cdot \! \smint_{S^{\scope{\ellring,m}}}
					a(x,t,y,\vecs_m; \grad u + \grad_y u_1)
						\, \rmd \vecs_{\scope{\ellring,m}}	&= 0, \\
			\d_{s_i} u_1(x,t,y,\vecs_m)						&= 0
			\qquad (i \in \scope{\ellring,m}).
 		\end{aligned}
 	\right.
\end{equation*}
\end{theo}
\begin{proof}[\boldproof]
Let $\hat{m} = m+1$ and $\hat{k} = k+1$. (Note that $\hat{k} \in \scope{\hat{m}}$ since $k \in \scope{m}_0$.)
Introduce the list $\{ \hat{\eps}_j' \}_{j=1}^{\hat{m}}$ of $\hat{m}$ new temporal scale functions defined
according to
\begin{equation*}
	\left\{
		\begin{aligned}
			&\hat{\eps}_1' = \eps, \; \hat{\eps}_j' = \eps_{j-1}' \textrm{ for } j \in \scope{2,\hat{m}}
				& &\quad \!\! \textrm{if } \hat{k} = 1, \\
			&\hat{\eps}_j' = \eps_j' \textrm{ for } j \in \scope{\hat{k}-1}, \; \hat{\eps}_{\hat{k}}' = \eps,
				\textrm{ and } \hat{\eps}_j' = \eps_{j-1}' \textrm{ for } j \in \scope{\hat{k}+1,\hat{m}}
				& &\quad \!\! \textrm{if } \hat{k} \in \scope{2,\hat{m}-1}, \\
			&\hat{\eps}_j' = \eps_j' \textrm{ for } j \in \scope{\hat{m}-1},
				\textrm{ and } \hat{\eps}_{\hat{m}}'	= \eps
				& &\quad \!\! \textrm{if } \hat{k} = \hat{m}.
		\end{aligned}
	\right.
\end{equation*}
Since $\bigl( \eps , \{ \eps_j' \}_{j=1}^m \bigr) \in \jointly_{\wsep}^{m \prec k}$ it must
thus equivalently hold that $\bigl( \eps , \{ \hat{\eps}_j' \}_{j=1}^{\hat{m}} \bigr) \in
\jointly_{\wsep}^{\hat{m} \sim \hat{k}}$. Define $\hat{a} \, : \, \bar{\Omega}_T \times
\real^{N+\hat{m}} \times \real^N \rightarrow \real^N$ according to
\begin{equation*}
	\hat{a} (x,t,y,\hat{\vecs}_{\hat{m}};q) = a (x,t,y,\vecs_m;q)
		\qquad \qquad ((x,t) \in \bar{\Omega}_T, (y,\vecs_m) \in \sety_{1m}, q \in \real^N),
\end{equation*}
where we define (provided $\hat{k} \in \scope{2,\hat{m}-1}$)
\begin{equation*}
			\hat{\vecs}_{\hat{m}}	= (\vecs_{\hat{k}-1},\hat{s}_{\hat{k}},\vecs_{\scope{\hat{k},\hat{m}-1}})
											\qquad	\qquad  (\vecs_{\hat{m}-1} = \vecs_m \in S^m = S^{\hat{m}-1})
\end{equation*}
for any $\hat{s}_{\hat{k}} \in \hat{S}_{\hat{k}} = (0,1)$. (The cases $\hat{k} = 1$ and $\hat{k} = \hat{m}$
require obvious respective modifications of the definition.) This means that $\hat{a}$ is in fact independent
of $\hat{s}_{\hat{k}} \in \hat{S}_{\hat{k}}$, though not manifestly so. Furthermore, define
$\hat{\sety}_{1\hat{m}} = Y \times \hat{S}^{\hat{m}}$ where (provided $\hat{k} \in \scope{2,\hat{m}-1}$)
\begin{equation*}
	\hat{S}^{\hat{m}} = S^{\hat{k}-1} \times \hat{S}_{\hat{k}} \times S^{\scope{\hat{k},\hat{m}-1}}.
\end{equation*}
(The cases $\hat{k} = 1$ and $\hat{k} = \hat{m}$ require obvious respective modifications of the definition.)

It is clear that since $a$ satisfies (B$_1$)--(B$_5$), so does $\hat{a}$. Let $\{ \hat{u}_{\eps} \}$ be the
sequence of weak solutions in $H^1 \bigl( \ZeToT; H_0^1(\Omega), H^{-1}(\Omega) \bigr)$ to the evolution problem
\eqref{eq:monoparaprob} with $\hat{a}$ instead of $a$. (Note that $\hat{u}_{\eps} = u_{\eps}$ since $\hat{a} = a$.)
By Theorem~\ref{th:maintheorem} (with ``hatted'' quantities) we then get
\begin{align*}
	\hat{u}_{\eps}		&\rightarrow	\hat{u}		& & \textrm{in } L^2(\Omega_T), \\
	\hat{u}_{\eps}		&\rightharpoonup  \hat{u}	& & \textrm{in } L^2 \bigl( \ZeToT; H_0^1(\Omega) \bigr), \\
\inter{and}
 	\grad \hat{u}_{\eps}	&\scaleconv{(2,m+2)} \grad \hat{u} + \grad_y \hat{u}_1, & &
\end{align*}
where $\hat{u} \in H^1 \bigl( \ZeToT; H_0^1(\Omega), H^{-1}(\Omega) \bigr)$ and $\hat{u}_1 \in L^2 ( \Omega_T
\times \hat{S}^{\hat{m}}; \wrum )$. Here $\hat{u}$ is the
unique weak solution to the homogenised problem \eqref{eq:defoflim} but with the homogenised flux
$\hat{b} \, : \, \bar{\Omega}_T \times \real^N \rightarrow \real^N$ given by
\begin{equation*}
	\hat{b}(x,t,\grad \hat{u}) = \smint_{\hat{\sety}_{1\hat{m}}} \hat{a}(x,t,y,\hat{\vecs}_{\hat{m}};
										\grad \hat{u} + \grad_y \hat{u}_1) \, \rmd \hat{\vecs}_{\hat{m}} \rmd y,
\end{equation*}
and $\hat{u}_1$ is the unique weak solution to the local problems
\begin{equation}
	\label{eq:firstlochat}
	- \grad_y \cdot \hat{a}(x,t,y,\hat{\vecs}_{\hat{m}}; \grad \hat{u}
		+ \grad_y \hat{u}_1)	= 0
\end{equation}
if $\bigl( \eps , \{ \hat{\eps}_j' \}_{j=1}^{\hat{m}} \bigr) \in \jointly_{\wsep,1}^{\hat{m}
\sim \hat{k}}$;
\begin{equation}
	\label{eq:seclochat}
	\d_{\hat{s}_{\hat{m}}} \hat{u}_1(x,t,y,\hat{\vecs}_{\hat{m}})
		- \grad_y \cdot \hat{a}(x,t,y,\hat{\vecs}_{\hat{m}}; \grad \hat{u} + \grad_y \hat{u}_1) = 0
\end{equation}
if $\bigl( \eps , \{ \hat{\eps}_j' \}_{j=1}^{\hat{m}} \bigr) \in \jointly_{\wsep,2}^{\hat{m}
\sim \hat{k}}$ and assuming $\hat{u}_1 \in L^2 \bigl( \Omega_T \times \hat{S}^{\hat{m}-1};
H_{\#}^1(\hat{S}_{\hat{m}};\wrum,\wrum') \bigr)$;
\begin{equation}
	\label{eq:thirdlochat}
	\left\{
		\begin{aligned}
			- \grad_y \cdot \smint_{\hat{S}^{\scope{\hat{\ell},\hat{m}}}}
				\hat{a}(x,t,y,\hat{\vecs}_{\hat{m}}; \grad \hat{u}
					+ \grad_y \hat{u}_1) \, \rmd \hat{\vecs}_{\scope{\hat{\ell},\hat{m}}}	&= 0, \\
			\d_{\hat{s}_{\hat{\imath}}} \hat{u}_1(x,t,y,\hat{\vecs}_{\hat{m}})				&= 0
				\qquad (\hat{\imath} \in \scope{\hat{\ell},\hat{m}})
		\end{aligned}
	\right.
\end{equation}
if $\bigl( \eps , \{ \hat{\eps}_j' \}_{j=1}^{\hat{m}} \bigr) \in
\jointly_{\wsep,2+\hat{\ell}-\hat{k}}^{\hat{m} \sim \hat{k}}$ for some
$\hat{\ell} \in \scope{\hat{k}+1,\hat{m}}$ provided $\hat{k} \in \scope{\hat{m}-1}$; and
\begin{equation}
	\label{eq:fourthlochat}
 	\left\{
 		\begin{aligned}
			\!\! \d_{\hat{s}_{\hatellring-1}} \hat{u}_1(x,t,y,\hat{\vecs}_{\hat{m}})
				 - \grad_y \cdot \!\!\! \smint\limits_{\hat{S}^{\scope{\hatellring,\hat{m}}}}
					\! \hat{a}(x,t,y,\hat{\vecs}_{\hat{m}}; \grad u \! + \! \grad_y u_1)
						\, \rmd \hat{\vecs}_{\scope{\hatellring,\hat{m}}}				&= 0, \\
			\d_{\hat{s}_{\hat{\imath}}} \hat{u}_1(x,t,y,\hat{\vecs}_{\hat{m}})	&= 0
			\quad (\hat{\imath} \in \scope{\hatellring,\hat{m}})
 		\end{aligned}
 	\right.
\end{equation}
if $\bigl( \eps , \{ \hat{\eps}_j' \}_{j=1}^{\hat{m}} \bigr) \in
\jointly_{\wsep,1+\hat{m}+\hatellring-2\hat{k}}^{\hat{m} \sim \hat{k}}$ for some
$\hatellring \in \scope{\hat{k}+2,\hat{m}}$ provided $\hat{k} \in \scope{\hat{m}-2}$
and assuming $\hat{u}_1 \in L^2 \bigl( \Omega_T \times \hat{S}^{\hatellring-2} \times \hat{S}^{\scope{\hatellring,\hat{m}}};
H_{\#}^1(\hat{S}_{\hatellring-1};\wrum,\wrum') \bigr)$.
(For the sake of notational simplicity, we consider the strongly rather than weakly
formulated versions of the local problems.)

Define $u = \hat{u} \in H^1 \bigl( \ZeToT; H_0^1(\Omega), H^{-1}(\Omega) \bigr)$ which depends only on
$(x,t) \in \Omega_T$. \\
$\bullet$ We can write \eqref{eq:firstlochat} as
\begin{equation*}
	- \grad_{y} \cdot a(x,t,y,\vecs_m; \grad u + \grad_{y} \hat{u}_1)	= 0,
\end{equation*}
which is the local problem if $\bigl( \eps , \{ \hat{\eps}_j' \}_{j=1}^{\hat{m}} \bigr) \in
\jointly_{\wsep,1}^{\hat{m} \sim \hat{k}}$, i.e., $\bigl( \eps , \{ \hat{\eps}_j' \}_{j=1}^{\hat{m}} \bigr)
\in \jointly_{\wsep}^{\hat{m} \sim \hat{k}}$ and $\tfrac{\eps^2}{\hat{\eps}_{\hat{m}}'} \rightarrow 0$,
which is equivalent to $\bigl( \eps , \{ \eps_j' \}_{j=1}^m \bigr) \in \jointly_{\wsep}^{m \prec k}$ and
$\tfrac{\eps^2}{\eps_m'} \rightarrow 0$, i.e., we have precisely $\bigl( \eps , \{ \eps_j' \}_{j=1}^m \bigr) \in
\jointly_{\wsep,1}^{m \prec k}$. Obviously, $\hat{u}_1$ must be independent of $\hat{s}_{\hat{k}}$,
i.e., we can write the unique solution as $u_1 = \hat{u}_1 \in L^2 ( \Omega_T \times S^m; \wrum )$
which depends only on $(x,t) \in \Omega_T$ and $(y,\vecs_m) \in \sety_{1m}$. We thus conclude that the
local problem when $\bigl( \eps , \{ \eps_j' \}_{j=1}^m \bigr) \in \jointly_{\wsep,1}^{m \prec k}$ is
\begin{equation*}
	- \grad_{y} \cdot a(x,t,y,\vecs_m; \grad u + \grad_{y} u_1)	= 0,
\end{equation*}
and the homogenised flux $b \, : \, \bar{\Omega}_T \times \real^N \rightarrow \real^N$ is defined by
\begin{align}
	b(x,t;\grad u)	&= \hat{b}(x,t,\grad \hat{u}) \nonumber \\
						&= \smint_{\hat{\sety}_{1\hat{m}}} \hat{a}(x,t,y,\hat{\vecs}_{\hat{m}};
							\grad \hat{u} + \grad_{\hat{y}} \hat{u}_1)
								\, \rmd \hat{\vecs}_{\hat{m}} \rmd y \nonumber \\
						&= \smint_{\sety_{1m}} a(x,t, y,\vecs_m; \grad u
							+ \grad_{y} u_1) \, \rmd \vecs_m \rmd y. \label{eq:homofluxnohat}
\end{align}
(This is because $a^{\eps}(\, \cdot \; ;\grad u_{\eps}) = \hat{a}^{\eps}(\, \cdot \; ;\grad \hat{u}_{\eps}) \rightharpoonup
\hat{b}(\, \cdot \; ;\grad \hat{u}) = b(\, \cdot \; ;\grad u)$ in $L^2(\Omega_T)$; see Definition~\ref{def:Hmp-conv}.) \\
$\bullet$ We can write \eqref{eq:seclochat} as
\begin{equation*}
	\d_{s_m} \hat{u}_1(x,t,y,\hat{\vecs}_{\hat{m}})
		- \grad_y \cdot a(x,t,y,\vecs_m; \grad u + \grad_y \hat{u}_1) = 0,
\end{equation*}
which is the local problem if $\bigl( \eps , \{ \hat{\eps}_j' \}_{j=1}^{\hat{m}} \bigr) \in
\jointly_{\wsep,2}^{\hat{m} \sim \hat{k}}$, i.e., $\bigl( \eps , \{ \hat{\eps}_j' \}_{j=1}^{\hat{m}} \bigr)
\in \jointly_{\wsep}^{\hat{m} \sim \hat{k}}$ and $\hat{\eps}_{\hat{m}}' \sim \eps^2$, which is equivalent
to $\bigl( \eps , \{ \eps_j' \}_{j=1}^m \bigr) \in \jointly_{\wsep}^{m \prec k}$ and $\eps_m' \sim \eps^2$,
i.e., we have precisely $\bigl( \eps , \{ \eps_j' \}_{j=1}^m \bigr) \in \jointly_{\wsep,2}^{m \prec k}$.
Obviously, $\hat{u}_1$ must be independent of $\hat{s}_{\hat{k}}$, i.e., we can write the unique solution as $u_1 =
\hat{u}_1 \in L^2 ( \Omega_T \times S^m; \wrum )$ which depends only on $(x,t) \in \Omega_T$ and $(y,\vecs_m) \in
\sety_{1m}$. The assumption $\hat{u}_1 \in L^2 \bigl( \Omega_T \times \hat{S}^{\hat{m}-1};
H_{\#}^1(\hat{S}_{\hat{m}};\wrum,\wrum') \bigr)$ is clearly equivalent to $u_1 \in L^2 \bigl( \Omega_T \times
S^{m-1}; H_{\#}^1(S_m;\wrum,\wrum') \bigr)$. We thus conclude that the local problem when $\bigl( \eps ,
\{ \eps_j' \}_{j=1}^m \bigr) \in \jointly_{\wsep,2}^{m \prec k}$ assuming $u_1 \in L^2 \bigl( \Omega_T
\times S^{m-1}; H_{\#}^1(S_m;\wrum,\wrum') \bigr)$ is
\begin{equation*}
	\d_{s_m} u_1(x,t,y,\vecs_m) - \grad_y
		\cdot a(x,t,y,\vecs_m; \grad u + \grad_y u_1) = 0,
\end{equation*}
and the homogenised flux $b$ is given by \eqref{eq:homofluxnohat} again. \\
$\bullet$ Let $\bar{\ell}$ and $i$ be defined through $\hat{\ell} = \bar{\ell}+1$ and $\hat{\imath} = i+1$,
respectively; we can then write \eqref{eq:thirdlochat} as
\begin{equation*}
	\left\{
		\begin{aligned}
			- \grad_y \cdot \smint_{S^{\scope{\bar{\ell},m}}}
				a(x,t,y,\vecs_m; \grad u + \grad_y \hat{u}_1) \, \rmd \vecs_{\scope{\bar{\ell},m}}	&= 0, \\
			\d_{s_i} \hat{u}_1(x,t,y,\hat{\vecs}_{\hat{m}})											&= 0
				\qquad (i \in \scope{\bar{\ell},m}),
		\end{aligned}
	\right.
\end{equation*}
which are the local problems if $\bigl( \eps , \{ \hat{\eps}_j' \}_{j=1}^{\hat{m}} \bigr) \in
\jointly_{\wsep,2+\hat{\ell}-\hat{k}}^{\hat{m} \sim \hat{k}}$, $\hat{\ell} \in \scope{\hat{k}+1,\hat{m}}$, i.e., $\bigl( \eps ,
\{ \hat{\eps}_j' \}_{j=1}^{\hat{m}} \bigr) \in \jointly_{\wsep}^{\hat{m} \sim \hat{k}}$
and $\tfrac{\hat{\eps}_{\hat{\ell}}'}{\eps^2} \rightarrow 0$ but, only necessary if and only if $\hat{\ell} \neq \hat{k}+1 \,
\Leftrightarrow \, (\hat{k},\hat{\ell}) \neq (\hat{k},\hat{k}+1)$, $\tfrac{\hat{\eps}_{\hat{\ell}-1}'}{\eps^2}
\rightarrow \infty$. This is in turn equivalent to $\bigl( \eps , \{ \eps_j' \}_{j=1}^m \bigr) \in
\jointly_{\wsep}^{m \prec k}$ and $\tfrac{\eps_{\bar{\ell}}'}{\eps^2} \rightarrow 0$ but, if and only if $(k,\bar{\ell})
\neq (0,1) \, \Leftrightarrow \, (\hat{k},\hat{\ell}) \neq (1,2)$,
$\tfrac{\eps_{\bar{\ell}-1}'}{\eps^2} \rightarrow \infty$, i.e., we have precisely $\bigl( \eps ,
\{ \eps_j' \}_{j=1}^m \bigr) \in \jointly_{\wsep,2+\bar{\ell}-k}^{m \prec k}$, $\bar{\ell} \in
\scope{k+1,m}$. Obviously, $\hat{u}_1$ must be independent of $\hat{s}_{\hat{k}}$, i.e., we can write the
unique solution as $u_1 = \hat{u}_1 \in L^2 ( \Omega_T \times S^m; \wrum )$ which depends only on
$(x,t) \in \Omega_T$ and $(y,\vecs_m) \in \sety_{1m}$. We thus conclude that the local problems when
$\bigl( \eps , \{ \eps_j' \}_{j=1}^m \bigr) \in \jointly_{\wsep,2+\bar{\ell}-k}^{m \prec k}$
for some $\bar{\ell} \in \scope{k+1,m}$ are
\begin{equation*}
	\left\{
		\begin{aligned}
			- \grad_y \cdot \smint_{S^{\scope{\bar{\ell},m}}}
				a(x,t,y,\vecs_m; \grad u + \grad_y u_1) \, \rmd \vecs_{\scope{\bar{\ell},m}}	&= 0, \\
			\d_{s_i} u_1(x,t,y,\vecs_m)														&= 0
				\qquad (i \in \scope{\bar{\ell},m}),
		\end{aligned}
	\right.
\end{equation*}
and the homogenised flux $b$ is given by \eqref{eq:homofluxnohat} again. \\
Let $\ellring$ and $i$ be defined through $\hatellring = \ellring+1$ and $\hat{\imath} = i+1$,
respectively; we can then write \eqref{eq:fourthlochat} as
\begin{equation*}
 	\left\{
 		\begin{aligned}
			\d_{s_{\ellring-1}} \hat{u}_1(x,t,y,\hat{\vecs}_{\hat{m}})
				 - \grad_y \cdot \smint_{S^{\scope{\ellring,m}}}
					 a(x,t,y,\hat{\vecs}_{\hat{m}}; \grad u + \grad_y u_1)
						\, \rmd \vecs_{\scope{\ellring,m}}				&= 0, \\
			\d_{s_i} \hat{u}_1(x,t,y,\hat{\vecs}_{\hat{m}})	&= 0
			\qquad (i \in \scope{\ellring,m}),
 		\end{aligned}
 	\right.
\end{equation*}
which are the local problems if $\bigl( \eps , \{ \hat{\eps}_j' \}_{j=1}^{\hat{m}} \bigr) \in
\jointly_{\wsep,1+\hat{m}+\hatellring-2\hat{k}}^{\hat{m} \sim \hat{k}}$, $\hatellring \in \scope{\hat{k}+2,\hat{m}}$, i.e., $\bigl( \eps ,
\{ \hat{\eps}_j' \}_{j=1}^{\hat{m}} \bigr) \in \jointly_{\wsep}^{\hat{m} \sim \hat{k}}$
and $\tfrac{\hat{\eps}_{\hatellring-1}'}{\eps^2} \rightarrow 1$, which is equivalent to $\bigl( \eps , \{ \eps_j' \}_{j=1}^m \bigr) \in
\jointly_{\wsep}^{m \prec k}$ and $\tfrac{\eps_{\ellring-1}'}{\eps^2} \rightarrow 1$, i.e., we have precisely $\bigl( \eps ,
\{ \eps_j' \}_{j=1}^m \bigr) \in \jointly_{\wsep,1+m+\ellring-2k}^{m \prec k}$, $\ellring \in \scope{k+2,m}$.
Obviously, $\hat{u}_1$ must be independent of $\hat{s}_{\hat{k}}$, i.e., we can write the
unique solution as $u_1 = \hat{u}_1 \in L^2 ( \Omega_T \times S^m; \wrum )$ which depends only on
$(x,t) \in \Omega_T$ and $(y,\vecs_m) \in \sety_{1m}$.
The assumption  $\hat{u}_1 \in L^2 \bigl( \Omega_T \times \hat{S}^{\hatellring-2} \times \hat{S}^{\scope{\hatellring,\hat{m}}};
H_{\#}^1(\hat{S}_{\hatellring-1};\wrum,\wrum') \bigr)$ is obviously equivalent to $u_1 \in
L^2 \bigl( \Omega_T \times S^{\ellring-2} \times S^{\scope{\ellring,m}};H_{\#}^1(S_{\ellring-1};\mathcal{W},\mathcal{W}') \bigr)$.
We thus conclude that the local problems when
$\bigl( \eps , \{ \eps_j' \}_{j=1}^m \bigr) \in \jointly_{\wsep,1+m+\ellring-2k}^{m \prec k}$
for some $\ellring \in \scope{k+2,m}$ assuming $u_1 \in
L^2 \bigl( \Omega_T \times S^{\ellring-2} \times S^{\scope{\ellring,m}};H_{\#}^1(S_{\ellring-1};\mathcal{W},\mathcal{W}') \bigr)$ are
\begin{equation*}
 	\left\{
 		\begin{aligned}
			\d_{s_{\ellring-1}} u_1(x,t,y,\vecs_m)
				 - \grad_y \cdot \smint_{S^{\scope{\ellring,m}}}
					 a(x,t,y,\vecs_m; \grad u + \grad_y u_1)
						\, \rmd \vecs_{\scope{\ellring,m}}	&= 0, \\
			\d_{s_i} u_1(x,t,y,\vecs_m)						&= 0
			\qquad (i \in \scope{\ellring,m}),
 		\end{aligned}
 	\right.
\end{equation*}
and the homogenised flux $b$ is given by \eqref{eq:homofluxnohat} again.
The proof is complete.
\end{proof}
Define $\jointly_{\wsep}^{m \sim 0} = \emptyset$ and $\jointly_{\wsep,j}^{m \sim 0}
= \emptyset$, $j \in \scope{1+2m}$. Let $k \in \scope{m}_0$ and introduce $\jointly_{\wsep}^{m \preq k}
= \jointly_{\wsep}^{m \sim k} \cup \jointly_{\wsep}^{m \prec k}$ and
$\jointly_{\wsep,i}^{m \preq k} = \jointly_{\wsep,i}^{m \sim k} \cup
\jointly_{\wsep,i}^{m \prec k}$, $i \in \scope{1+2(m-k)}$. (The notation ``$\preq k$'' could be read as
``the spatial scale is asymptotically equal to or less than the $k$-th temporal scale''. The asymptotic equality
to the $0$-th temporal scale is meaningless which explains why we define the corresponding sets of pairs of lists
of scale functions as being empty.) From Theorems~\ref{th:maintheorem}~and~\ref{th:maintheorem2} we immediately
arrive in the corollary below, which is the main result of this e-print paper.
\begin{coro}
\label{cor:mainresult}
Let $k \in \scope{m}_0$. Suppose that the pair $e = \bigl( \eps , \{ \eps_j' \}_{j=1}^m \bigr)$ of lists of
spatial and temporal scale functions belongs to $\bigcup_{i=1}^{1+2(m-k)} \jointly_{\wsep,i}^{m \preq k}$.
Let $\{ u_{\eps} \}$ be the sequence of weak solutions in $H^1 \bigl( \ZeToT; H_0^1(\Omega),
H^{-1}(\Omega) \bigr)$ to the evolution problem \eqref{eq:monoparaprob} with $a \, : \, \bar{\Omega}_T
\times \real^{N+m} \times \real^N \rightarrow \real^N$ satisfying the structure conditions
\emph{(B}$_1$\emph{)}--\emph{(B}$_5$\emph{)}. Then
\begin{align*}
	u_{\eps}		&\rightarrow	u	&	&\textrm{in } L^2(\Omega_T), \\
	u_{\eps}		&\rightharpoonup  u	&	&\textrm{in } L^2 \bigl( \ZeToT; H_0^1(\Omega) \bigr), \\
\inter{and}
 	\grad u_{\eps}	&\scaleconv{(2,m+1)} \grad u + \grad_{y} u_1, & &
\end{align*}
where $u \in H^1 \bigl( \ZeToT; H_0^1(\Omega), H^{-1}(\Omega) \bigr)$ and
$u_1 \in L^2 ( \Omega_T \times S^m; \wrum )$. Here $u$ is the unique weak solution to the
homogenised problem \eqref{eq:defoflim} with the homogenised flux $b \, : \, \bar{\Omega}_T
\times \real^N \rightarrow \real^N$ given by
\begin{equation*}
	b(x,t;\grad u) = \smint_{\sety_{1m}} a(x,t, y,\vecs_m; \grad u + \grad_{y} u_1) \, \rmd \vecs_m \rmd y.
\end{equation*}
Moreover, we have the following characterisation of $u_1$:

$\bullet$ If $e \in \jointly_{\wsep,1}^{m \preq k}$ then the function $u_1$
is the unique weak solution to the local problem
\begin{equation*}
	- \grad_y \cdot a(x,t,y,\vecs_m; \grad u + \grad_y u_1) = 0.
\end{equation*}

$\bullet$ If $e \in \jointly_{\wsep,2}^{m \preq k}$, assuming $u_1 \in L^2 \bigl( \Omega_T \times S^{m-1};
H_{\#}^1(S_m;\wrum,\wrum') \bigr)$, then the function $u_1$ is the unique weak solution to the system
of local problems
\begin{equation*}
	\left\{
		\begin{aligned}
			\d_{s_m} u_1(x,t,y,\vecs_m) - \grad_y \cdot a(x,t,y,\vecs_m; \grad u + \grad_y u_1)	&= 0, \\
			- \grad_y \cdot \smint_{S_m} a(x,t,y,\vecs_m; \grad u + \grad_y u_1) \, \rmd \vecs_m 	&= 0.
		\end{aligned}
	\right.
\end{equation*}

$\bullet$ If $e \in \jointly_{\wsep,2+\bar{\ell}-k}^{m \preq k}$ for some $\bar{\ell} \in \scope{k+1,m}$,
provided $k \in \scope{m-1}_0$, then the function $u_1$ is the unique weak solution to the system of local problems
\begin{equation*}
	\left\{
		\begin{aligned}
			- \grad_y \cdot \smint_{S^{\scope{\bar{\ell},m}}}
				a(x,t,y,\vecs_m; \grad u + \grad_y u_1) \, \rmd \vecs_{\scope{\bar{\ell},m}}	&= 0, \\
				\d_{s_i} u_1(x,t,y,\vecs_m)	&= 0 \qquad (i \in \scope{\bar{\ell},m}).
		\end{aligned}
	\right.
\end{equation*}

$\bullet$ If $e \in \jointly_{\wsep,1+m+\ellring-2k}^{m \preq k}$ for
some $\ellring \in \scope{k+2,m}$, provided $k \in \scope{m-2}_0$ and assuming $u_1 \in
L^2 \bigl( \Omega_T \times S^{\ellring-2} \times S^{\scope{\ellring,m}};H_{\#}^1(S_{\ellring-1};\mathcal{W},\mathcal{W}') \bigr)$,
then the function $u_1$ is the unique weak solution to the system of local problems
\begin{equation*}
 	\left\{
 		\begin{aligned}
			\d_{s_{\ellring-1}} u_1(x,t,y,\vecs_m)
				- \grad_y \! \cdot \! \smint_{S^{\scope{\ellring,m}}}
					a(x,t,y,\vecs_m; \grad u + \grad_y u_1)
						\, \rmd \vecs_{\scope{\ellring,m}}	&= 0, \\
			\d_{s_i} u_1(x,t,y,\vecs_m)						&= 0
			\qquad (i \in \scope{\ellring,m}).
 		\end{aligned}
 	\right.
\end{equation*}
\end{coro}
\begin{rema}
Corollary~\ref{cor:mainresult} can only handle the subset $\bigcup_{i=1}^{1+2(m-k)}
\jointly_{\wsep,i}^{m \preccurlyeq k}$ of $\jointly_{\wsep}^{m \preccurlyeq k}$.
The conclusion of Proposition~\ref{prop:collispart} is true also in the setting of Corollary~\ref{cor:mainresult} though,
i.e., the collection $\bigl\{ \mathcal{P}_{i}^{m \preccurlyeq k} \bigr\}_{i=1}^{1+2(m-k)}$ forms a partition of
$\mathcal{P}^{m \preccurlyeq k}$ where $\mathcal{P}^{m \preccurlyeq k}$ is the subset of
$\jointly_{\wsep}^{m \preccurlyeq k}$ with temporal scale functions expressed as power functions,
and $\mathcal{P}_{i}^{m \preccurlyeq k}$ is the corresponding subset of
$\jointly_{\wsep,i}^{m \preccurlyeq k}$ for every $i \in \scope{1+2(m-k)}$.
\end{rema}


\begin{thebibliography}{45}
	\bibitem{Alla92}	Allaire, G.
						``Homogenization and two-scale convergence'',
						\emph{SIAM J.~Math.~Anal.} 23 (1992),  no.~6, 1482--1518.
	\bibitem{AllBri96}	Allaire, G.; Briane, M.
						``Multiscale convergence and reiterated homogenisation'',
						\emph{Proc.~Roy.~Soc.~Edinburgh Sect.~A} 126 (1996), no.~2, 297--342.
	\bibitem{ArDoHo90}	Arbogast, T.; Douglas, J., Jr.; Hornung, U.
						``Derivation of the double porosity model of single phase flow
							via homogenization theory'',
						\emph{SIAM J.~Math.~Anal.} 21 (1990), no.~4, 823--836.
	\bibitem{BeLiPa78}	Bensoussan, A.; Lions, J.-L.; Papanicolaou, G.
						``Asymptotic analysis for periodic structures'',
						Studies in Mathematics and its Applications, 5. \emph{North-Holland Publishing Co.,
							Amsterdam-New York}, 1978.
	\bibitem{ChDaDe90}	Chiad\`{o} Piat, V.; Dal Maso, G.; Defranceschi, A.
						``$G$-convergence of monotone operators'',
						\emph{Ann.~Inst.~H.~Poincar\'{e} Anal.~Non Lin\'{e}aire} 7 (1990), no.~3, 123--160.
	\bibitem{ChiDef90}	Chiad\`{o} Piat, V.; Defranceschi, A.
						``Homogenization of monotone operators'',
						\emph{Nonlinear Anal.} 14 (1990), no.~9, 717--732.
	\bibitem{CiDaGr02}	Cior\u{a}nescu, D.; Damlamian, A.; Griso, G.
						``Periodic unfolding and homogenization'',
						\emph{C.~R.~Math.~Acad.~Sci.~Paris} 335 (2002), no.~1, 99--104.
	\bibitem{CioDon99}	Cior\u{a}nescu, D.; Donato P.
						``An introduction to homogenization'',
						\emph{Oxford University Press}, New York, 1999.
	\bibitem{ColSpa77}	Colombini, F.; Spagnolo, S.
						``Sur la convergence de solutions d'\'{e}quations paraboliques'',
						\emph{J.~Math.~Pures Appl.~(9)} 56 (1977), no.~3, 263--305.
	\bibitem{ConFoi88}	Constantin, P.; Foia\c{s}, C.
						``Navier-Stokes equations'',
						\emph{University of Chigaco Press}, Chicago, 1988.
	\bibitem{Evan89}	Evans, L.~C.
						``The perturbed test function method for viscosity solutions of nonlinear PDE'',
						\emph{Proc.~Roy.~Soc.~Edinburgh Sect.~A} 111 (1989), no.~3-4, 359--375.
	\bibitem{Evan92}	Evans, L.~C.
						``Periodic homogenization of certain fully nonlinear
							partial differential equations'',
						\emph{Proc.~Roy.~Soc.~Edinburgh Sect.~A} 120 (1992), no.~3-4, 245--265.
	\bibitem{FloOls06}	Flod\'{e}n, L.; Olsson, M.
						``Reiterated homogenization of some linear and nonlinear
							monotone parabolic operators'',
						\emph{Can.~Appl.~Math.~Q.} 14 (2006), no.~2, 149--183.					
	\bibitem{FloOls07}	Flod\'{e}n, L.; Olsson, M.
						``Homogenization of some parabolic operators with several time scales'', 
						\emph{Appl.~Math.} 52 (2007), no.~5, 431--446.
	\bibitem{FHOS07}	Flod\'{e}n, L.; Holmbom, A.; Olsson, M.; Svanstedt, N.
						``Reiterated homogenization of monotone parabolic problems'',
						\emph{Ann.~Univ.~Ferrara Sez.~VII Sci.~Mat.} 53 (2007), no.~2, 217--232.
	\bibitem{Holm96}	Holmbom, A.
						``Some modes of convergence and their application to
							homogenization and optimal composites design'',
						Doctoral thesis 1996:208~D, Department of mathematics, Lule{\aa} university, 1996.
	\bibitem{Holm97}	Holmbom, A.
						``Homogenization of parabolic equations: an alternative
							approach and some corrector-type results'',
						\emph{Appl.~Math.} 42 (1997), no.~5, 321--343.
	\bibitem{HolSil06}	Holmbom, A.; Silfver, J.
						``On the convergence of some sequences of oscillating functionals'',
						\emph{WSEAS Trans.~Math.} 5 (2006), no.~8, 951--956.
	\bibitem{HSSW06}	Holmbom, A.; Silfver, J.; Svanstedt, N.; Wellander, N.
						``On two-scale convergence and related sequential compactness topics'',
						\emph{Appl.~Math.} 51 (2006), no.~3, 247--262.
	\bibitem{HoSvWe05}	Holmbom, A.; Svanstedt, N.; Wellander, N.
						``Multiscale convergence and reiterated homogenization of parabolic problems'',
						\emph{Appl.~Math.} 50 (2005),  no.~2, 131--151.
	\bibitem{Kunch88}	Kun'ch, R.~N.
						``$G$-convergence of nonlinear parabolic operators'',
						Thesis, Donetsk, 1988.
	\bibitem{KunPan86}	Kun'ch, R.~N.; Pankov, A.~A.
						``$G$-convergence of monotone parabolic operators'',
						\emph{Dokl.~Akad.~Nauk Ukrain.~SSR Ser.~A} 1986,  no.~8, 8--10.
	\bibitem{LLPW01}	Lions, J.-L.; Lukkassen, D.; Persson, L.-E.; Wall, P.
						``Reiterated homogenization of nonlinear monotone operators'',
						\emph{Chinese Ann.~Math.~Ser.~B} 22  (2001),  no.~1, 1--12.
	\bibitem{LuNgWa02}	Lukkassen, D.; Nguetseng, G.; Wall, P.
						``Two-scale convergence'',
						\emph{Int.~J.~Pure Appl.~Math.} 2 (2002), no.~1, 35--86.
	\bibitem{MasToa01}	Mascarenhas, M.~L.; Toader, A.-M.
						``Scale convergence in homogenization'',
						\emph{Numer.~Funct.~Anal.~Optim.} 22 (2001), no.~1-2, 127--158.
	\bibitem{Mura78}	Murat, F. ``H-convergence'',
						S\'{e}minaire d'analyse fonctionelle et num\'{e}rique
							de l'Universit\'{e} d'Alger, 1978.
	\bibitem{Mura78b}	Murat, F.
						``Compacit\'{e} par compensation'',
						\emph{Ann.~Scuola Norm.~Sup.~Pisa Cl.~Sci.~(4)} 5 (1978), no.~3, 489--507.
	\bibitem{Mura79}	Murat, F.
						``Compacit\'{e} par compensation.~II'',
						\emph{Proceedings of the International Meeting on Recent Methods in
							Nonlinear Analysis (Rome, 1978)}, pp.~245--256, Pitagora, Bologna, 1979.
	\bibitem{Nech04}	Nechv\'{a}tal, L.
						``Alternative approaches to the two-scale convergence'',
						\emph{Appl.~Math.} 49 (2004), no.~2, 97--110.
	\bibitem{Ngue89}	Nguetseng, G.
						``A general convergence result for a functional
							related to the theory of homogenization'',
						\emph{SIAM J.~Math.~Anal.} 20 (1989),  no.~3, 608--623.
	\bibitem{Ngue03}	Nguetseng, G.
						``Homogenization structures and applications.~I'',
						\emph{Z.~Anal.~Anwendungen} 22 (2003),  no.~1, 73--107.
	\bibitem{Ngue04}	Nguetseng, G.
						``Homogenization structures and applications.~II'',
						\emph{Z.~Anal.~Anwendungen} 23 (2004),  no.~3, 483--508.
	\bibitem{NguWou04}	Nguetseng, G.; Woukeng, J.~L.
						``Deterministic homogenization of parabolic monotone
							operators with time dependent coefficients'', 
						\emph{Electron.~J.~Differential Equations} 2004, No.~82, 23~pp.
	\bibitem{NguWou07}	Nguetseng, G.; Woukeng, J.~L.
						``$\Sigma$-convergence of nonlinear parabolic operators'',
						\emph{Nonlinear Anal.} 66 (2007), no.~4, 968--1004.
	\bibitem{Pank97}	Pankov, A.
						``$G$-convergence and homogenization of
							nonlinear partial differential operators'',
						Mathematics and its Applications, 422.~\emph{Kluwer
							Academic Publishers}, Dordrecht, 1997.
	\bibitem{Silf07}	Silfver, J.
						``$G$-convergence and homogenization involving
							operators compatible with two-scale convergence'',
						Mid Sweden University Doctoral Thesis 23, 2007.
	\bibitem{Spag67}	Spagnolo, S.
						``Sul limite delle soluzioni di problemi
							di Cauchy relativi all'equazione del calore'',
						\emph{Ann.~Scuola Norm.~Sup.~Pisa (3)} 21 (1967), 657--699.
	\bibitem{Spag68}	Spagnolo, S.
						``Sulla convergenza di soluzioni di equazioni paraboliche ed ellittiche'',
						\emph{Ann.~Scuola Norm.~Sup.~Pisa (3)} 22 (1968), 571--597.
	\bibitem{Spag76}	Spagnolo, S.
						``Convergence in energy for elliptic operators'',
						\emph{Numerical solution of partial differential equations, III
							(Proc.~Third Sympos.~(SYNSPADE), Univ.~Maryland, College Park, Md., 1975)},
								pp.~469--498. \emph{Academic Press, New York}, 1976.
	\bibitem{Svan92}	Svanstedt, N.
						``G-convergence and homogenization of sequences
							of linear and nonlinear partial differential operators'',
						Doctoral thesis 1992:105~D, Department of mathematics, Lule{\aa} university, 1992.
	\bibitem{Svan99}	Svanstedt, N.
						``$G$-convergence of parabolic operators'',
						\emph{Nonlinear Anal.} 36 (1999), no.~7, Ser.~A: Theory Methods, 807--842.
	\bibitem{Tart77}	Tartar, L.
						``Cours peccot'', Coll\`{e}ge de France, 1977,
							unpublished, partially written in \cite{Mura78}.
	\bibitem{Tart78}	Tartar, L.
						``Quelques remarques sur l'homog\'{e}n\'{e}isation'',
						In: M.~Fujita, Editor, \emph{Functional Analysis and Numerical Analysis,
							Proc.~Japan–France Seminar 1976},
								Japanese Society for the Promotion of Science (1978), 468-–482.
	\bibitem{Tart79}	Tartar, L.
						``Compensated compactness and applications to partial differential equations'',
						\emph{Nonlinear analysis and mechanics: Heriot-Watt Symposium, Vol.~IV},
							136--212, Res.~Notes in Math., 39, \emph{Pitman, Boston, Mass.-London}, 1979.
	\bibitem{Weic98}	Weickert, J.
						``Anisotropic diffusion in image processing'',
						\emph{Teubner Verlag}, Stuttgart, 1998. Available online from: \\
							{\tt \href{http://www.mia.uni-saarland.de/weickert/book.html}{
							http://www.mia.uni-saarland.de/weickert/book.html}}.
	\bibitem{Wouk10}	Woukeng, J.~L.
						``Periodic homogenization of nonlinear non-monotone
							parabolic operators with three time scales'',
						\emph{Ann.~Mat.~Pura Appl.~(4)} 189 (2010),  no.~3, 357--379.
	\bibitem{ZeidIIA}	Zeidler, E.
						``Nonlinear functional analysis and its applications. II/A.
							Linear monotone operators'',
						\emph{Springer-Verlag}, New York, 1990.
	\bibitem{ZeidIIB}	Zeidler, E.
						``Nonlinear functional analysis and its applications. II/B.
							Nonlinear monotone operators'',
						\emph{Springer-Verlag}, New York, 1990.
\end{thebibliography}
\end{document}